\input amstex

\documentstyle{amsppt}
\magnification=1100
\loadbold
\input epsf
\mathsurround=1.2pt

\overfullrule=0pt
%
\def\op {\operatorname}
\define\GR#1#2{{\text{\bf #1}}_{#2}}
\define\GRt#1{ {\tilde {\text{\bf A}}}_{#1}}
\def\V{\Bbb V}
\def\co{\op{O}}
\def\el{\frak l}
\def\cc{{\Cal C}}

\def\sel#1{{\frak s\frak l}_{#1}}
\def\tri{{\frak s\frak l}_2}
\def\ka{\frak k}
\def\dis{\displaystyle  }
\def\eps{{\varepsilon}}
\def\iso{\,{\overset {\sim} \to\longrightarrow}\,}
\def\ZZ{\Cal Z}
\def\bs{{\bold S}}
\def\C{\Bbb C}
\def\cd{\!\cdot\!}
\def\B{{\Cal B}}

\def\F{{\Cal F}}

\def\tlim{{\text{\it lim}}}
\def\into{{\;\hookrightarrow\;}}

\def\nw{_{\!_{{\bold {nw}}}}}
\def\sw{_{\!_{{\bold {sw}}}}}
\def\ne{_{{\bold {ne}}}}
\def\north{^{\!^{{\bold {north}}}}}
\def\east{^{\!^{{\bold {east}}}}}
\def\gnw{\pmb{\g}_{_{{\bold {nw}}}}}
\def\gnww{\g_{_{{\bold {nw}}}}^}
\def\nnw{\n_{_{{\bold {nw}}}}}
\def\znw{\z_{_{{\bold {nw}}}}}
\def\nsw{\n_{_{{\bold {sw}}}}}


\def\PP{{\Cal P}}
\def\O{\bold O}
\def\0{{\bold 0}}

\def\P{{\boldsymbol\Delta}}
\def\pe{{\boldsymbol\Delta}_{\bold e}}

\def\x{\bold x}
\def\y{\bold y}

\def\m{{\frak m}}
\def\l{{\frak l}}

\def\e{\bold e}

\def\bgm{\pmb{\frak g}_{_{\pmb{\bold{-}}}}}

\def\Rsmall{{R_{_\clubsuit}}}

\def\pp{{\frak P}}

\def\uu{{\bold u}}
\def\vv{{\bold v}}
\def\ll2{{\frak s\frak l}_2}

\def\bnu{{\boldsymbol\nu}}
\def\blambda{{\boldsymbol\lambda}}

\def\Z{\Bbb Z}

\def\z{\frak z}
\def\h{\frak h}
\def\ZZ{\Cal Z}

\def\sln{{\frak s\frak l}_n}
\def\sll{{\frak s\frak l}}
\def\gl{{\frak g\frak l}}
\def\so{{\frak s\frak o}}

\def\g{{\frak g}}
\def\u{\frak V}
\def\s{{\frak s}}
\def\c{\frak c}
\def\bc{\pmb{{\frak c}^\circ}}
\def\p{{\frak p}}

\def\n{{\frak n}}
\def\u{{\frak u}}

\def\d{{\bold d}}
\define \SS{\operatorname{S}}
\define \Ad {{\operatorname{Ad}^{}}}
\define \Ker {{\operatorname{Ker}^{}}}
\define \ad {{\operatorname{ad}^{}}}
\define \Lie {{\operatorname{Lie}^{}}}
\def\bh{\bold h}

\define \Taut {\text{\it Taut}}
\define \Par {\text{\it Par}}
\define \GL {\operatorname{GL}}  \define \SL {\operatorname{SL}}
\define \Hilb {\operatorname{Hilb}}  \define \Proj {\operatorname{Proj}}
\define \rk {\operatorname{rk}}
\define \reg {{\text{\rm reg}}}  
\define \sign {\pmb{\boldsymbol\varepsilon}}
\define \triv {\pmb{\bold1}}
\define \Specm {\operatorname{Specm}}  
\define \Hom {\operatorname{Hom}}
\define \Hol {\operatorname{\Cal Hol}}
\define \Sol {\operatorname{\Cal Sol}}
\define \Alt {\operatorname{\Cal Alt}}
\define \Gr {\operatorname{Gr}}
\define \ind {\operatorname{Ind}}
\define \gr {{\operatorname{gr}^{}}}
\define \bll {\bold {\bigl(\!\bigl(}} 
\define \rr  {\bold {\bigr)\!\bigr)}}
\define \quot {\operatorname{\Cal Q\text{\it uot}}}
\define \Ind {\operatorname{Ind}}

\TagsOnRight
\nopagenumbers


\topmatter

\title
Principal Nilpotent pairs in a semisimple Lie algebra $\,\pmb{1}$.
\endtitle
\author
VICTOR GINZBURG
\endauthor\footnote"${}$"{Address: The University of Chicago, Mathematics Department,
Chicago, IL 60637, USA; \newline
$\hphantom{x}\quad$E-mail: {\bf {ginzburg\@math.uchicago.edu}}}

\abstract
This is the first of a series of papers devoted to
certain pairs of commuting  nilpotent elements in a
semisimple
Lie algebra that enjoy quite remarkable properties and which are
expected to play a major role in Representation theory.
The properties of these pairs and their role is similar
to those of the principal nilpotents. 
To any principal nilpotent pair  we
associate a two-parameter analogue of the Kostant partition function,
and propose the corresponding
 two-parameter analogue of the weight multiplicity formula.

In a different direction, each principal
nilpotent pair gives rise to a harmonic
polynomial on the {\it Cartesian square}
of the Cartan subalgebra, that transforms under an irreducible
representation of the Weyl group. In the special case of $\sln$,
 the conjugacy classes of principal
nilpotent pairs and the irreducible representations of the Symmetric
group,
$S_n$,
are both parametrised (in a compatible way) by Young diagrams. In general,
 our theory 
provides a natural  generalization to  arbitrary Weyl groups
of the classical  construction of simple $S_n$-modules in terms of
{\it Young's symmetrisers}.

First results towards a complete classification of all principal
nilpotent pairs in a simple Lie algebra are presented at the end
of this paper in an
Appendix, written by A. Elashvili and D. Panyushev.
\endabstract

%

\endtopmatter
\document

\medskip

\centerline{\bf Table of Contents}
\vskip 5mm

\qquad${}^{\text{\bf  1.{ $\;$} Principal 
 nilpotent pairs}}_{\text{\bf  
2.{$\;\;\;$} Bi-filtration associated to a nilpotent pair}}$

\qquad${}^{\text{\bf  3.{ $\;$} Two parabolics}}_{\text{\bf 
4.{$\;\;\;$} Harmonic polynomial 
attached to a principal nilpotent pair}}$
\vskip 1pt
\qquad${}^{\text{\bf  5.{ $\;_{}$} Distinguished nilpotent pairs;
$\pmb{\frak s}\pmb{\frak l}_{\bold n}$-case}}_{\text{\bf 
6.{$\;\;\;$}  Some 
cohomology and generating functions}}$
\vskip 1pt

\qquad${}^{\text{\bf  7.{ ${}_{_{}}\;$} Partial slices}}_{\text{\bf  
8.{ $\;{}_{}$} Appendix by A. Elashvili and 
D. Panyushev:}}$
\vskip 1pt
\qquad${}^{\;\qquad\text{\bf "Towards a classification
of principal nilpotent pairs"}}$

\bigskip
\vskip 1pt

\head{0. Introduction.} \endhead
\medskip

Let $\g$ be a complex  semisimple Lie algebra,
and $G$ the corresponding adjoint group, i.e.,
the identity component
of $\text{\it Aut}(\g)\,$. 

Recall that an element
$x\in \g$ is called {\it regular} if $\z_\g(x)$, the centralizer of $x$
in $\g$, has the minimal possible dimension,
i.e., $\dim \z_\g(x) = \rk\g$. 
The most interesting,
in a sense, among  regular elements of $\g$ are
regular {\it nilpotent} elements, called "principal nilpotents".
These elements form a single $\Ad G$-orbit
in $\g$.
We refer to the  beautiful papers [K1], [K2]
for a comprehensive  study of principal nilpotents.

We 
are going to "double" the above setup and
replace the Lie algebra $\g$ by $\,\ZZ\subset \g\oplus\g\,,$ 
the set of
all pairs $(x_1,x_2)\in \g\oplus\g$, such that $[x_1,x_2]=0$,
called 
the {\it commuting
variety} of $\g$, 
see [R]. 
It turns out that there
is a remarkable  "doubled" counterpart
of the notion  of a principal nilpotent,
that we call a {\it principal nilpotent pair}.
The underlying idea is best explained as follows.
Let $\ZZ^{\text{reg}}$
be the union of all $\Ad G$-diagonal orbits of maximal dimension in
$\ZZ$. 
This is a smooth Zariski open, dense subset in $\ZZ$.
The natural $\C^*$-action on $\g$ by dilations gives rise
to a $\C^*\times \C^*$-action on $\g\oplus\g$, and
$\ZZ^{\text{reg}}$ is a $\C^*\times \C^*$-stable subvariety.
A principal nilpotent pair is the one whose 
$\Ad G$-diagonal orbit is a  fixed point of
the induced $\C^*\times \C^*$-action on $"\ZZ^{\text{reg}}/\Ad G^{} "$.
Here $"\ZZ^{\text{reg}}/\Ad G^{} "$ denotes the naive set of orbits;
it has no structure of an algebraic variety.
This construction was motivated in part
by our work on Hilbert schemes (joint work
with R. Bezrukavnikov currently in progress).

The set of principal nilpotent pairs does  {\it not}
form a single orbit under $\Ad G$-diagonal action on $\ZZ$,
but it consists of only finitely many such orbits. 
In the case
 $\g=\sln$, for instance,  conjugacy classes of 
principal nilpotent pairs are parametrized essentially
(up to transposing matrices)
by  Young diagrams with $n$-boxes. This comes about as follows.
Given a Young diagram:
$$ \tag 0.1 $$
\hskip30mm \epsffile{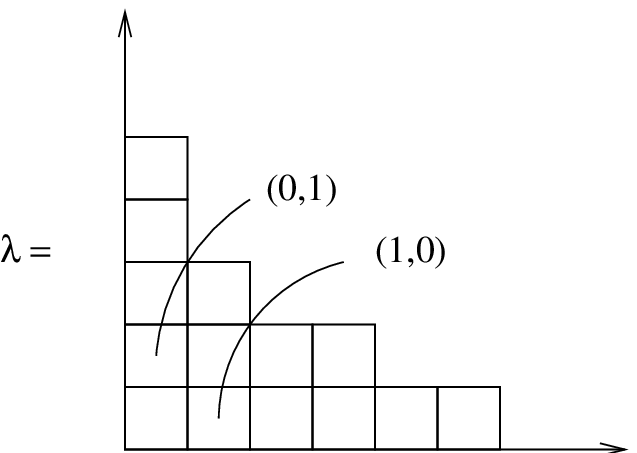}
$$   $$
we enumerate its boxes in some order, and
label the standard base vectors in $\C^n$ by
the box with the corresponding number. We  define
an endomorphism  $e_1\in\sln$  by letting
it act  "along the rows of the diagram",
i.e., by sending the base vector labelled by a box
to the base vector labelled by the
next right box, if this box
 belongs to $\lambda$, and to $0$ otherwise.
Similarly, we define $e_2\in\sln$ by letting
it act "along the columns",
from bottom to top. It is easy to see that the operators
$e_1, e_2$ thus defined commute, and form
a principal nilpotent pair. Note that if $\lambda$ consists
of either a single row or a single column, then the corresponding
principal nilpotent pair is either of the form $(e, 0)$ or $(0,e)$,
where $e$ is a  principal nilpotent in $\sln$ in the ordinary sense.
Moreover, any principal nilpotent pair
for $\g=\sln$ is associated to a Young diagram in the above way
(up to conjugation and transposing matrices).

For Lie algebras of
 types other than $\sln$ the number of principal nilpotent pairs
is typically less than one might have expected, because the
the set   $"\ZZ^{\text{reg}}/\Ad G\,"$ is even
farther away from being
an algebraic variety (in the $\g=\sln$-case the corresponding algebraic
variety  exists, in a sense;
it is the Hilbert scheme of $n$ points on $\C^2$).

In Section 2, to each principal nilpotent pair in
an arbitrary semisimple Lie algebra $\g$, and each simple finite-dimensional
$\g$-module, we associate a certain two-variable
analogue
of  Kostant's partition function,
and propose the corresponding $(s,t)$-weight multiplicity formula.
This formula may be related to
double analogues of Kostka polynomials that arise
in Macdonald theory [M2]. The geometry underlying our
multiplicity formula will be discussed in the second paper of this series.
Later on, we are also going
to express the two-variable weight multiplicity function
in terms of Intersection cohomology of an appropriate double-loop
 Grassmannian,
similar to the way, the $q$-analogue of weight multiplicity introduced
by Lusztig  is related to perverse sheaves on the loop Grassmannian,
see [L1], [Gi].

In section 4 we associate to a 
principal nilpotent pair in $\g$
a harmonic
polynomial on the Cartesian square of a Cartan subalgebra
of $\g$. This polynomial has very interesting properties,
in particular transforms under an irreducible 
representation of the  Weyl group.
In the special case $\g=\sln$ the polynomial in question
has been known in Combinatorics as a double-analogue 
of the Vandermonde determinant.
In the latter case, 
the irreducible representation of the Symmetric group
generated by the 
 polynomial turns out to be
parametrized by the
Young diagram labelling the  principal nilpotent pair.
Our construction is closely related to 
the theory of Springer representations.

In section 6 we compute the tangent space to
$"\ZZ^{\text{reg}}/\Ad G^{} "$ at a principal nilpotent conjugacy
class and use this computation to derive some
combinatorial identities involving a double-analogue
of the notion of exponents of a semisimple Lie algebra.
 In the last section we
introduce a distinguished "partial transverse slice"
to each principal nilpotent conjugacy
class, which is a sort of  generalisation 
of Kostant-Slodowy slices to nilpotent orbits in $\g$.
The parallelism with Kostant-Slodowy theory goes 
surprisingly far, see notably Theorem 7.4.
This parallelism suggests the existence of a yet unknown
double-analogue of Whittaker modules theory, similar to the
one developed in [K3]. We hope to return to this issue
elsewhere.
\bigskip

\noindent
{\it Acknowledgements.} I am grateful to A. Broer for helping me
with the proof of Proposition 4.11, to M. Finkelberg for helping me
with Theorem 7.4,
and to A. Joseph and V. Ostrik
for very useful suggestions 
at the end of \S4.
\medskip
\bigskip

\head{1. Principal nilpotent pairs.} \endhead
\bigskip

Given two elements $x_1,x_2$ in a Lie algebra ${\frak a}$, we 
write $\x=(x_1,x_2)$, and let $\z_{{\frak a}}(x_i)$, resp.
$\z_{{\frak a}}(\x)$, denote the centralizer, resp. the
the simultaneous
centralizer of $x_1$ and $x_2$,  in ${\frak a}$.

Fix a complex semisimple Lie algebra $\g$ with adjoint group $G$.
We will  write 
$\z(\x)$ instead of $\z_\g(\x)$, for simplicity.
Let 
$\ZZ\subset \g\oplus\g$
be the commuting variety of $\g$,
as defined in the Introduction. By a theorem of Richardson [R],
the pairs $\x=(x_1,x_2)$, where $x_1,x_2$ are in the same  Cartan
subalgebra
of $\g$, form a dense subset in $\ZZ$. In particular, for any 
$\x\in\ZZ$ one has: $\dim\z(\x) \ge \rk\g$. 
 We say that 
$\x=(x_1,x_2)\in \ZZ$ is a {\it regular pair} if the equality
$\dim\z(\x)=\rk\g$ holds.
\medskip

\noindent
{\bf Definition 1.1.} A pair $\e=(e_1,e_2) \in {\frak g}\oplus {\frak g}$
is called a {\it principal nilpotent pair} if the following conditions
$\op{(Reg)}$ and $\op{(Nil)}$ hold:
$$\align
&{\bold {(Reg)}} \quad\e=(e_1,e_2) \in\ZZ \text{ is a regular pair, i.e.},
[e_1, e_2]=0 \text{ and }\quad \dim \z(\e)=\rk {\frak g};\\
&{\bold {(Nil)}}\quad
\text{For any }(t_1,t_2)\in \C^*\times \C^*, \text{ there exists}\quad
g=g(t_1,t_2)\in G\quad\text{
 such that  }\\
&\quad\quad\quad\quad\quad
(\,t_1\!\cdot\! e_1,\,t_2\!\cdot\! e_2\,) = (\,\operatorname{Ad} g(
e_1),\, \operatorname{Ad} g (e_2)\,).
\endalign
$$
\medskip

Note that any elements $e_1,e_2$ satisfying the (Nil)-condition  are
necessarily nilpotent, since $\Ad G$-conjugacy classes of $e_1$ and $e_2$
are both stable under dilation. 

As a first example, take $e\in\g$ to be a regular (=principal)
nilpotent
in $\g$ (in the ordinary sense). 
Then $(e,0)$ and $(0,e)$ are two principal nilpotent 
pairs in $\g$. In general, if $(e_1,e_2)$ is a principal nilpotent 
pair then so is $(e_2,e_1)$.

Write $Z_G^{\bold {unip}}(\e)$ for the  unipotent radical
of the centralizer of the pair
$\e$ in $G$.

\proclaim{Theorem 1.2} Given a principal nilpotent 
pair $\e=(e_1,e_2)$, there exists a pair $\bh=(h_1, h_2)\in \ZZ$
 formed by semisimple elements of $\g$
and such that one has:

\vskip 2pt
$\op{(i)}\qquad\qquad\dis [h_i, e_j]= \delta_{i,j}\cdot e_i\quad,\;
i,j\in \{1,2\}\quad,\quad[h_1,h_2]=0.$
\vskip 2pt

$\op{(ii)}$ The pair $\bh=(h_1, h_2)$ is regular, i.e.,
$\z(\bh)$ is a Cartan subalgebra in $\g$.

$\op{(iii)}$ All eigenvalues of the operators
$\ad h_i : \g\to\g$ ($i=1,2$) are integral.

\noindent
Furthermore, the semi-simple pair $\bh$ is 
determined by the nilpotent pair $\e$
uniquely up to conjugacy, more precisely:

$\op{(iv)}$ For any two semisimple pairs $\bh=(h_1,h_2)$ and $\bh'=(h'_1 ,h'_2)$,
each satisfying commutation relations $\op{1.2(i)}$, there exists
$u\in Z_G^{\bold {unip}}(\e)$ such that $\bh'= \Ad u (\bh)$.
\endproclaim

\noindent
A pair $\bh=(h_1, h_2)$
of  semisimple elements of $\g$ satisfying commutation relations
1.2(i)
will be referred to as  an {\it associated
semisimple} pair.
\smallskip

Proof of the theorem will be done in stages and will occupy most of
this
section.
We  first derive a few consequences
of Definition 1.1. 

\proclaim{Lemma 1.3} For any 
commuting pair $\e=(e_1, e_2)$ 
 satisfying the (Nil)-condition in $\op{(1.1)}$, 
 there exist integers $m_1,\,m_2 >0$, and
 an algebraic group homomorphism
$\gamma: \C^*\times\C^* \to G$ such that
$$\bigl(\operatorname{Ad} \gamma(t_1, t_2)\bigr)\,e_\imath
= t_\imath^{m_\imath}\!\cdot\!e_\imath\quad,
\enspace\forall (t_1, t_2)\in\C^*\times\C^* \,,\,
\imath=1,2.
$$
\endproclaim
\smallskip

We will see later in Corollary 3.6 that in the case of a principal
nilpotent pair $\e=(e_1, e_2)$
the integers $m_1,\,m_2$ can be in effect chosen to be equal to $1$.
\smallskip

{\it Proof of Lemma.} Let  
$M\subset \C^*\times\C^*\times G$ be the algebraic group
formed by all triples:
$$ M=\{(t_1,t_2, g)\in \C^*\times\C^*\times G
\enspace | \enspace
\operatorname{Ad} g (e_1)= t_1\!\cdot\!e_1\;,\;
\operatorname{Ad} g (e_2)= t_2\!\cdot\!e_2 \}\;.$$
The assumptions of the lemma imply that the  projection
$\,p: M\to \C^*\times\C^*$ on the
first two factors
 is surjective. This gives a short exact sequence
$$
1\to Z_G(e_1, e_2)\to M
\overset {p} \to 
\longrightarrow  \C^*\times\C^*\to 1\tag 1.4
$$
Let $M^\circ$ be 
the identity component of $M$, and
$M^\circ = M_{\bold {red}}\cdot M_{\bold {unip}}$  the
decomposition of $M^\circ$ as a semidirect product of a reductive
subgroup, $M_{\bold {red}}$, and the unipotent radical.
The restriction of the homomorphism
$p$ to $M^\circ$ remains surjective,
and clearly vanishes both on $M_{\bold {unip}}$ and on the semisimple
part of $M_{\bold {red}}$. Writing 
$C$ for the connected center of  $M_{\bold {red}}$,
we see that the projection $p: M^\circ \twoheadrightarrow
 \C^*\times\C^*$ restricts to a surjective
homomorphism  $C\twoheadrightarrow \C^*\times\C^*$. Note that
$C$ is a torus, and any surjective homomorphism of tori
admits a quasi-splitting, i.e.,
there are integers ${m_1,\,m_2 >0}$ and an algebraic homorphism
 $s : \C^*\times\C^* \hookrightarrow C$ such that,
$\forall\,(t_1,t_2)\in
\C^*\times\C^*$, one has:
$\,p\circ s (t_1,t_2) = (t_1^{m_1},t_2^{m_2})\,.$
 The composition
$$\gamma :\, \C^*\times\C^* \overset {s} \to 
\longrightarrow\, C\, \hookrightarrow \, M\,\hookrightarrow\,
\C^*\times\C^*\times G\,
\;
 \overset {_{\operatorname{pr}_{_3}}} \to \longrightarrow \;G$$  
gives a homomorphism with the desired properties.
\qed

\medskip

{\it Proof of parts $\op{(i),(ii)}$ of Theorem 1.2.}
Using the homomorphism $\gamma$ of the Lemma,
we define two elements, $h_1, h_2\in\g$, by
$\;h_\imath := \frac{1}{m_\imath}\cdot\frac{\partial\gamma}{\partial t_\imath} 
\big|_{t_1=t_2=1}\;,\, \imath=1,2.$ 
These elements are semisimple, since the image of $\gamma$ is
contained
in a torus in $G$, and by construction
satisfy the commutation relations 1.2(i). This proves part (i) of Theorem 1.2.

Next observe that the commutation relations imply readily
that, for the
$\Ad$-diagonal
action, we have
$$ \bigl(\Ad \exp t\!\cdot\!(e_1+e_2)\bigr)\, \bh = \bh -
t\cdot\e\quad,\quad\forall t\in \C\;.\tag 1.5
$$

Define an algebraic family $\{\bh_t\in \ZZ\}$ parametrised by  points
$\, t\in \C^*\sqcup\{\infty\}=\C{\Bbb P}^1\smallsetminus \{0\}$ as follows. 
If $t\in\C^*$ put
 $\bh_t := \e-t^{-1}\cdot\bh$
$=-t^{-1}\cdot\bigl(\Ad
\exp t\cdot(e_1+e_2)\bigr)\, \bh ;$ and for $t=\infty$ put $\bh_\infty:=\e$.
Formula (1.5) shows that  $\bh_t
\to \bh_\infty$, as $t\to\infty$. 

To prove part (ii) of Theorem 1.2, we consider the family of vector spaces
$\z(\bh_t),
 \, t\in \C{\Bbb P}^1\smallsetminus \{0\}$. We claim that the dimensions of 
all the spaces $\z(\bh_t)$ in our family are independent of $t$.
For $t\in \C^*$, formula (1.5) shows that 
$\z(\bh_t)=\bigl(\Ad \exp t\!\cdot\!(e_1+e_2)\bigr)\, \z(\bh).$
Hence the centralizers $\z(\bh_t), t\in \C^*$, 
all have the same dimension
equal to $\dim\z(\bh)$. Now, by semicontinuity, the dimension of a special
member of any algebraic family can not be less than the dimension of
the general member. It follows that: 
$\dim \z(\bh)\leq \dim \z(\bh_\infty)$.
On the other hand, we know that $\dim \z(\bh_\infty)=
\dim \z(\e)
=\rk\g$, where the last equality is
due to the (Reg)-condition in (1.1). 
The claim follows.

Thus, we have proved that $\dim \z(\bh)=\rk\g$.
It remains to observe that the centralizer of the pair $(h_1,h_2)$
 of commuting semisimple
elements
is a Levi subalgebra in $\g$. Being of dimension $\rk\g$
this  Levi subalgebra must be a Cartan subalgebra.
\qed\medskip

Observe that our proof implies also that
in 
the Grassmannian of $\rk\g$-dimensional subspaces in $\g$ we have
$$\lim_{t\to\infty} \bigl(\Ad
\exp t\!\cdot\!(e_1+e_2)\bigr)\; \z(\bh)= \z(\e)\;.\tag 1.6
$$

Fix a principal nilpotent pair $\e=(e_1,e_2)$ and
an associated semisimple pair $\bh=(h_1,h_2)$.
 The adjoint action of commuting elements $(h_1,h_2)$
gives rise to a bigrading on $\g$:
$$ \g=\bigoplus_{(p,q) \in {\Bbb Q}\, 
\oplus\, {\Bbb Q}}\,\, \g_{p,q}\quad,\quad
\g_{p,q}=\{x\in\g\enspace |\enspace [h_1, x] =
p\cdot x\,,\; [h_2, x]= q\cdot
x\}\,.\tag 1.7
$$

We have: $e_1\in \g_{_{1,0}}\,,$ and
$e_2\in \g_{_{0,1}}$. Further, regularity of $\bh$ implies that
$\g_{_{0,0}}=\z(\bh)$
is a Cartan subalgebra in $\g$.
A routine argument shows now that:
$$ [\g_{p,q}\, ,\,\g_{p',q'}] = \g_{p+p', q+q'}\quad,
\quad \g_{p,q} \perp \g_{p',q'}\enspace\text{unless}\enspace
p=-p'\,\;\&\;\, q=-q',
$$
where `$\perp$' is taken with respect to the Killing form on $\g$.
This implies

\proclaim{Corollary 1.8} The Killing form induces a perfect pairing:
$\g_{p,q} \times \g_{-p,-q} \to \C.$\qed
\endproclaim
\smallskip

Observe that 
 $[h_i, \z(\e)] \subset \z(\e)\,,$
$i=1,2,$ hence the bigrading on $\g$ induces a bigrading
$\displaystyle \z(\e)=\bigoplus_{p,q}\,\z_{p,q}(\e)$. 

\proclaim{Proposition 1.9 (positive quadrant)} 
The algebra $\z(\e)$ is graded by the
`positive
quadrant', more precisely we have: $\;\displaystyle \z(\e)=\bigoplus_{p,q
\in\Z_{\,\ge 0},\,
(p,q)\neq(0,0)}\,\;\z_{p,q}(\e)$.
\endproclaim
\smallskip

{\it Proof.}
Observe 
that, for any $h\in \g_{_{0,0}}$ and any $t\in\C^*$,
all non-vanishing bigraded components of the element
$\bigl(\Ad
\exp t\!\cdot\!(e_1+e_2)\bigr)\,h$ are concentrated in
bidegrees $(p,q)\in\Z\oplus\Z$ such that $p,q\ge 0$. Hence, the same holds 
for any element in $\z(\e)$, due to (1.6).

It remains to prove that $\z_{0,0}(\e)=0$, equivalently,
that there is no non-zero element $h\in\g_{_{0,0}}=\z(\bh)$
that commutes with $(e_1,e_2)$. Assume that such an 
$h$ exists. 

We claim that $\z(\e) \subset \z(h)$. To prove this, consider
the commuting pair $(h+e_1,$
$ h+e_2)\in \ZZ$. The Richardson theorem [R] mentioned at the
beginning of this section 
yields
 the inequality: $\dim \z(h+e_1, h+e_2) \ge \rk\g$. On the other hand,
since $h$ and $e_i$ are respectively
the semisimple and nilpotent components
of $h+e_i,\, i=1,2,$ we have $\z(h+e_1, h+e_2)=\z(h)\cap
\z(e_1,e_2)$.
Therefore, if $\z(h)$ does not contain $\z(\e)$, we get
$$\dim \z(h+e_1, h+e_2) = \dim \bigl(\z(h)\cap
\z(e_1,e_2)\bigr) < \dim \z(e_1,e_2) = \rk\g\,.$$
This contradicts the Richardson's inequality above, and the claim
follows.

Next, we have a direct sum decomposition
$\g=\bigoplus_{\nu\in\C}\, \g(\nu)$ into $\ad h$-eigenspaces.
Clearly, for  each $\nu\in\C$,
the eigen-space $\g(\nu)$ is stable under $\ad e_1$ and $\ad e_2$,
since $e_1, e_2$ commute with $h$ by assumption. Hence, the operators
$\ad e_1$ and $\ad e_2$ form a commuting pair of nilpotent
endomorphisms of $\g(\nu)$. It follows that there is
a non-zero vector $v\in \g(\nu)$ annihilated by both endomorphisms.
Choosing $\nu$ here to be a non-zero eigenvalue, we conclude
that $\z(\e)\cap\g(\nu)\neq 0 $. This contradicts the inclusion
$\z(\e) \subset \g(0)=\z(h)$ proved in the preceding paragraph,
and the proposition follows.\qed
\medskip

From formula (1.6) and Proposition 1.9 we obtain:

\proclaim{Corollary 1.10} $\enspace\z(\e)$ is an abelian
Lie algebra consisting of
nilpotent elements.\qed
\endproclaim 
\medskip

For example, let $\g=\sln$ and let the
pair $\e_\lambda=(e_1, e_2)$ be associated to a Young diagram
$\lambda$ as explained in the Introduction.
As will be shown in
\S5 below, the Lie algebra $\z(\e)$ has a basis formed by the matrices:
$e_1^{\,\,p}\cd e_2^{\,\,q}$, 
where the pairs $(p,q)$ run over  the coordinates of
all boxes of $\lambda$, except for $(p,q)= (0,0)$, since $e^{\,\,0}_1\cd
e_2^{\,\,0}
=\op{Id}\not\in\sln$. Thus, we 
observe that 
the set of points $(p,q)$ 
on the 2-plane  such that $\dim\z_{p,q}(\e_\lambda) \neq 0$
 forms a figure of  shape
 $\lambda$. 

Motivated by this observation, and also by [K1, K2],
given a principal nilpotent pair $\e$ in any
semisimple Lie algebra $\g$, 
we  define {\it biexponents} of $\g$ relative
to $\e$ as the pairs of non-negative integers
$(p,q)$, cf. Proposition 1.9, such that $\dim\z_{p,q}(\e) \neq 0$.
In more detail, fix
 a bi-homogeneous base $z_1,\ldots,z_r\,,\,
r=\rk\g$, of the centralizer $\z(\e)=\bigoplus\;\z_{p,q}(\e)$.
To each $i= 1,\ldots,\rk\g$, assign
the pair $(p_i, q_i)\in \Z^2_{\geq 0}$
such that $z_i\in \z_{p_i,q_i}(\e)$. We introduce the following

\medskip

\noindent
{\bf Definition 1.11.} The subset ${\op{Exp}}_\e(\g)=
\{(p_1, q_1), \ldots, (p_r, q_r)\,,\,
r=\rk\g\}\;\subset\,\Z_{\geq 0}^2$ 
will be referrred to as the collection of {\it biexponents}
of $\g$ relative to $\e$.

\medskip

A new feature of the  "doubled" setup under investigation,
 making it quite
different from the "classical" one, is 
that it is impossible, in general, to
find $\ll2$-triples associated with  $e_1$ and
$e_2$ in such a way that they commute
with each other. As a result,
"{\it hard} Lefschetz" type equations, like $\dim \g_{p,q}
=\dim \g_{-p,q}$, are typically false in the bigraded setup.
The proposition below says that an analogue of
the "{\it weak} Lefschetz theorem" still holds in our situation.

\proclaim{Proposition 1.12 (weak Lefschetz)} For any principal pair $\e=(e_1,e_2)$,
the map $\ad e_1 : \g_{p,q} \to \g_{p+1,q}$ is { injective}
whenever $p < 0$, and is { surjective} whenever $p\geq 0$.
Similarly, the map
$\ad e_2: \g_{p,q} \to \g_{p, q+1}\,$ is { injective}
whenever $q< 0$, and is { surjective} whenever $q\geq 0$.
\endproclaim

{\it Proof.} For any $x\in \g$, the operator $\ad x: \g\to\g$ is
skew-adjoint with respect to the Killing form. Taking $x=e_1$
and using Corollary 1.8 we see that, for any $(p,q)$,
the adjoint of the operator on the left (below) equals, up to sign,
the one on the right:
$$\ad e_1 : \g_{p,q} \to \g_{p+1,q}\qquad,\qquad
\ad e_1 : \g_{-p-1,-q} \to \g_{-p,-q}\,.
$$
Since a linear map is surjective if and only if the dual map
is injective, it suffices to prove  the injectivity
part of the proposition only.

Set $\g_{p,*}=\bigoplus_q\,\g_{p,q}$.
We must prove that the
map $\ad e_1 : \g_{p,*} \to \g_{p+1,*}$ is injective for any
$p < 0$. If this map has a non-trivial kernel,
$\Ker_p(\ad e_1)\subset \g_{p,*}$, then the adjoint action
of $e_2$ gives an operator
$\ad e_2:\Ker_p(\ad e_1)\to\Ker_p(\ad e_1)$. This operator is nilpotent,
hence,
has a non-trivial kernel $= \Ker_p(\ad e_1)\cap \Ker(\ad e_2)$
$\neq 0$.
It follows that $\z(e_1,e_2)$ has a non-zero intersection
with $\g_{p,*}$. But this contradicts Proposition
1.9 ("positive quadrant"), because we assumed that $p<0$.\qed
\medskip

{\it Proof of parts $\op{(iii), (iv)}$ of Theorem 1.2.} 
Proving integrality of the eigenvalues of the operators
$\ad h_i$ amounts to showing that, in the bigrading
$\g=\bigoplus\,\g_{_{p,q}}$,
we have $\g_{_{p,q}}=0$ unless $(p,q)\in \Z\oplus\Z$.
Assume there exists a pair 
$(p_\circ, q_\circ)\not\in \Z\oplus\Z$
such that $\g_{p_\circ, q_\circ}\neq 0$.
Set $\Lambda= \{(p,q)\;|\; p\in p_\circ +\Z\;\,\&\;\, q\in q_\circ +
\Z\}\,,$
 and let $\g_{_\Lambda} =$
$\bigoplus_{(p,q)\in\Lambda}\,\,\g_{p,q}$. Clearly, $\g_{_\Lambda}$
is a subspace in $\g$ stable under the action of
$\{\ad e_i\}_{i=1,2}$. Hence, there is a non-zero vector $x\in
\g_{_\Lambda}$
annihiladed by both operators (which are nilpotent and commute).
Thus,
$\z(\e)\cap \g_{_\Lambda}\neq 0$. This contradicts Proposition 1.9,
saying that $\z_{p,q}(\e)$ may be non-zero only for integral $(p,q)$,
and part (iii)  of  Theorem 1.2 is proved.

To prove part (iv) fix an associated semisimple pair $\bh=(h_1,h_2)$,
and let $\langle h_1,h_2 \rangle$ be the $\C$-linear span
of $h_1,h_2$.

We consider the algebraic group $M$ 
from the exact sequence (1.4).
Let $\overline{M}\subset G$ be the image of $M$  under the 3-d
projection: $\C^*\times\C^* \times G \to G$,
and  $\overline{\m}:=\Lie\overline{M}$. Thus, we have
$\overline{\m}=\{x\in\g\;|\; \ad x(e_i) \in\C\cdot e_i\,,\;i=1,2\}.$
Exact sequence (1.4) yields a Lie algebra semi-direct product
decomposition
$\overline{\m} = \langle h_1,h_2 \rangle \ltimes \z(\e)$. Proposition 1.9 implies that
$\z(\e)$ is the nilradical (= Lie algebra of the unipotent radical)
 of the algebraic Lie algebra $\overline{\m}$.
Thus, $\overline{\m}$ is a solvable Lie algebra, and
$\langle h_1,h_2 \rangle$ is a maximal diagonalizable subalgebra in $\overline{\m}$.

Now, given another associated semisimple pair $\bh'=(h'_1,h'_2),$
we will get the same way another maximal diagonalizable subalgebra
$\langle h'_1,h'_2 \rangle$ in $\overline{\m}$. Recall that any
two maximal diagonalizable subalgebras of an algebraic solvable
Lie algebra
are conjugate by a unipotent element. We deduce that there exists a unipotent
element
$u\in Z_G^\circ(\e)$ 
such that $\Ad u \langle h_1,h_2 \rangle =\langle h'_1,h'_2 \rangle$. 
Thus, we may assume without
loss of generality that $\langle h'_1,h'_2 \rangle$
$=\langle h_1,h_2
\rangle$. But then the commutation relations 1.2(i) uniquely
determine the
positions of $h_1,h_2$ inside the 2-dimensional space
$\langle h_1,h_2 \rangle$, hence
$h'_i=h_i$, for $i=1,2$. This completes the proof of the Theorem.
\qed
\medskip

\medskip
The following result providing an alternative characterisation of
principal nilpotent pairs will be proved later in Section 2.

\proclaim{Theorem 1.13} A commuting pair $\e=(e_1,e_2)\in \ZZ$
is a principal nilpotent pair if and only if the following two
conditions
hold:

$\op{(a)}$ There exists a {\it regular} semisimple pair $\bh=(h_1,
h_2)\in \ZZ$ 
such that:
$$
 [h_i, e_j]= \delta_{i,j}\cdot e_i\quad,\;
i,j\in \{1,2\}\,,\quad\text{cf. 1.2(i)}\;;$$

$\op{(b)}$ The subalgebra $\z(\e)$ is graded by the
`positive
quadrant' 
with respect to the corresponding bigrading $\dis\g=\oplus_{p,q}\,
\g_{_{p,q}}$; 
that is:$ \quad\dis\z(\e)\;\subset\;
\bigoplus_{p,q
\in\Z_{\,\ge 0},\,
(p,q)\neq(0,0)}\,\;\g_{_{p,q}}\,.$
\endproclaim

Note that regularity of the pair $\e$ is {\it not} assumed in 
conditions (a)-(b) above, and that the grading on $\z(\e)$ is
required in (b) to be {\it integral}. 

\bigskip

\head{2. Bi-filtration associated to a nilpotent pair.}
\endhead
\bigskip

Given a vector space $V$ equipped with an increasing $\Z_{\geq 0}
\bigoplus
\Z_{\geq 0}$-filtration
(referred to as a "bi-filtration"),
$F_{\!_{i,j}}V$,
 such that $F_{\!_{i',j'}}V \subset F_{\!_{i,j}}V$,
whenever $i'<i$ or $j'<j$, we
define the corresponding associated bigraded space as follows:
$$
\operatorname{gr}^F\!V := \bigoplus_{i,j}\;\gr_{i,j\!}V\quad,
\quad
\gr_{i,j\!}V=\frac{F_{_{i,j}}V}{F_{_{i-1 , j}}V + F_{_{i , j-1}}V}\;.\tag 2.1
$$

To any commuting pair $(e_1,e_2)$ of nilpotent endomorphisms
of the vector space $V$ one associates the bi-filtration:
$F_{\!_{i,j}}(V) :=$
$\displaystyle \{v\in V\enspace|\enspace
e_1^{i+1}e_2^j(v)=0=
e_1^ie_2^{j+1}(v)\}.$
This bi-filtration on $V$ induces, by restriction, a bifiltration on
any vector subspace $E\subset V$, and we write
$\operatorname{gr}\!E$ for
the corresponding associated bigraded space.

Let $\dim E=r$, so that $E$ gives a point in
$\Gr_r(V)$, the Grassmannian of $r$-planes in $V$.
The natural action on $\Gr_r(V)$ of the one-parameter subgroup
$t\mapsto $
$\exp t\cd(e_1+e_2) \in \GL(V)$
gives 
an algebraic path $\gamma: \C\to \Gr_r(V)$ through the point $E$;
explicitly, $\gamma(t) = \exp t\cdot(e_1+e_2)(E)$. The Grassmannian being
a compact variety, the path $\gamma$ has a well-defined limit
$\gamma(\infty)\in \Gr_r(V)$, as $t\to\infty$.
The limit point  $\gamma(\infty)$ corresponds to
an $r$-dimensional vector subspace in $V$, to be denoted
$\tlim_\e E$, where $\e$ indicates the commuting pair $e_1,e_2$. 

Exploiting an idea of R. Brylinski 
one can give a more concrete description of
 the space $\tlim_\e E$ in terms of the bi-filtration 
$F_{_{p,q}}E$,
 which reduces in the special case $e_2=0$ to [Br, Lemma 2.5]:

\proclaim{Lemma 2.2} Assume that the sum
 $\sum_{p,q\geq 0}\;e_1^pe_2^q(F_{_{p,q}}E)\,$ is a direct sum
of vector subspaces
in $V$. Then

$\op{(i)}$ We have: $\tlim_\e E\;=\;\bigoplus_{p,q\geq 0}\;e_1^pe_2^q(F_{_{p,q}}E)\,;$
in
particular the space $\tlim_\e E$ is annihilated by the
operators $e_1,e_2$.

$\op{(ii)}$
The direct sum (over all $p,q\geq 0$) of maps
$e_1^pe_2^q: F_{_{p,q}}E \to \tlim_\e E$  induces
a vector space isomorphism
$$\text{\it lim}_\e\,\,:\enspace\,\operatorname{gr} E 
=\bigoplus_{p,q}\;  \frac{F_{_{p,q}}E}{F_{_{p-1 , q}}E + 
F_{_{p , q-1}}E}\;\enspace\overset {\sim} \to 
\longrightarrow\;\enspace 
\tlim_\e E\,.
$$

\endproclaim
\smallskip

To prove the lemma, note first that,
for any $x\in F_{_{p,q}}E$, one has a finite 
expansion
$$\bigl(\exp t\cd(e_1+e_2)\bigr)x =\sum_{i=0}^p\sum_{j=0}^q\;
\frac{t^{i+j}}{i!j!}\cdot e_1^ie_2^j(x)\,.$$
Observe next that
the action on $V$ of the operator
$e_1^pe_2^q$ kills the subspace
$F_{_{p-1 , q}}E + F_{_{p , q-1}}E$.
From this, one proves by induction on $(p,q)$
that $e_1^pe_2^q(F_{_{p,q}}E) \subset \tlim_\e E$.
Furthermore, the induced map $ e_1^pe_2^q:$
$ {F_{_{p,q}}E}/
\bigl(F_{_{p-1 , q}}E + F_{_{p , q-1}}E\bigr)\;$
$\to\,\tlim_\e E$ is injective.
Using  
that $\dim E=\dim (\tlim_\e E)$, one completes the proof
by a similar induction. \qed\medskip

Now let $\g$ be a semisimple Lie algebra.
Fix a principal nilpotent pair $\e=(e_1,e_2)$, and an associated
semisimple pair $\bh=(h_1,h_2)$. The construction above gives a
bi-filtration
on $\g$, associated to the commuting pair of nilpotent endomorphisms
$\ad e_1\,,\ad e_2$ on $V=\g$. We first consider the subspace
 $E:=\h=\z(\bh)\subset \g$.
Our construction gives a subspace $\tlim_\e
\h\subset\g$. Note that the assumption of Lemma 2.2 holds
since: $e_1^p e_2^q(\h)\subset \g_{_{p,q}}$,
and
the sum $\bigoplus  \g_{_{p,q}}$ is direct. Hence
the subspace $\tlim_\e\h$ is annihilated by $\e$, by Lemma 2.2(i).
We conclude that
$\tlim_\e
\h\subset\z(\e)$. But $\dim\h=\dim\z(\e)$, since $\e$ is a regular pair.
Thus, the inclusion implies an equality: $\tlim_\e
\h=\z(\e)$, which is  just another form of (1.6).

Next, we take $E:=\z(h_1,e_2)$.
Since the operator
$\ad e_2$ annihilates the space $\z(h_1,e_2)$, the bifiltration
$F_{_{i,j}}\z(h_1,e_2)$ reduces to an ordinary filtration so that,
for any $i,j\geq 0$,
we have: $F_{_{i,j}}\z(h_1,e_2)=
\text{Ker}\bigl(\ad^{i+1}e_1:\z_{_{0,*}}(e_2) \to 
\z_{_{i,*}}(e_2)\bigr)\,.$ The latter space being independent
of $j$, we will simply write $F_i:=F_{_{i,j}}\z(h_1,e_2)$. 

As we know, the map $\ad^{i}e_1$ gives an imbedding
$\ad^{i}e_1:$ 
$ F_i/F_{i-1}\hookrightarrow\,\z(\e)$. 
Furthermore, the assumption of Lemma 2.2 trivially holds.
Hence,  in the Grassmannian of
$\rk\g$-planes,
there exists a limit of the family of spaces
$\{(\Ad\exp t\cdot e_1)\,\z(h_1,e_2)\}_{t\to\infty}$,
 and by Lemma 2.2(i) one has
$\tlim_\e\, \z(h_1,e_2)\subset \z(\e)$. 
We note further that $(h_1,e_2)\in\ZZ$ is a commuting pair.
Hence
the Richardson inequality and the inclusion above yield:
$$\rk\g\;\leq\; \dim\z(h_1,e_2)\; =\;
\dim\bigl(\tlim_\e\,\z(h_1,e_2)\bigr)\; \leq \;\dim\z(\e)= \rk\g\,.$$
Thus, the inclusion: $\tlim_\e\, \z(h_1,e_2)\subset \z(\e)$ must be an equality, 
and we obtain 
the following analogue of (1.6):
$$\!\!\tlim_\e\,\z(h_1,e_2)=\lim_{t\to\infty} \bigl(\Ad
\exp \,t\!\cdot\!e_1\bigr)\;\z(h_1, e_2) 
\;=\;\bigoplus_{j\ge 0}\;\ad^{j}e_1\bigl(\z(h_1,e_2)\bigr)
\,= \z(\e).\tag 2.3
$$
\medskip

{\it Proof of Theorem 1.13.} We already know that if $\e$ is
a principal nilpotent pair, then properties (a)-(b) of Theorem 1.13
hold (by Theorem 1.2(i) and Proposition 1.9).

Thus, we must only
show that any pair $\e\in\ZZ$ satisfying the conditions
of Theorem 1.13 is a principal nilpotent pair. 
Condition 1.1(Nil) follows
readily from the very existence of an associated semisimple pair $\bh$
with the commutation relations
as in Theorem 1.13(a). Thus, it suffices to show that $\e$ is a  regular
pair.

To this end, we note first that the only ingredient used in the
proof of Proposition 1.12 was the "positive quadrant" property,
which is just condition (b) of Theorem 1.13. Therefore, the
weak Lefschetz property holds for our pair $\e$, even though we don't
know
yet that it is a principal nilpotent pair. 

Write $e_i^r$ for
the $r$-th power of the adjoint action of $e_i$ on $\g$.
The surjectivity part of the weak Lefschetz
yields: $\bigoplus_{p,q\in\Z_{\ge 0}}\;\g_{_{p,q}}=
\bigoplus_{p,q\in\Z_{\ge 0}}\;e_1^pe_2^q(\g_{_{0,0}})\,.$
Using the
"positive quadrant" condition of Theorem 1.13(b) once again
we deduce that
$\dis\,\z_{_{p,q}}(\e)\subset
e_1^pe_2^q(\g_{_{0,0}})\;,\;\forall p,q\in\Z_{\ge 0}\,.$
Hence, given
$x\in \z_{_{p,q}}(\e)$, one can find $h\in \g_{_{0,0}}$ such that
$x=e_1^pe_2^q(h)$. The  condition $x\in\z_{_{p,q}}(\e)$ thus
reads: $0=e_1(x)=e_1^{p+1}e_2^q(h)$, and
$0=e_2(x)=e_1^pe_2^{q+1}(h)$.

By the assumptions of Theorem 1.13, the pair $\bh$ is
 regular, hence $\h= \g_{_{0,0}}$ is a Cartan subalgebra in $\g$.
We have the bi-filtration $F_{_{i,j}}\h$ on $\h$ arising
from the adjoint action of the commuting pair $(e_1,e_2)$ on $\g$.
The equations at the end of the previous paragraph
 say that $h\in F_{_{p,q}}\h$. Hence, we have proved
that $\,\z(\e)=\bigoplus_{p,q\geq 0}\;e_1^pe_2^q(F_{_{p,q}}\h)\,.$
On the other hand, in the special case $V=\g$ and $E=\h$,
Lemma 2.1(i) reads:
$\tlim_\e \h =$
$\bigoplus_{p,q\geq 0}\;e_1^pe_2^q(F_{_{p,q}}\h)\,.$
Thus, $\,\z(\e)=\bigoplus\,e_1^pe_2^q(F_{_{p,q}}\h)
=\tlim_\e \h,$ and we deduce that
$\dim \z(\e)= \dim(\tlim_\e \h)=\dim\h\,.$
This proves regularity of $\e$, and the theorem follows.
$\quad\square$
\medskip

Now let $V$ be a finite dimensional $\g$-module.
We write $V^{\m}$ for the subspace of $V$ annihilated by
a Lie subalgebra $\m\subset\g$.
Let $T$ be the maximal torus in $G$ corresponding
to the Cartan subalgebra $\h=\z(\bh)$.
Given
a weight $\mu\in X^*(T):=\Hom_{_{\text{alg group}}}(T,\C^*)$, let
$V(\mu)$ denote the corresponding weight space.
 We claim, inspired by [Br, Proposition 2.6], that
$$\tlim_\e\,V(\mu) \subset V^{\z(\e)}\quad,\quad\forall \mu\in X^*(T)\,.\tag 2.4
$$
To prove this inclusion observe that every non-zero vector in
$V(\mu)$ gives  a point in ${\Bbb P}\bigl(V(\mu)\bigr)$ (= projectivisation of
 $V(\mu)$) fixed by the $T$-action,
hence, by the infinitesimal $\h$-action. It follows, by continuity,
that each point of ${\Bbb P}\bigl(\tlim_\e V(\mu)\bigr)$ is fixed 
by the infinitesimal $\tlim_\e\h$-action. Hence
$\tlim_\e V(\mu)$ is a weight space for the subalgebra
$\tlim_\e\h$. But since the subalgebra $\tlim_\e\h=\z(\e)$ consists
of nilpotent elements, the weight in question must vanish,
and our  claim follows.

As an important consequence of (2.4), assume $V$ is
a rational
$G$-module, so that all the weights of $V$ belong to the root lattice of
$\g$,
in particular, zero is a weight of $V$ and $V^{\z(\bh)}\neq 0$.
Then, the subspace $\tlim_\e\, V^{\z(\bh)}\subset V^{\z(\e)}$
is {\it independent} of the choice of an associated semisimple pair
$\bh$.
Indeed, by Theorem 1.2(iv),
any two such semisimple pairs are conjugate by  the unipotent group
$Z_G^\circ(\e)$.
Hence, the corresponding Cartan subalgebras $\h=\z(\bh)$
are conjugate by $Z_G^\circ(\e)$. It follows 
that the corresponding zero-weight spaces $V^{\z(\bh)}$, whence the limit-spaces
$\tlim_\e\, V^{\z(\bh)}$, are conjugate by the $Z_G^\circ(\e)$-action.
But the latter action is trivial by (2.4), thus, the two limit-spaces
coincide. Note that if $\dim V^{\z(\bh)} =\dim V^{\z(\e)}$, then (2.4) becomes
an equality.
We emphasize that unlike the case of an ordinary principal nilpotent,
the inclusion $\tlim_\e\,
 V^{\z(\bh)}\subset V^{\z(\e)}$ may be {\it strict},
in general. This happens, for instance, for the 
pair $\e$ in $\g={\frak s}{\frak l}_3$, corresponding to the "hook"
Young diagram, and
$V=S^3(\C^3)$, the third Symmetric power of the fundamental
representation.

In spite of this counter-example we expect that one of the most
important results of [K2] continues to hold true in the double
setup. Specifically, let $\overline{\O}_\e$
denote the closure of the $\Ad G$-diagonal orbit of $\e$,
an affine subvariety in $\g\oplus\g$.
Write $\C[X]$ for the coordinate ring
of an affine algebraic variety $X$. Then we have

\proclaim{Conjecture 2.5} There is a $G$-module isomorphism:
$\C[\overline{\O}_\e] \simeq \C[G/T]$.
\endproclaim

\medskip

Recall next that, associated to $\e$, we have defined the bifiltration
$F_{\!_{\bullet,\bullet}}V$
 on a finite dimensional $\g$-module
$V$.
Given
a weight $\mu\in X^*(T)$, we endow the weight space $V(\mu)$
with the induced filtration
 $F_{\!_{i,j}}V(\mu) = V(\mu) \cap F_{\!_{i,j}}(V)$,
and let $\,
\operatorname{gr}V(\mu) = \bigoplus_{i,j}\;\gr_{i,j\!}V(\mu)\,$
denote the corresponding associated bigraded space, see (2.1).
By Theorem 1.2(iv)
the integers $\dim \gr_{i,j\!}V(\mu)$, called
{\it bi-exponents} of $V$ relative to $\e$, cf. Definition 1.11,
are canonically associated to
$\e$ and $V$, i.e., do not depend on the choice of an
associated pair $(h_1,h_2)$.

We define an $(s,t)$-{\it weight multiplicity} of $V(\mu)$ as the Poincar\'e
polynomial:
$$ P_\mu(V,\e) =\sum_{i,j\geq 0}\,s^it^j \cdot\dim \gr_{i,j\!}V(\mu)
\,.\tag 2.6$$
Our goal is to produce
an explicit formula for the polynomials $P_\mu(V,\e)$,
analogous to Lusztig's $q$-weight multiplicity [L1], see also
[Br].  

First, write $R$ for the root system of $\g$ with respect to the
Cartan subalgebra $\h=\z(h_1,h_2)$. There is a distinguished
subset
$R_{\ne}\subset R$ (where {\bf {ne}} = {\it northeast})
formed by the roots
 that occur in the ``positive quadrant'': 
$\,R_{\ne}=\{\alpha\in R\;|\;
\alpha(h_1) \ge 0\;\&\;\alpha(h_2) \ge 0\}\,.$ The set $R_{\ne}$
is clearly contained in an open half-space of $\h^*$.
Thus, we can (and will)
make a choice of the
set $R_+\subset R$ of positive roots in such a way that
$R_{\ne}\subset R_+$. By construction, we have a decomposition
$R_+ = R_{\ne}\,\sqcup\,(R_+\smallsetminus R_{\ne})$.
Furthermore, the set $R_{\ne}$ has a natural decomposition
into a disjoint union of  three subsets:
$$
\align
&R_{\ne}= R_+^1 \,\sqcup\,R_+^2\,\sqcup\,\Rsmall
\quad\text{where}_{_{}}\quad
\Rsmall= \{\alpha\in R\;|\;
\alpha(h_1) >0\;\&\;\alpha(h_2) >0\},\tag 2.7\\
&R_+^1= \{\alpha\in R\;|\;
\alpha(h_1) >0 = \alpha(h_2)\}^{^{}}\enspace,\enspace
R_+^2= \{\alpha\in R\;|\;
\alpha(h_1) =0 < \alpha(h_2)\}\,.
\endalign
$$

Recall that since the group $G$ is of adjoint type, the
weight lattice $X^*(T)$ may be identified with the root
lattice. Let  
$Q_{\ne}$ be the sub-semigroup of $X^*(T)$
generated by the set $R_{\ne}$.

Write $e^\beta : T\to \C^*$ for the group
character corresponding to a root $\beta \in R$.
For each $s,t \in \C^*$, 
we define a rational function
$\pp_\e$ on the torus $T$
by  the following product: 
$$\pp_\e=\!\!
\!\prod_{\beta\in R_+ \smallsetminus R_{\ne}}\!\!\!\!\!(1-e^{\beta})\,\,
\cdot\prod_{\mu\in R^{_1}_+}\!(1-s\, e^{\mu})\,\cdot 
\prod_{\nu\in R^{_2}_+}\!(1-t\, e^{\nu})\,\cdot
\prod_{\alpha\in \Rsmall}\!
\frac{(1-s\, e^{\alpha})(1-t\, e^{\alpha})}{(1-st\cd e^{\alpha})}\,.
\tag 2.8
$$ 
Thus, the assignment: $s,t,z \mapsto \pp_\e(s,t,z)$
gives a rational function on $\C^*\times\C^*\times T$.
It is clear that the function $1/\pp_\e$
 has the following expansion:
$$\frac{1}{\pp_\e}\,=\,
\sum_{\alpha\in Q_{\ne}}\;\wp_{\!_\e}(\alpha)\cdot e^\alpha\,,
\quad\text{where}\quad\wp_{\!_\e}(\alpha)\in \Z_{_{\geq 0}}[s,t]\,.
$$
Note further that if $s=t=1$ then the product (2.8)
specializes to the classical Weyl denominator:
$ \pp_\e\big|_{_{s=t=1}} = \prod_{\alpha\in R_+}\;
(1-e^\alpha)\,,$ so that
the function $\wp_{\!_\e}$ becomes the Kostant's partition function.
Motivated by the classical case,
 we propose  
\medskip

\noindent
{\bf Definition 2.9.} The function $\wp_{\!_\e}:\,
Q_{\ne} \to \Z_{\geq 0}[s,t]\,,$
$\alpha\mapsto \wp_{\!_\e}(\alpha)\,,$ is called  double-analogue
of the Kostant partition function associated to the principal pair~$\e$.
\medskip

Note that for a principal nilpotent pair of the form
$\e=(e,0)$ or $\e=(0,e)$ the function $\wp_{\!_\e}$
reduces to the $q$-analogue of Kostant's partition function,
introduced in [L1].
 
Let $\rho$ be the half-sum of positive roots in $R_+$.
Write $W$ for the Weyl group of $(\g,\h)$, and let
${\sign}: w\mapsto \eps(w)= (-1)^{\text{length} (w)}$ denote
the determinant of the $w$-action on $\h$.
In the sequel to this paper we will prove
the following double-analogue of Kostant's weight multiplicity formula,
 expressing the
$(s,t)$-weight multiplicities  in terms
of the partition function $\wp_{\!_\e}$:

\proclaim{Theorem 2.10} Let $\lambda$
be a dominant weight  such that
$\lambda\in Q_{\ne}$, and $V_\lambda$ the simple $G$-module with highest weight
$\lambda$.
Then, one has:
$$P_\mu(V_\lambda,\e) = \sum_{w\in W}\; 
\varepsilon(w)\cd\wp_{\!_\e}(w(\lambda\!+\!\rho)-\mu-\!\rho)\;,\;\;
\text{for any weight}\;\mu\in X^*(T)\,.
$$
\endproclaim

In the case of $\e=(e,0)$ where $e$ is an ordinary principal
nilpotent,  the theorem above reduces to a result
of R. Brylinski [Br].
Our proof consists of several steps
which are analogous to the corresponding steps in [Br],
but with considerable complications due to the difference between
the classical geometry of conjugacy classes in $\g$ and that
of double-orbital varieties to be introduced and
studied in the subsequent paper of this series.
 Later on, we are also going to give an alternative
geometric interpretation of the polynomials $P_\mu(V_\lambda,\e)$
in terms of intersection cohomology of certain varieties
in a "double-loop" Grassmannian, cf. [L1], [Gi] for the affine case. 
\medskip

\bigskip

\head{3. Two parabolics.}
\endhead
\bigskip

Fix a principal nilpotent pair $\e=(e_1,e_2)$ and an associated
semisimple pair $\bh=(h_1,h_2)$. 
We introduce two Levi subalgebras
$$\g^{_1} := \z_\g(h_2) = \bigoplus_{p \in \Z}\, \g_{p,0}\;,\quad
\text{and}\quad
\g^{_2} := \z_\g(h_1) = \bigoplus_{q \in \Z}\, \g_{0,q}\,.\tag 3.1
$$
We will analyze the first one, $\g^{_1}$, the structure
of the second being entirely similar.

Let $\g^{_1}=\c(\g^{_1})\oplus\s^{_1}$ be the decomposition of $\g^{_1}$ into its center
and the semisimple derived Lie algebra. It is clear that
$e_1\in\s^{_1}$ and $h_1\in\g^{_1}$. Write $h_1=c_1+s_1$, where
$c_1\in \c(\g^{_1})$ and $s_1\in\s^{_1}$.

Let $\h=\z(h_1,h_2)= \g_{_{0,0}}$
be the  Cartan subalgebra in $\g$ associated to $\bh$.
It is clear that $\c(\g^{_1})\subset \h$.
Furthermore, the space $\h^{_1}:=\s^{_1}\cap\h =\z_{\s^{_1}}(s_1)$
is a Cartan subalgebra in $\s^{_1}$ since we have:
$\c(\g^{_1}) + \h^{_1}=\c(\g^{_1}) + \z_{\s^{_1}}(s_1) =\z_{\g^{_1}}(s_1) =$
$\z_{\g^{_1}}(h_1) =\h,$ is a Cartan subalgebra in $\g^{_1}$.
In particular, $s_1$ is a regular element of $\s^{_1}$.

Write $R\subset \h^*$ for the root system of $(\g,\h)$,
and $R^{1}, R^{2}\subset R$ for the root systems of
$(\g^{_1}, \h)$ and $(\g^{_2}, \h)$, respectively. 
For each $\alpha\in R$, let $\alpha^\vee\in\h$ denote the
corresponding
coroot, and
$x_\alpha\in \g$  the corresponding root vector.
Let
$R^{1}_+$ be the set of  roots $\alpha \in R^{1}$ 
such that $\alpha(h_1) > 0$,
see $(2.7)$. Then $R^{1}_+$
is a system of positive roots for $R^{1}$, and we have the corresponding
triangular decomposition
$\g^{_1}= \g^{_1}_{_+} \,\oplus\, \h\,\oplus\,\g^{_1}_{_-}\,.$
Let $\Delta^{1}\subset R^{1}_+$ be the 
set of simple roots.

\proclaim {Proposition 3.2}$\,\op{(i)}\;
\g^{_1}_{_\pm} $ are the nilradicals of two opposite Borel 
 subalgebras in~$\s^{_1};$

$\op{(ii)} $ The element $s_1$ equals
the half-sum of positive coroots of $\h^{_1}$, that is:\par
 $\; s_1= \frac{1}{2}\sum_{\alpha\in R^{1}_+}\,{\alpha^\vee}\,;$

$\op{(iii) }\enspace e_1\in \s^{_1}_+$ is a principal nilpotent in $\s^{_1}$;
furthermore, 
$e_1 = \sum_{\alpha\in \Delta^{^{\!1}}}\,t_\alpha\cdot x_\alpha\,,$\par
for certain $t_\alpha\in \C^*\,;$

$\op{(iv) }$ The elements $s_1, e_1$ are 
members of a principal
$\ll2$-triple $\langle  s_1, e_1, f_1\rangle\;$ in~$\s^{_1}$.
\endproclaim

{\it Proof.} Part (i) has been already proved.
To prove that $e_1$ is a regular element in $\s^{_1}$ we note that
$\rk\g = \rk\g^{_1}$ and
$\dim \z_{\g^{_1}}(e_1) =
\dim \z(e_1, h_2)$. But now isomorphism (2.3) yields
$\dim \z(e_1, h_2) = \rk\g$. Thus, $e_1$ is regular in
$\g^{_1}$, hence, in $\s^{_1}$. 

To get the expression for $e_1$ given in
(iii), we  write $e_1$ as a linear
combination of root vectors:
$e_1 = \sum_{\alpha\in R^{_1}_+}\, t_\alpha\!\cdot\!x_\alpha$,
for some $t_\alpha\in \C$.
It is known that, since $e_1$ is  a regular nilpotent,  
the coefficient $t_\alpha$ corresponding to every
simple  root $\alpha\in\Delta^{_1}$ must be non-zero. Furthermore,
 $e_1$ is an eigenvector for $\ad s_1$. But since  $s_1 \in \h^{_1}$ is
regular, any $\ad s_1$-eigenvector having a non-zero component for each simple
root vector must be a linear combination of simple
root vectors (with non-vanishing coefficients). This
proves (iii).

The commutation relation $[s_1, e_1]=e_1$, where
$e_1 = \sum_{\alpha\in R^{_1}_+}\, t_\alpha\cdot x_\alpha$,
with all $t_\alpha\neq 0$, forces
$s_1\in \h^{_1}$ to be equal to the half-sum of positive coroots:
$\,s_1 = \frac{1}{2}\sum_{\alpha\in R^{_1}_+}\,{\alpha^\vee}\,.$
This proves (ii). Finally, it is known
(and straightforward to check directly) 
that the elements $s_1, e_1$ given by the expressions above
are members of a principal $\ll2$-triple, and part (iv)
follows.
\qed
\bigskip

Write $\displaystyle \g_{p,*} :=
\oplus_q\,\,\g_{p,q}\,$, and  $\displaystyle \g_{*,q} :=
\oplus_p\,\,\g_{p,q}\,$. Thus, $\g^{_1}=\g_{_{*,0}}$, and
$\g^{_2}=\g_{_{0,*}}$. 
We introduce two parabolic subalgebras in $\g$ 
with Levi factors $\g^{_2}$ and $\g^{_1}$, respectively:
$$
\p\east=\,\g^{_2}\oplus \bigl(\bigoplus_{p>0}\,\g_{p,*}\bigr)
\, =\,\bigoplus_{p\geq 0}\,\g_{p,*}
\quad,\quad
\p\north=\g^{_1}\oplus \bigl(\bigoplus_{q>0}\,\g_{*,q}\bigr)
\, =\,\bigoplus_{q\geq 0}\,\g_{*,q}\,.\tag 3.3
$$
Thus, the parabolic
$\p\east$ is the sum of  $\g_{p,q}$ over all $(p,q)$ in  the right
half-plane of the $(p,q)$-plane, resp.
$\p\north$ is the sum of  $\g_{p,q}$ over all $(p,q)$ in
 the upper half-plane.

Let $G^{2}\,\subset\,P\east\,\subset\,G$ denote the subgroups  corresponding to
the subalgebras $\g^{_2}\,\subset\,\p\east\,\subset\,\g$, respectively. 
Note that $\z(\e)\subset
\p\east\cap\p\north$. Therefore, the
parabolics $\p\east,\,\p\north$ are completely determined by $\e$ and do {\it
not}
depend on the choice of an associated semisimple pair $\bh$, because
the latter is determined up to conjugacy by an element of
$Z_G^\circ(\e)=
\exp \z(\e)
\subset P\east\cap P\north$.

Recall that an element $e$ in   the Lie algebra $\p$
of a parabolic subgroup $P\subset G$ is called {\it Richardson}
for $\p$ 
if the orbit $\Ad P(e)$ is open dense in the nilradical of $\p$.
It is clear that
the space $\bigoplus_{p\ge 1}\,\g_{p,*}$ is the nilradical of $\p\east$, hence is
 stable under the $\Ad P\east$-action. Further,
the space $\g_{_{1,*}}$ contains the element $e_1$ and
is stable under the $\Ad G^{2}$-action.

\proclaim{Proposition 3.4} 
$\op{(i)}$ The $\Ad G^{2}$-orbit of $e_1\in\g_{_{1,*}}$ is Zariski open, dense
in $\g_{_{1,*}}$;

$\op{(ii)}$ The $\Ad P\east$-orbit of $e_1$ is Zariski open, dense
in $\bigoplus_{p\ge 1}\,\g_{p,*}$, 
in particular, \par
$e_1$ is a Richardson element for $\p\east$.
Similar results hold for $e_2$.
\endproclaim

{\it Proof.} To prove (i), it suffices to show that the tangent
space at $e_1$ to the $\Ad G^{2}$-orbit of $e_1$ equals $\g_{_{1,*}}$,
that is to show that $[e_1, \g^{_2}]=\g_{_{1,*}}$. But this is immediate
from the surjectivity claim of the weak Lefschetz. Part (ii)
is proved in exactly the same way.\qed
\medskip

\proclaim{Lemma 3.5}
$\op{(i)}\,\c(\g^{_1})\cap\c(\g^{_2})=0\,.\;\op{(ii)}$
The subspaces $\g^{_1}$ and $\g^{_2}$ generate $\g$ as
a Lie algebra. 
\endproclaim

{\it Proof.} Recall that the centralizer of
a principal $\ll2$-triple in a semisimple Lie algebra is trivial.
Hence, Proposition 3.2(iv) yields $\z_{\g^{_1}}(e_1, h_1)=\c(\g^{_1})$,
and $\z_{\g^{_2}}(e_2, h_2)=\c(\g^{_2})$. Part (i) now follows from Proposition
1.9 (positive quadrant)
and the equation
$$\c(\g^{_1})\cap\c(\g^{_2})\,=\, \z_{\g^{_1}}(e_1, h_1) \cap \z_{\g^{_2}}(e_2, h_2)\,=\,
\z(e_1, e_2) \cap \z(h_1, h_2)\,=\, \z(\e)\cap\h = 0\,.$$

To prove (ii), let $\widetilde\g$ denote the Lie subalgebra in $\g$
generated by $\g^{_1}$ and $\g^{_2}$.
The weak Lefschetz surjectivity result
implies that the parabolic subalgebra $\p\east$ is generated,
as a Lie algebra, by the subalgebra $\g^{_2}$ and the element $e_1$.
It follows that $\widetilde\g\supset \p\east$. Further, if $\sigma$ is
the Cartan involution of $\g$ with respect to $\h$, then 
clearly $\sigma(\g^{_1})=\g^{_1}$ , $\sigma(\g^{_2})=\g^{_2}$,
and $\sigma(\p\east)$ is the parabolic opposite to $\p\east$.
Therefore, the algebra
${\widetilde\g}=\sigma({\widetilde\g})$ contains both $\p\east$ and
$\sigma(\p\east)$, and we are done.\qed\medskip

\proclaim{Corollary 3.6}  $\op{(i)}$ For any semisimple pair $(h_1, h_2)$
associated to a principal nilpotent pair $\e$,
 the integers $m_1,m_2$ in Lemma $\op{1.3}$
can be  chosen to be equal to~$1$; 

$\op{(ii)}$  The group $Z_G(\e)$, the simultaneous
centralizer of $e_1$ and $e_2$ in $G$, is connected.
\endproclaim

{\it Proof.} Write $T$ for the maximal torus in $G$
corresponding to $\h$, and let $R\subset\h^*$ be the
root system of $(\g,\h)$. Recall that
the lattice $X_*(T)=\Hom_{alg}(\C^*, T)$ is dual to the
root lattice of $\g$,
since the group $G$ is of adjoint type. Theorem 1.2(iii)
implies that, for any root $\alpha\in R$, we have $\alpha(h_i)\in \Z$.
Thus, if we view $X_*(T)$ as a lattice in $\h$ via taking
the differentials of homomorphisms: $\C^*\to T$ at the identity,
then: $h_i\in X_*(T)\,,\,i=1,2$. It follows that
there exists $\gamma:\C^*\times\C^*\to T$ such that
$\frac{\partial\gamma}{\partial t_i}\big|_{t_1=t_2=1}$
$=h_i$, for $i=1,2.$ This proves part (i).

 To prove (ii) note that the connected component 
of the group $Z_G(\e)$ is a unipotent group, due to
Proposition 1.9 ("positive quadrant").
Hence, any maximal reductive subgroup in  $Z_G(\e)$ is a finite group.

Recall now the short exact sequence
(1.4). This sequence is split, by part (i) of the corollary. It follows,
that $\overline{M}$, the image of $M$ under the projection to $G$
considered in the
proof of Theorem 1.2(iv), has the form of a semidirect product
$\overline{M}=(\C^*\times\C^*)\ltimes Z_G(\e)$. Hence, if
$\overline{M}_{\bold r}\subset \overline{M}$ is a
maximal
 reductive subgroup containing  $\C^*\times\C^*$,
then $Z_{\bold r} := Z_G(\e) \cap \overline{M}_{\bold r}$ is a maximal
 reductive subgroup  in $Z_G(\e)$. In particular, $Z_{\bold r}$
is  a normal subgroup of $\overline{M}_{\bold r}$, and there
is a semi-direct product decomposition $\overline{M}_{\bold r} =$
$(\C^*\times\C^*)\ltimes Z_{\bold r}$. The adjoint 
$\C^*\times\C^*$-action on $Z_{\bold r}$ must be trivial, since
the group $Z_{\bold r}$ is finite. Hence, $Z_{\bold r}$
commutes
with the image of the homomorphism $\gamma:\C^*\times\C^*\to T$
constructed in the proof of part (i).
It follows that $Z_{\bold r} \subset T$, since $T=Z_G(h_1,h_2)$.

Let $z\in Z_{\bold r}$. 
Since $z\in T$ commutes with $e_1$
and $e_2$,
Proposition 3.2(ii) shows that, for any root $\alpha\in 
\Delta^{1} \cup \Delta^{2}$ (sets of simple roots for
$\g^{_1}$ and $\g^{_2}$),
we have $\alpha(z)=1$. By Lemma 3.5(ii) the set
$\Delta^{1}\cup\Delta^{2}$ generates the root lattice of $\g$.
It follows that $z$ commutes
 with $\g$. Thus,
$z$ belongs to the center of $G$, which is trivial because
$G$ is of adjoint type.\qed 
\medskip

\proclaim{Theorem 3.7} 
Principal nilpotent pair $\e=(e_1, e_2)$ is uniquely
determined, up to conjugacy, by the associated semisimple
pair $\bh=(h_1,h_2)$, that is, 
any two principal nilpotent pairs $\e$ and $\tilde{\e}$ that have
the same associated semisimple
pair $\bh$ are conjugate to each other by the maximal torus 
$T= Z_G(\bh)$.
\endproclaim
\medskip

We will use the following lemma; its proof is postponed until \S6
(Corollary 6.10).

\proclaim{Lemma 3.8} Set $\z_{_{p,q}}(e_1):=\g_{_{p,q}}\cap
\z_\g(e_1).$
 Then, $\;\z_{_{0,1}}(e_1)= [e_2\,,\,\z_{_{0,0}}(e_1)]\;.$
\endproclaim
\medskip

{\it Proof of Theorem 3.7.} Let $\e=(e_1, e_2)$
and $\tilde{\e}=(\tilde{e}_1, \tilde{e}_2)$ be
two principal nilpotent pairs with
the same associated semisimple
pair $\bh=(h_1,h_2)$.
Put $\h=\z(\bh)$ and define the Levi subalgebras
$\g^{_1}$ and $\g^{_2}$ as at the beginning of this section.
We have $\h=\c(\g^{_1})\oplus\h^{_1}$. Let $T^{1}=\exp(\h^{_1}) 
\subset T$ be the
torus  corresponding to the Lie
subalgebra $\h^{_1}$. Using Proposition 3.2(iii) we may, conjugating by
$T^{1}$ if necessary, achieve that 
$e_1=\tilde{e}_1$. We assume this, from now on.

Let $C=\exp(\c(\g^{_1}))\subset T$ be the torus corresponding
to the center of $\g^{_1}$. Conjugating by $C$ acts trivially on $\g^{_1}$,
hence does not affect $\bh$ and the equality $e_1=\tilde{e}_1$.
Therefore, it suffices to prove that $\tilde{e}_2$ is $\Ad C$-conjugate to
$e_2$. Note that, by construction,
$e_2,\,\tilde{e}_2 \in\z_{_{0,1}}(e_1)$, and that
the space $\z_{_{0,1}}(e_1)$ is $\Ad C$-stable. Thus, the theorem
will
follow provided we show that 
$\Ad C$-orbits of $e_2$ and $\tilde{e}_2$ are both Zariski
open in $\z_{_{0,1}}(e_1)$. To this end, we observe that
the $\Ad C$-orbit of $e_2$ is Zariski open in  $\z_{_{0,1}}(e_1)$
if and only if its tangent space, $\bigl(\ad\,\c(\g^{_1})\bigr) e_2$,
equals the whole space, i.e.: $ [e_2, \c(\g^{_1})] =\z_{_{0,1}}(e_1)$.
Since $\c(\g^{_1})=\z_{_{0,0}}(e_1)$,
the latter equation is insured by Lemma 3.8. Similar
argument applies to  $\tilde{e}_2$.
\qed
\bigskip
 
\medskip
\proclaim{Theorem 3.9} The number of $\Ad G$-orbits of 
principal nilpotent pairs in $\g$ is finite.
\endproclaim

{\it First proof.} By Theorem 3.7 we only need to show that there are
finitely
many $\Ad G$-conjugacy classes of all possible associated semisimple
pairs $\bh=(h_1,h_2)$. To this end, we fix a Cartan subalgebra $\h\subset \g$
and prove that the number of associated semisimple
pairs $\bh=(h_1,h_2)$ such that $\h=\z(\bh)$ is finite.

Given $\h$, let $\g^{_1},\,\g^{_2}\supset \h$ be
two Levi subalgebras that generate $\g$
as a Lie algebra and such that $\g^{_1}\cap\g^{_2}=\h$.
Choose $\Delta^{1}$ and $\Delta^{2}$, bases of simple roots for
$(\g^{_1},\h)$ and $(\g^{_2},\h)$, respectively. Note that there are
finitely many choices of quadruples
$(\g^{_1},\,\g^{_2},\,\Delta^{1},\,\Delta^{2})$ as above.
We now prove that there is at most
one associated semisimple
pair $\bh=(h_1,h_2)$ compatible with
such a quadruple.

To this end, assume that the quadruple $(\g^{_1},\,\g^{_2},\,\Delta^{1},\,\Delta^{2})$ 
and the semisimple
pair $\bh=(h_1,h_2)$ come from a
principal nilpotent pair $\e$, as at the beginning
of this section. Then $\g^{_1}=\c(\g^{_1})\oplus\s^{_1}$ and  $h_1=c_1+s_1$,  where 
$s_1 = \frac{1}{2}\sum_{\alpha\in R^{_1}_+}\,{\alpha^\vee}\,,$
by Proposition 3.2(ii).
Further, using part (iii) of Proposition 3.2,
 we can write $e_2=\sum_{\alpha\in
\Delta^{2}}\,t_\alpha\cdot
x_\alpha\,,$  for some $t_\alpha\in \C^*$.
Thus, equation $0=[h_1, e_2]= \ad h_1(e_2)$ yields:
$$ 
\alpha(c_1) +\alpha(s_1) = 0\quad,\quad\forall \alpha\in \Delta^{2}\,.\tag
3.10
$$

We reinterpret (3.10) without explicit mentioning of roots.
Let $p_{\s}$ and $p_{\c}$ denote the restrictions to the 
Lie algebras $\s^{_1}$ and $\c(\g^{_1})$, respectively,
of the first projection:
$\h=\s^{_2}\oplus\c(\g^{_2})\twoheadrightarrow \s^{_2}$.
Then equation (3.10) says: $p_{\c}(c_1)=-p_{\s}(s_1)=$
$-p_{\s}(\frac{1}{2}\sum_{\alpha\in R^{_1}_+}\,{\alpha^\vee})\,.$
But this last equation determines
$c_1$ uniquely, since the map $p_{\c}$ is injective, because
$\Ker(p_{\c}) = \c(\g^{_1})\cap\c(\g^{_2})=0$ (Lemma 3.5). Thus, $h_1=c_1+s_1$
is determined uniquely.\qed
\medskip

\noindent
{\bf Remark 3.11.} The above argument implies in particular
the following constraint on the
relative position of the  parabolics $\p\north$ and $\p\east$,
which may be helpful in a practical search for
principal nilpotent pairs. Let $\p''$ and $\p'$
be any two parabolics with Levi factors $\g^{_1}$ and
$\g^{_2}$, respectively. Assume  that $\g^{_1}\cap\g^{_2}=\h$
is a Cartan subalgebra of $\g$, and for $i=1,2$,
choose $\Delta^{i}$,
the set of simple roots for $(\g^{_i}, \h)$.
Further, let ${\frak b}^{_1}$ be the Borel subalgebra
in $\g^{_1}$ corresponding to the set $\Delta^{1}$.
Write $\Delta\supset \Delta^{1}$ for the set
of simple roots for $(\g, \h)$ corresponding to
the
unique Borel subalgebra ${\frak b}$
in $\g$ such that 
${\frak b}^{_1}\subset {\frak b}\subset\p''$.
Let $h_2$ be the sum of the
fundamental coweights that correspond to the simple roots in $\Delta
\smallsetminus\Delta^{1}$.
Then $h_2$ takes
the value 1 on every simple root for the Levi factor $\g^{_2}$ of $\p'$.  At the
same time, $h_2$, defined similarly,
 takes the value 1 on every simple root for $\g^{_1}$.
\bigskip

We now give an alternative more geometric proof of Theorem 3.9,
based on the following construction.

Let $\p=\l\oplus \u$ be a parabolic subalgebra in $\g$ 
with Levi factor $\l$ and the nilradical $\u$.
Let $P= L\cdot U \subset G$ denote the
corresponding 
(connected) parabolic subgroup.
Fix a principal
nilpotent $e_1 \in \l$.
Write $\z_\u(e_1)$ and $Z_{_P}(e_1)$ for the
centralizers of $e_1$ in $\u$ and $P$, respectively.

\proclaim{Proposition 3.12} $\op{(\bold i)}$ Assume that
$e_2\in\z_\u(e_1)\subset\u$ is a Richardson element for $\p$.
Set $\e=(e_1,  e_2)$. Then the following conditions
$\op{(a)}$ and $\op{(b)}$
are equivalent:

$\hphantom{x}\,\,\quad\qquad\op{(a)}$
 The orbit $\Ad Z_P(e_1)e_2$ is Zariski open and dense in $\z_\u(e_1)$.

$\qquad\qquad\op{(b)}\quad\e=(e_1,  e_2)$ is a principal nilpotent pair in $\g$.

$\op{{\bold {(ii)}}}$ 
Any two pairs $(e_1,  e_2)$ and $(e_1,  e'_2)$ that satisfy the
equivalent
conditions $\op{(a)-(b)}$ above are conjugate
by the group $Z_P(e_1)$.

$\op{{\bold {(iii)}}}$ 
If $\e=(e_1,  e_2)$ is a principal nilpotent pair in $\g$
and $\p\north=\bigoplus_{q\geq 0}\,\g_{*,q},$
the associated parabolic, see (3.3),
then part $\op{(i)}$ holds for $\l:=\g^{_1}=\g_{*,0}$ and $\p:=\p\north$.
\endproclaim

{\it Proof.} Fix $e_1$ as in (i).
Observe that we have $\z_\p(e_1)=\z_\l(e_1)\oplus\z_\u(e_1)$.
Furthermore, $\dim\z_\l(e_1)=\rk \l$, since $e_1$ is
a principal nilpotent in $\l$. Hence,
$\dim\z_\p(e_1)=\rk \l + \dim\z_\u(e_1)$. 

If $e_2\in \u$ is Richardson, then one has $\z_\g(e_2)=\z_\p(e_2)$.
Therefore, if $e_2$ is a Richardson element in $\z_\u(e_1)$ then:
 $\z_\g(e_1, e_2)=\z_\p(e_1, e_2)$, and
the Richardson inequality for the pair $\e=(e_1, e_2)$ yields:
$\dim \z_\p(\e) = \dim\z_\g(\e) \geq \rk\g =\rk \l$. Thus
we find
$$
\align
\dim \z_\u(e_1) = \dim\z_\p(e_1) - \rk\l &\geq \dim \z_\p(e_1) - \dim
\z_\p(\e)\tag 3.13\\
& = \dim Z_P(e_1)/Z_P(\e)
= \dim \bigl(\Ad Z_P(e_1)\cdot e_2\bigr)\,.
\endalign
$$
We see that the equality in (3.13) holds if and only if 
$\dim \z_\g(\e) =\rk\g$. Thus, the $\Ad Z_P(e_1)$-orbit
of $e_2$ is Zariski open in $\z_\u(e_1)$ if and only if
the
pair $\e=(e_1, e_2)$ is regular. Since $\z_\u(e_1)$ is 
a vector space, hence an irreducible variety, it may
contain at most one Zariski open orbit, which is then dense
in $\z_\u(e_1)$. This proves that i(b)$\Rightarrow$ i(a),
and also part (iii) of the Proposition.

To prove  i(a) $\Rightarrow$ i(b)
assume that the $\Ad Z_P(e_1)$-orbit
of $e_2$ is Zariski open dense
in $\z_\u(e_1)$.
 By (3.13) and the discussion following it,
 we must only show that the pair
$\e=(e_1,e_2)$  satisfies the Nil-condition 1.2(ii).
To that end, consider the natural $\C^*$-action on $\g$
by homotheties, and the induced
$\C^*$-action on $\z_\u(e_1)$.
The $\C^*$-action clearly commutes with the
$\Ad Z_P(e_1)$-action, hence, preserves the open $\Ad Z_P(e_1)$-orbit
in $\z_\u(e_1)$. It follows that, for any
$t_2\in \C^*$, the pairs 
$(e_1, t_2\cdot e_2)$ and $(e_1, e_2)$ are
$\Ad Z_P(e_1)$-conjugate. Now, the orbit
$\Ad L(e_1)$ is  nilpotent, hence, is a
$\C^*$-stable subvariety in $\l$. We deduce
that the $\Ad P$-diagonal orbit of $\e=(e_1, e_2)$ is a
$\C^*\times\C^*$-stable
subvariety in $\p\oplus\u$. This implies condition 1.2(ii),
and part i(b) follows.

To prove part (ii) assume $e'_2$ is another element satisfying
conditions
of part (i). Then, by i(a), the $\Ad Z_P(e_1)$-orbits of
$e_2$ and $e'_2$ are both Zariski open dense in $\z_\u(e_1)$.
Hence, these orbits coincide. Therefore,  $e'_2 \in \Ad Z_P(e_1)\cdot e_2$,
 and part (ii) 
is proved.\qed\medskip

Let $e_1$ be a principal nilpotent in a Levi
subalgebra $\l$, and
let $\p^{(\imath)}\,,\,(\imath=1,2),$ be two parabolics with nilradicals
$\u^{(\imath)}$
and the
same
Levi factor $\l$. Assume that $e^{(\imath)}_2\in
\z_{\u^{(\imath)}}(e_1)\,,\,
\imath=1,2,$ is a
Richardson
element for $\p^{(\imath)}$ satisfying conditions (a)-(b) of Proposition
3.12.
\smallskip

\noindent
{\bf Question 3.14.} 
Is it true that the principal nilpotent pairs $(e_1, e^{(1)}_2)$ and
$(e_1, e^{(2)}_2)$ are always $\Ad G$-conjugate ? 
\medskip

{\it Second proof of Theorem 3.9.} Part (iii) of Proposition 3.12
says that every  principal nilpotent pair in $\g$  arises from a
certain
parabolic subalgebra $\p\subset\g$ via the construction of
 Proposition 3.12(i). Note that the
$\Ad G$-conjugacy class of the pair $\e$ arising from
the construction does {\it not}
depend on the choice of a Levi factor $\l\subset \p$, since all
such factors are $\Ad P$-conjugate.
Moreover, part (ii) of Proposition 3.12
insures that any two principal nilpotent pairs arising from the same
parabolic $\p$ are $\Ad G$-conjugate. Hence, the result follows from
the finiteness of the number of $\Ad G$-conjugacy classes of parabolic
subalgebras in $\g$.
\qed\medskip

In the next section (see Corollary 4.13) we will prove the following

\proclaim{Theorem 3.15} The conjugacy class $\Ad G(e_2)\subset\g$  of the second
member of a principal nilpotent pair $(e_1,e_2)$ is totally determined
by the conjugacy class of $e_1$.
\endproclaim

\medskip
\bigskip

\head{4. Harmonic polynomial attached to a principal nilpotent pair.}
\endhead
\bigskip

 The main result of this section, Theorem 4.4 below, holds
for all principal nilpotent pairs with one exception.
The exception occurs for a particular
 principal nilpotent pair $\e_{\text{except}}=
(e_1, e_2)$ in the simple Lie algebra $\g$
of type $\GR{E}{7}$. 
The elements  $e_1$ and $e_2$ in the pair both belong to the
same nilpotent orbit $\co_{\text{except}} \subset \g$, the
regular nilpotent conjugacy class in the Levi subalgebra of type 
$\GR{A}{4} +
\GR{A}{1}$ (according to the classification of \S8,
there is only one, up to conjugacy,
 principal nilpotent pair in $\g=\GR{E}{7}$ with this property).
The nilpotent orbit $\co_{\text{except}}$ is known to have
an  "exceptional" behavior in other respects as well, see Remark 4.12
below. 
\medskip

\noindent
{\bf Definition 4.1.} A commuting
 nilpotent pair $\e=
(e_1, e_2) \in \ZZ$, resp. a nilpotent orbit $\co$,
in a semisimple Lie algebra
$\g$ is said to be {\it non-exceptional} if none of the
components of $\e$, resp. $\co$,
 corresponding to the simple factors of $\g$ of type
$\GR{E}{7}$ are conjugate to the pair $\e_{\text{except}}$,
resp. are equal to $\co_{\text{except}}$.
\medskip

Given a vector space $V$, we
write $\SS\!V=\bigoplus_i\,\SS^i\!V$ for the Symmetric algebra of
$V$,
 identified with $\C[V^*]$, the polynomial algebra on the dual
space. Abusing the language we will often refer to  elements of
$\SS\!V$ as "polynomials". We also have the completed algebra
$\, {\widehat {\SS\!V}}= \prod_{i\ge 0} \,\SS^i\!V \;\simeq$
$\C[[V^*]]\,$ of formal power series.

We keep the notations of the previous sections,
in particular, we have fixed a principal nilpotent pair
$\e$, and an associated semisimple pair $\bh=(h_1,h_2)$.
Let $\h=\z(\bh)$ be the corresponding Cartan subalgebra, and
$W$ the Weyl group. We write $\triv$ for the trivial 1-dimensional
representation of $W$, and $\sign$ for the {\it sign}-representation.
For any simple $W$-module $E$,
the space $E\; \bigotimes\; (E^*\otimes\,\sign)$
contains a unique  1-dimensional
subspace, $\sign\cdot\operatorname{Id}_E \subset \sign\otimes\Hom(E,E) =$
$E\; \bigotimes\; (E^*\otimes\,\sign)\,,$
that transforms as the sign-representation under the
diagonal $W$-action. Recall further that any simple $W$-module
is defined over ${\Bbb Q}$, hence is
self-dual, $E^*\simeq E$. Thus, there is also a distinguished line
in $E\; \bigotimes\; (E\otimes\,\sign)$ that will be referred to
 as  the {\it sign-subspace}.

To each element $\x=(x_1,x_2)\in\h\oplus\h$ we associate the following
"doubled" analogue of a familiar Weyl character type alternating expression:
$$
\align
\P_\x:= \sum_{w\in W}\;& \eps(w)\cdot
e^{w(\x)}\,
=\,
\sum_{w\in W}\; \eps(w)\cdot w(e^{x_1}\otimes e^{x_2})\tag 4.2\\
&=\,\sum_{i,j\geq 0}\, \frac{1}{i!\cdot j!}\cdot
\left(\sum_{w\in W}\; \eps(w)\cdot w(x_1^i\otimes x_2^j)\right)\,
\;\in\;\prod_{i,j\ge 0} \,\SS^{i+j}(\h\oplus \h)\;,
\endalign
$$
where the group $W$ acts diagonally on 
$\SS\!\h\otimes \SS\!\h = \SS(\h\oplus \h)$. Note that
the element $\P_\x$ is non-zero whenever $\x$ is regular, i.e.,
when all the points $w(\x)\,,\, w\in W$ are distinct.
For each non-negative integers $(d_1,d_2)$, we
consider the  polynomial: $\;
\Delta_\x(d_1,d_2) :=$
$(d_1!\cdot d_2!)^{-1}\cdot
\sum_{w\in W}\; \eps(w)\cdot w(x_1^{d_1}\otimes x_2^{d_2})\,,$
the $(d_1,d_2)$-bihomogeneous component of the
Taylor expansion of (4.2).
\smallskip

Recall 
 the  Levi suablgeras $\g^{_1}=\z_\g(h_2)$ and $\g^{_2}=\z_\g(h_1)$.
Write $\c(\g^{_i})\subset \h\,,\, (i=1,2),$ for the center of $\g^{_i}$,
and set $\c(\g^{_i})^{\text{reg}}:= \{ h\in \c(\g^{_i})\;|\;
\z_\g(h)=\g^{_i}\}\,.$ We define:
$$ \align
&\c(\g^{_1})^\circ:= \{ h \in \h\;|\;
\alpha(h)= 0,\,\forall \alpha \in R^1\;,\;
\beta(h)\neq 0,\,\forall \beta \in R^2\}\\
&\c(\g^{_2})^\circ:= \{ h \in \h\;|\;
\beta(h)= 0,\,\forall \beta \in R^2\;,\;
\alpha(h)\neq 0,\,\forall \alpha \in R^1\}\,.
\endalign
$$
We have: $\c(\g^{_i})^{\text{reg}}\,
\subset\,\c(\g^{_i})^\circ \,
\subset\,\c(\g^{_i}).$ Put $\bc := \{\x=(x_1,x_2)\;|\;
x_1 \in \c(\g^{_2})^\circ\,,\,x_2 \in \c(\g^{_1})^\circ\}\,,$
(note the {\it flip} involved in the definition).
 It is clear from the definitions that:
 $h_1\in\c(\g^{_2})^{\text{reg}}$, and $h_2\in\c(\g^{_1})^{\text{reg}}$,
hence $\bh\in\bc$.

We introduce two integers:
$\dis
\d_i\, :=\, \sharp\, R^{_i}_+\, =\, \dim \g^{_i}_{_+}\,,$
$(i=1,2)$, see 2.2(i).
By Corollary 6.10 (of section 6  below) one can rewrite these numbers
 in terms of
bi-exponents:
$$
\d_1 =\sum_{i,j}\,i\cdot\dim\z_{i,j}(\e) =\!
\sum_{(p,q)\in {\text{Exp}}_\e(\g)}\!\!p
\quad,\quad
\d_2=\sum_{i,j}\,j\cdot\dim\z_{i,j}(\e)
=\!\sum_{(p,q)\in {\text{Exp}}_\e(\g)}\!\!q\;.
$$

The following result will be proved shortly.

\proclaim{Lemma 4.3} For any $\x= (x_1, x_2)\in \bc$,  we have:

$\op{(i)} \quad\Delta_\x(d_1,d_2) = 0$ whenever $d_1<\d_1$ or
$d_2<\d_2$.

$\op{(i)}$ The 
polynomial $\Delta_\x(\d_1,\d_2)$ is proportional to
$\Delta_\bh(\d_1,\d_2)$, hence is independent, up to a
constant factor,  of the choice of
$\x\in \bc$.
\endproclaim

We set
$$\P_\e :=\Delta_\bh(\d_1,\d_2)= \frac{1}{\d_1!\cdot \d_2!} \cdot
 \sum_{w\in W}\; \eps(w)\cdot w(h_1^{\d_1}\otimes h_2^{\d_2}) 
\;\,\in\; \SS^{\d_1}\!\h\otimes \SS^{\d_2}\!\h .
$$
According to the lemma, this polynomial of bi-degree
$(\d_1,\d_2)$
is (if non-zero) the first non-vanishing term
in the Taylor expansion of $\Delta_\bh$
(the the notation $\P_\e$ is legitimate since $\e$ is not in
$\h\oplus\h$,
so that definition (4.2) doesn't apply).

We recall a few standard results and notation concerning  Weyl groups.
First, associated to the
Levi subalgebras $\g^{_1}$ and $\g^{_2}$, respectively,
there are polynomials:
$$ \pi_1 := \prod_{\alpha\in R^{_1}_+}\alpha^\vee\;\in\;
\SS^{\d_1}\!\h,
\quad\text{and}\quad
\pi_2 := \prod_{\alpha\in R^{_2}_+}\alpha^\vee\;\in\;
\SS^{\d_2}\!\h\,.
$$

Let $\C[W]$ denote the group algebra of $W$. Write
$E_i:=\C[W]\!\cdot\!\pi_i \subset \SS^{\d_i}\!\h$,
for the $W$-submodule generated by
${\pi_i}_{\,\{i=1,2\}}$.
It was shown in [M1]
that $E_i$ is a simple $W$-module that occurs in
$\SS^{\d_i}\!\h$ with multiplicity 1, and does not occur in
$\SS^d\!\h$, for any $d<\d_i$.
Furthermore, $\pi_i$ is a $W$-harmonic
polynomial, cf. e.g., [CG, \S6.3].

\medskip

One of the main results of this paper is the following

\proclaim{Theorem 4.4} If $\e$ is a non-exceptional (Def. 4.1)
principal nilpotent pair, then

$\op{(i)}$ We have:$\enspace E_2\simeq E_1\otimes\sign\,.$

$\op{(ii)}\quad \P_\e$ is a non-zero $W$-harmonic
polynomial  with respect to the
diagonal \par
$W$-action on $\h\oplus\h$;

$\op{(iii)}\enspace$ The $W\times W$-submodule in $\SS^{\d_1}\!\h\otimes \SS^{\d_2}\!\h$
generated by $\P_\e$ equals $E_1\otimes E_2$, \par
 and $\P_\e$ is a generator of the one-dimensional
$\boldkey{s}\boldkey{i}\boldkey{g}\boldkey{n}$-subspace
in  $E_1\otimes E_2$.
\endproclaim

\smallskip

{\it Proof of Lemma 4.3.} Recall the decomposition
$\g^{_i}=
\c(\g^{_i})\oplus\s^{_i}$. Fix 
 $s\in \h^{_1}=\h\cap\s^{_1}$ such that $\alpha(s) \neq 0$,
for all $\alpha\in R^{_1}$. It is well-known that there
is a non-zero constant  {\it const}$_s$ such that
the following identity holds:
$$
\sum_{w\in W^{1}}\,\eps(w)\cdot w(s^{d})\enspace=\enspace
\cases \kappa_s\cdot
\prod_{\alpha\in R^{_1}_+}\alpha^\vee\; =\kappa_s\cdot
\pi_1\quad\text{if}\quad
d=\d_1\\
\quad 0\qquad\quad\text{if}\quad d<\d_1
\endcases
\tag 4.5
$$
Moreover, $\kappa_s=\d_1!\,\,,$
if $s$ equals the half-sum of  positive coroots of $\s^{_1}$.
We claim first
that the expression on the left of (4.5) remains unaffected
if $s$  is replaced there by $x\in s+ \c(\g^{_1})$. Indeed, 
for any $x=s +c$ where $c\in \c(\g^{_1})$,
we find: $\sum_{w\in W^{1}}\,\eps(w)\cdot w(x)^d= $
$\sum_{w\in W^{1}}\,\eps(w)\cdot w(s+c)^d=$ 
$\sum_{w\in W^{1}}\,\eps(w)\cdot w(s)^d + \ldots\,,$
where dot-terms belong
to the components $\SS^k(\c(\g^{_1}))\otimes\SS^{d-k}(\s^{_1})$
with $k>0$. But the sign-representation 
of $W^{1}$ does not occur in $\SS^d(\s^{_1})$
for any $d<\d_1$. Hence, all the dot-terms vanish,
and the claim follows. Observe, that we have proved in particular
that, for any $s\in \c(\g^{_2})^\circ\subset \h$, the identity
(4.5) holds and, moreover, $\kappa_s\neq 0$.

Now fix $\x=(x_1,x_2)\in \bc$ and
some $d_1, d_2\geq 0$, where  $d_1\leq\d_1$.
We have: $x_2\in \c(\g^{_1})$, hence
the element $x_2^{d_2}\in\SS^{d_2}\!\h$ is fixed
by the $W^{1}$-action. Therefore, using (4.5) and the claim proved
in the preceding paragraph  for $s=x_1\in \c(\g^{_2})^\circ$
we get
$$
\sum_{w\in W^{1}}\,\eps(w)\cdot w(x_1^{d_1}\otimes x_2^{d_2})=
\cases \kappa_{x_1}\cdot
\pi_1\otimes x_2^{d_2}\;\in\,\SS^{d_1}(\s^{_1}) \otimes
\SS^{d_2}\!\h
\quad\text{if}\quad
d_1=\d_1\\
\quad 0\qquad\quad\text{if}\quad d_1<\d_1
\endcases
\tag 4.6
$$
where $\kappa_{x_1}\neq 0$.
The element $\Delta_\x(d_1,d_2)$ is clearly obtained by alternating the
expression on the LHS of (4.6) with respect to the diagonal
$W$-action on $\SS^{d_1}\!\h\otimes\SS^{d_2}\!\h$. 
By (4.6), this expression vanishes whenever $d_1<\d_1$, and the lemma
follows.
\qed
\medskip

\proclaim{Lemma 4.7} If the polynomial $\Delta_\e$ is non-zero,
then all other claims of Theorem 4.4 hold.
\endproclaim

{\it Proof.} If $\pe$ is nonzero, it is really 
the first non-vanishing term of $\Delta_\bh$, by  Lemma 4.3.
But the first non-vanishing term of
any linear
combination of exponents like in (4.2),
i.e., a combination of $W$-conjugate exponents, is known to
be $W$-harmonic,
see e.g. [CG, Prop.6.4.4], hence part (ii) of the theorem.

To prove part (iii), we observe, as in the proof of Lemma 4.3,
that $\pe$ is obtained by alternating the
expression on the LHS of (4.6) with respect to the diagonal
$W$-action. 
 The RHS of (4.6) 
shows that the result of such an alternating procedure 
clearly belongs (for $d_i=\d_i\,,\, i=1,2$)
 to the subspace $E_1\otimes \SS^{\d_2}\!\h
\,\subset$ $\SS^{\d_1}\!\h\otimes\SS^{\d_2}\!\h$. Similar arguments,
with the roles of $h_1$ and $h_2$ reversed, imply that
$\pe \in \SS^{\d_1}\!\h\otimes E_2$. Now, part (iii) follows from the equality:
$$\left(E_1\otimes \SS^{\d_2}\!\h\right)\;\cap\;
\left(\SS^{\d_1}\!\h\otimes E_2\right)\;
=\;E_1 \otimes E_2\;.
$$

Finally, 
using selfduality of all $W$-modules we get:
$$\Hom_W(E_2,\,E_1\otimes\sign)\, =\,
\bigl(E_2\;\otimes\;(E_1\otimes\sign)\bigr)^W \,=\,
\Hom_W(\sign,\, E_2\otimes E_1) \,,$$
and $\P_\e$ is a non-zero element in the RHS.
Hence, the LHS is non-zero, and part (i) of the theorem follows
from Schur lemma.\qed\medskip

Thus, most of the remaining part of this section is devoted to the proof of 

\proclaim{Proposition 4.8} If $\e$ is non-exceptional (Definition 4.1)
then, for any $\x\in\bc$, we have:$\;$ $\Delta_\x(\d_1,\d_2)\neq 0$, in
particular, $\Delta_{\e}\neq 0$. On the contrary,
 $\Delta_{\e_{\text{except}}}=0$.
\endproclaim

Let ${}^{\m\!}\B$ denote the Flag variety of all Borel 
subalgebras in a Lie algebra $\m$.
We begin the proof 
by reinterpreting the representations $E_i$ geometrically,
by means of Springer theory, cf. e.g., [CG, ch.3]. 

 Given a nilpotent
$x\in \g$, set
$C(x):= Z_G(x)/Z^\circ_G(x)$,
and let ${}^{\g\!}\B_x\subset {}^{\g\!}\B$ denote the subvariety of the
Borel 
subalgebras in $\g$ containing $x$. By  Springer theory,
 there is a natural $W$-action on each homology
group $H_i({}^{\g\!}\B_x)= H_i({}^{\g\!}\B_x, \C)$. Moreover, the
$W$-action commutes with the natural $C(x)$-action,
and one has a $W$-module isomorphism:
$$ E_i=H_{\text{top}}({}^{\g\!}\B_{e_i})^{C(e_i)}\,,
\quad\text{where}\quad \text{"top"}:=\dim_{_{\Bbb R}}\B_{e_i} =
2\d_i\;\;,\; i=1,2. \tag 4.9$$
In general, given any nilpotent $x\in \g$,
we will write $H(\B_x)$ for the representation of $W$ 
in the subspace of $C(x)$-invariants of the top homology of
${}^{\g\!}\B_x$.

Next, recall the concept
 of {\it induction} of nilpotent orbits.
Given a Levi subalgebra $\el\subset \g$, we write
$\Ind_\el^\g(0)$ for the nilpotent orbit in $\g$ that intersects
the nilradical of some (hence, every) parabolic subalgebra of $\g$
with Levi factor $\el$ by an open dense subset. 
We refer to [LS] for results, and [Spa] and [Ke] (resp. [El2]) for more
concrete information concerning induced
orbits in classical (resp. exceptional) Lie algebras.
The results of [El2] were announced without proofs in [Spa].

It will be convenient for us to introduce the following
\medskip

\noindent
{\bf Definition 4.10.} Two nilpotent orbits $\co_1$ and $\co_2$ in $\g$
are said to be
{\it reciprocal}, if there exist Levi subalgebras $\g^{_1}$ and $\g^{_2}$
such that
the following conditions  hold:

(i)$\quad\co_1 =\Ind_{\g^{_2}}^\g(0)\quad$
 and $\quad\co_2 =\Ind_{\g^{_1}}^\g(0)$. \par
(ii)$\quad\g^{_i}\cap\co_i$ 
is the regular nilpotent orbit in $\g^{_i}$, for $i=1,2$.
\medskip

This notion is relevant because of  Proposition 3.4 saying that,
for a principal nilpotent pair $\e=(e_1,e_2)$, the corresponding
nilpotent orbits $\co_1=\Ad G(e_1)$ and $\co_2=\Ad G(e_2)$ are
reciprocal orbits in $\g$. We note that, for $\g$ simple of type $E_7$,
 the exceptional orbit,
$\co_{\text{except}}$, is reciprocal to itself, see
 beginning of this section.

A crucial ingredient in our proof of the
non-vanishing of $\pe$ is the following 

\proclaim{Proposition 4.11} If $\co_1$ and   $\co_2$ are reciprocal
non-exceptional (Definition 4.1)
nilpotent orbits in $\g$, then for any $x_1\in \co_1\,,\,
x_2\in \co_2$, there is a $W$-module isomorphism:
 $H(\B_{x_2}) \simeq \sign \otimes H(\B_{x_1}).$
\endproclaim

This result is false if $\co_1=\co_2=\co_{\text{except}}$.
\medskip

\noindent
{\bf Remark 4.12.} Lusztig has introduced (see [L2] or [Ca]) the notion of
a {\it special} nilpotent orbit in $\g$. 
He introduced also the notion of
a {\it special} irreducible representation of the Weyl group [L3], and
proved  that the nilpotent orbit $\Ad G(x)$ is special
if and only if the Springer representation $H(\B_{x})$ is
special. The orbit $\co_{\text{except}}$ is special.
The corresponding
representation $H(\B_{x})\,,\,x\in \co_{\text{except}},$
 is among very few of the special
representations $E$, such that the representation $\sign \otimes E$
is not special, see [L3].
Furthermore, Spaltenstein has defined (on a case-by-case
basis) an involution $\sigma$ on the set of all special
nilpotent orbits in $\g$, see [Spa]. It is known, see e.g., [Ca, pp.373-374,
389] that in all but
 two exceptional cases in types $E_7 ,E_8$ the following holds. 
If $\co_1$ and   $\co_2$ are special nilpotent orbits
such that $\co_2=\sigma(\co_1)$, then
$H(\B_{x_2}) \simeq \sign \otimes H(\B_{x_1})$, for any
$x_i \in \co_i\,,\, i=1,2$. This property does not hold
for $\co_1=\co_{\text{except}}$, in which case we have:
$\sigma(\co_{\text{except}})=\co_{\text{except}}$,
but $ H(\B_{x_1}) \not\simeq \sign \otimes H(\B_{x_1})$,
see [Ca].
\medskip

{\it Sketch of first proof of Proposition 4.11.} For each simple Lie
algebra $\g$, there are available tables of all
induced orbits of the form $\Ind_\el^\g(0)$, see [El1] for exceptional
Lie algebras and  [Spa] (or [CM]) for classical Lie algebras.
From this, it is a streightforward matter to derive a
list all reciprocal pairs. We do {\it not} know at the moment
which of the reciprocal pairs arise from principal nilpotent pairs.
However, using explicit tables for the Springer correspondence,
see [AL], [Sh] and references in [L2], one verifies case-by-case that
the result of the proposition holds for  every reciprocal pair
in any simple Lie algebra, but
the exceptional pair for $\g$ of type $E_7$ (the special nilpotent
orbit in $\g$ of type $E_8$ with an "exceptional"
behavior in the sense of Remark 4.12 does not give rise
to any reciprocal pair).\qed
\medskip

\medskip 
Recall that every
induced orbit of the form $\Ind_\el^\g(0)$ is known to be special.
In particular, each member of any reciprocal pair $(\co_1, \co_2)$
is a special orbit. We have

\proclaim{Lemma 4.13} If $(\co_1, \co_2)$ is a reciprocal pair, then
$\co_2=\sigma(\co_1)$, where $\sigma$ is the 
Spaltenstein involution
on the set of special orbits, see Remark 4.12.
\endproclaim 

{\it Proof.} Let $\el\subset \g$ be a proper Levi subalgebra. It was shown by
Barbasch-Vogan [BV], see also [CM, Theorem 8.3.1], that 
if $\co \subset \g$ is a special nilpotent orbit such that
$\co_\el:=\co \cap\el$ is a special orbit in $\el$, then
$\sigma(\co) = \Ind_\el^\g\bigr(\sigma(\co_\el)\bigl)$.
It follows from Definition 4.10(ii) that, for
a reciprocal pair $(\co_1, \co_2)$, the set
$\co_{\g^{_1}}:=\co_2\cap \g^{_1}$ is a principal nilpotent orbit in $\g^{_1}$.
Hence, $\sigma(\co_{\g^{_1}})= \{0\}$ is the zero-orbit
in $\g^{_1}$.
Therefore, by the above we get 
$\dis \sigma(\co_1) = \Ind_{\g^{_1}}^\g(0) = \co_2\,.$\qed

\medskip
{\it Second proof of Proposition 4.11.}
If $(\co_1, \co_2)$ is a reciprocal pair, then
$\co_2=\sigma(\co_1)$ by the Lemma. The claim now follows from
the result mentioned at the end of Remark 4.12, 
saying that $\co_2=\sigma(\co_1)$ implies
$H(\B_{x_2}) \simeq \sign \otimes H(\B_{x_1})$,
except for the two
cases in types $E_7,E_8$, which are handled separately.
\qed

\medskip

{\it  Proof of Proposition 4.8.} By Proposition 4.11, we know that
the vector space $E_1\otimes E_2 \subset
\SS^{\d_1}\!\h\bigotimes\SS^{\d_2}\!\h$
contains a distinguished  1-dimensional sign-subspace. Let
$\Delta = \sum\;P'_k\otimes P''_k\,,$
$\,P'_k\in \SS^{\d_1}\!\h,\, P''_k\in \SS^{\d_2}\!\h,
$ be a non-zero element of this
sign-subspace. We will show that, up to a non-zero constant
factor, we must have: $\Delta = \P_\x(\d_1, \d_2)$, hence 
$\P_\x(\d_1, \d_2)\neq 0$.

To this end, observe first that $\Delta$ is clearly
obtained by alternating the expression $\sum\;P'_k\otimes P''_k$
with respect to the $W$-diagonal action. Recall further, that
the module $E_2$ is $\C[W]$-generated by $\pi_2$. Hence,
$P''_k= u_k(\pi_2)\,,\, u_k\in \C[W],$ and it is then easy to
see that alternating $\sum\;P'_k\otimes P''_k$ gives the
same result as alternating an expression of the form
$Q\otimes\pi_2,$ where $Q\in \SS^{\d_1}\!\h.$ 
Using 
the identity $\pi_\imath= \kappa_{x_\imath}
\cdot \sum_{w_\imath
\in W^{\imath}}\;\eps(w_\imath)\cdot
w(x_\imath^{\d_\imath})$,
see
(4.5), we can further rewrite the latter alternating
expression as follows:
$$
\align
\sum_{w\in W}\;
\eps(w)\cdot w(Q\otimes\pi_2)& =\kappa_{x_2}
\cdot
\sum_{w\in W}\;
\eps(w)\cdot w(R\otimes x_2^{\d_2})\\
& = \kappa_{x_2}
\cdot\frac{1}{\sharp W^{_1}}
\sum_{w\in W}\;
\eps(w)\cdot w\Bigl(\left(\sum_{w_1\in W^{_1}}\;\eps(w_1)\cdot w_1(R)\right)
\otimes x_2^{\d_2}\Bigr)\;.
\endalign
$$
But since $\pi_1$ is the only $W^1$-skew-invariant in $\SS^{\d_1}\!\h$,
for any $R\in \SS^{\d_1}\!\h,$ the polynomial
$\sum_{w_1\in W^{_1}}\;$
$\eps(w_1)\cdot w_1(R)$ must either vanish,
or else be proportional to $\pi_1$. It cannot vanish
in our case, since $\Delta\neq 0$. We conclude that
$\Delta = \sum_{w\in W}\;
\eps(w)\cdot w(R\otimes x_2^{\d_2})$ equals 
(up to a non-zero constant factor, cf. (4.6)) to
$\sum_{w\in W}\;
\eps(w)\cdot w(\pi_1\otimes x_2^{\d_2})=\P_\x(\d_1, \d_2)$.\qed\bigskip

It is instructive to give a conceptual proof of Proposition 4.11,
at least in a special case. To this end, fix two Levi subalgebras
$\g^{_1}, \g^{_2}$ in $\g$ with (abstract) Weyl groups
$W^{1},\,W^{2}\subset W$, respectively. 
Let $(\co_1, \co_2)$
be a reciprocal pair of nilpotent orbits in $\g$, as in Definition 4.10.
Choose $e_i \in \co_i\cap \g^{_i}$, 
a principal nilpotent in $\g^{_i}\,,\,i=1,2.$ By definition, there are
parabolic subgroups $\p\north, \p\east\subset \g$ with Levi factors
$\g^{_1}$
and  $\g^{_2}$, respectively, and such that
the element $e_2$ is Richardson for $\p\north$, while
the element $e_1$ is Richardson for $\p\east$. In general,
given a parabolic $\p$ and  any nilpotent $x\in\g$, let
$\PP_x$ denote the variety of all parabolics of type $\p$
that contain $x$ in their nilradical. In particular, we consider
the sets $\PP\north_{e_2}$ and $\PP\east_{e_1}$. These sets are finite,
since the elements in the corresponding subscript are Richardson.

\proclaim{Proposition 4.14} $\op{(i)}$ There are natural
$W$-module isomorphisms:

\centerline{$\dis
\Ind_{_{W^{_i}}}^{^W}\!\triv\;\simeq\;
H_*({}^{\g\!}\B_{e_i})\quad,\quad i=1,2$.}
\vskip 2pt

$\op{(ii)}$ There are natural  vector space isomorphisms:
\vskip 2pt
\centerline{$\dis H_0(\PP\north_{e_2}) \;\simeq\;
\Hom_W\bigl(\Ind_{_{W^{_1}}}^{^W}\!\sign\,,\,
\Ind_{_{W^{_2}}}^{^W}\!\triv\bigr)\;\simeq\;
H^0(\PP\east_{e_1})\;.$}\vskip 3pt

$\op{(iii)}$  The following conditions are
equivalent:
$$ Z_G(e_2) \subset P\north\;\enspace\Longleftrightarrow\;\enspace
Z_G(e_1) \subset P\east\;\enspace\Longleftrightarrow\;\enspace
\sharp \PP\north_{e_2}\,=\, 1\;\enspace\Longleftrightarrow\;\enspace
\sharp \PP\east_{e_1}\,=\, 1\,.$$
\endproclaim

{\it Proof.}  By  Springer theory, for
any nilpotent $x\in{\g^{_2}}$, there is
a $W^{_2}$-action on $H_*({}^{{\g^{_2}}\!}\B_x)$, 
and a result of Borho-MacPherson [BM, 3.4] 
says that there is a natural $W$-module isomorphism:
$$\Ind_{_{W^{_2}}}^{^W}\!H_*({}^{{\g^{_2}}\!}\B_x)\;
\simeq\;H_*({}^{\g\!}\B_x)\,.
$$

In the special case of the principal
nilpotent $x=e_2$ in ${\g^{_2}}$, the set ${}^{{\g^{_2}}\!}\B_{e_2}$
consists of a single point,
hence, $H_*({}^{{\g^{_2}}\!}\B_{e_2}) = \triv$.
Applying the Borho-MacPherson  isomorphism in this case yields the 
isomorphism of part
(i) of the Proposition.

Next, given a $W$-module $M$, write
$M^{\sign(W^{_1})} := \{m\in M\;|\; w_1\cdot m = \eps(w_1)\cd m\,,$
$\,
\forall w_1\in W^{1}\}\,,$ the 
$\sign\big|_{W^{_1}}$-isotypic component of
$M\big|_{W^{_1}}$.
 By Frobenius reciprocity one has: $M^{\sign(W^{_1})}$
$\,\simeq\,Hom_W\bigl(\Ind_{_{W^{_1}}}^{^W}\!\sign\,,\,M\bigr)\,.$
Borho-MacPherson have proved that, for any $i\geq 0$
and any nilpotent $x\in\g$, 
the $\sign(W^{_1})$-isotypic component of 
 $M=H_i({}^{\g\!}\B_x)$ is given by the formula:
$\,\dis H_i({}^{\g\!}\B_x)^{\sign(W^{_1})} \;\simeq\;
 H_{i-2\d_1}(\PP\north_x)\,,$
$\d_1=\dim_{_\C}{}^{{\g^{_2}}\!}\B,$ see [BM, 3.4].
From this formula and part (i) we find:
 $$
\Hom_W\bigl(\Ind_{_{W^{_1}}}^{^W}\!\sign\,,\,
\Ind_{_{W^{_2}}}^{^W}\!\triv\bigr) 
\;\simeq\;
\bigl(\Ind_{_{W^{_2}}}^{^W}\!\triv\bigr)^{\sign(W^{_1})}
\;\simeq\;\bigl(H_*({}^{\g\!}\B_{e_2})\bigr)^{\sign(W^{_1})}
\;\simeq\;H_*(\PP\north_{e_2})\,.
$$
Since the set $\PP\north_{e_2}$ is finite, its homology
 reduces to 
$H_0(\PP\north_{e_2})$, and the first isomorphism
of part (ii)  follows. 

Further, using that
 $\Ind_{_{W^{_i}}}^{^W}\!\sign\; \simeq\;
\sign \bigotimes (\Ind_{_{W^{_i}}}^{^W}\!\triv)\,$ 
we obtain a chain of 
 canonical isomorphisms:
$$
\align
\Hom_W\bigl(\Ind_{_{W^{_1}}}^{^W}\!\sign\,,\,
\Ind_{_{W^{_2}}}^{^W}\!\triv\bigr) 
\;&\simeq\;\Hom_W\bigl(\sign \otimes (\Ind_{_{W^{_1}}}^{^W}\!\sign)\,,\,
\sign \otimes (\Ind_{_{W^{_2}}}^{^W}\!\triv)\bigr)\\ 
&\simeq\;\Hom_W\bigl(\Ind_{_{W^{_1}}}^{^W}\!\triv\,,\,
\Ind_{_{W^{_2}}}^{^W}\!\sign\bigr)\tag 4.15\\
&\simeq\;\Bigl(
\Hom_W\bigl(\Ind_{_{W^{_2}}}^{^W}\!\sign\,,\,
\Ind_{_{W^{_1}}}^{^W}\!\triv\bigr)\Bigr)^*
\;.
\endalign
$$
The $\Hom$-space in the last line is isomorphic,
by an analogue for $e_1$
of the first isomorphim in (ii), to
$H_0(\PP\east_{e_1}).$
Therefore we get $\Bigl(
\Hom_W\bigl(\Ind_{_{W^{_2}}}^{^W}\!\sign\,,\,
\Ind_{_{W^{_1}}}^{^W}\!\triv\bigr)\Bigr)^* \simeq
\bigl(H_0(\PP\east_{e_1})\bigr)^*$ 
$\simeq H^0(\PP\east_{e_1}),$ and
part (ii) follows.

To prove (iii) we recall that if $x$ is Richardson for a parabolic
$\p$, then the group $Z_G(x)$ acts transitively on the set of
parabolics of type $\p$ that contain $x$ in their nilradical, see [Ca].
This fact, combined with the isomorphism $H_0(\PP\north_{e_2}) \simeq
H^0(\PP\east_{e_1})$ of part (ii) yields the equivalences in (iii).
\qed
\medskip

The following result gives a conceptual proof of Proposition 4.11
in a special case. Note that the characterisation of the simple
module $E_2=H(\B_{e_2})$ arising from the last sentence of
the corollary below
provides
a natural generalisation to arbitrary Weyl groups of the 
Young's
classical construction of irreducible representations of the Symmetric
group. 

\proclaim{Corollary 4.16} Assume
the equivalent conditions of Proposition 4.14(iii) hold, e.g., 
$\g=\sln$. Then
we have: $H(\B_{e_2})\simeq H(\B_{e_1})\otimes\sign\,.$ This $W$-module is
 the only common irreducuble
constituent of $\Ind_{_{W^{_2}}}^{^W}\!\triv$
and $\Ind_{_{W^{_1}}}^{^W}\!\sign$; it occurs with multiplicity one
in either of these modules.
\endproclaim

{\it Proof.} 
By Proposition 4.14(ii) we have: 
$\dim \Hom_W\bigl(\Ind_{_{W^{_1}}}^{^W}\!\sign\,,\,
\Ind_{_{W^{_2}}}^{^W}\!\triv\bigr) \;=\;\sharp \PP_{e_2}$
$\,=\, 1\,.$ 
The above formula implies
that $\Ind_{_{W^{_2}}}^{^W}\!\triv$
and $\Ind_{_{W^{_1}}}^{^W}\!\sign$ have 
only one  irreducuble
constituent in common. Moreover, 
since $\d_1=\dim_{_\C}({}^{{\g^{_2}}\!}\B)
=\dim_{_\C}\B_{e_2}$, the Borho-MacPherson formula yields:
$H(\B_{e_2})^{\sign(W^{_1})} = H_{2\d_1}(\B_{e_2})^{\sign(W^{_1})} \simeq$
$H_0(\PP\north_{e_2}) \neq 0\,.$ Thus, $H(\B_{e_2})$ is 
this common irreducuble
constituent.

To complete the proof  we note that (4.15) implies  that
$\,\Hom_W\bigl(\Ind_{_{W^{_2}}}^{^W}\!\sign\,,\,
\Ind_{_{W^{_1}}}^{^W}\!\triv\bigr)$
is also 1-dimensional. Hence, the two induced
representations involved in the $\Hom$
also have only one irreducuble
constituent in common. By the symmetry, this irreducuble
constituent is $H(\B_{e_1})$. Thus, we must have
$H(\B_{e_1})\otimes\sign \simeq H(\B_{e_2})$, and the proof is complete.
\qed\medskip

In the special case considered above one has the following
strengthening of Theorem 4.4 with the element $\P_\e$
being replaced by $\P_\bh$. Moreover, during the
proof below, we will effectively
identify the element $\P_\bh$ with an intertwiner
analogous to the  classical "{\it Young symmetriser}" in the case of
Symmetric groups.

\proclaim{Proposition 4.17}
Let $\e=(e_1,e_2)$ be a principal nilpotent pair
such that the equivalent conditions of 
Proposition 4.14(iii)
 hold. Then, the $W\times W$-module generated by $\Delta_\bh$ is
isomorphic to $E_1\otimes E_2$.
\endproclaim 

{\it  Proof.} Choose an element $\theta\in( \h^*)^{W^1}$ generic
enough so that the linear functions $\{w(\theta)\,,\, w\in W/W^{_1}\}$
separate points of the orbit $W\cd h_1 \subset \h$. We fix such
a $\theta$, and set $\F_1=\{f\in \C[[\h]]\;|\;
f=\sum_{w\in W}\; a_w\cdot e^{w(\theta)}\;,\, a_w\in\C\}\,.$  
The space $\F_1$ has an obvious $W$-module structure and,
since $\theta$ is
fixed by
$W^1\subset W$,
there is a natural $W$-module isomorphism: 
$\F_1 \simeq \Ind_{W^{_1}}^W\triv$, that sends $e^\theta$ to $\triv$.
We will view elements of $\F_1$ as holomorphic functions on $\h$
(as opposed to the formal power series).
Similarly, we define
$\F_2=\{\Psi\in \widehat{\SS\!\h}\;|\;
\Psi=\sum_{w\in W}\; a_w\cdot e^{w(h_2)}\}\,,$
so that one has: $\F_2\simeq \Ind_{W^{_2}}^W\triv$.

Observe  next that we may (and will) identify
the element $\Delta_\bh$ with a map $\delta_\bh: \F_1 \to \F_2$
given by the formula: 
$$f \mapsto \delta_\bh(f) =
\sum_{w\in W}\; \epsilon(w)\cdot  f(w(h_1))\cdot
e^{w(h_2)}\; \in\; \F_2.$$
The map $\delta_\bh$ is non-zero, due to our choice of $\theta$.
Furthermore, it is clear from the definition
that this map  gives a $W$-module morphism
$\delta_\bh: \F_1 \to \F_2\otimes\sign$. 

Recall that $\F_i \simeq \Ind_{W^{_i}}^W\triv$,
and that $E_1$ is 
 the only common irreducible
constituent of $\Ind_{_{W^{_1}}}^{^W}\!\triv$
and $\Ind_{W^{_2}}^{^W}\!\sign$, by
Corollary 4.16(ii). It follows  that 
$\delta_\bh$, viewed as an element
of $\Hom(\F_1, \F_2)=\F_1^* \bigotimes \F_2$, belongs to
$E_1\otimes E_2$.
\qed\bigskip

We conclude this section by indicating a link
between our polynomial $\P_\e$ and canonical bases.
Recall that Kazhdan and Lusztig have constructed 
a canonical basis $\{c_w,\, w\in W\}$
in the group algebra $\C[W]$, cf. [KL], with remarkable
properties. In the same paper, Kazhdan and Lusztig 
have partitioned the set $W$
into two-sided cells $\cc_{LR}$, and have further partitioned
each two-sided cell into left cells $\cc_L$. 
To any left cell $\cc_L$, Kazhdan and Lusztig have attached
a $\C[W]$-module $E(\cc_L)$ with a canonical basis formed by
the elements $\{c_w, w\in \cc_L\}$. Furthermore, Lusztig has
constructed a bijection: $\co\, \leftrightarrow\,
\cc_{LR}(\co)$, between the set of special nilpotent
orbits in $\g$ and the set of two-sided cells in $W$.
This bijection has the following property: if $E$ is the
special irreducible representation of $W$ associated
to a  special nilpotent orbit $\co$ then,
for every left cell $\cc_L \subset \cc_{LR}(\co)$, the module
$E$ occurs with multiplicity one in $E(\cc_L)$,
and does not occur in any left-cell representations arising
from any other two-sided cell.

\proclaim{Conjecture 4.20} Let $\e=(e_1,e_2)$
be a non-exceptional principal nilpotent pair, and
$\co=\Ad G(e_1)$, a (special) nilpotent orbit in $\g$. Then,
the corresponding two-sided cell representation
of $W\times W$ is irreducible, equivalently,
for every left cell $\cc_L \subset \cc_{LR}(\co)$,
the $W$-module $E(\cc_L)$ is irreducible.
\endproclaim

The conjecture trivially holds for $\g=\sln$, since
 all left cell representations for $W=S_n$ are known 
to be irreducible, [KL, Theorem 1.3]. A. Elashvili has
verified, using
the tables of left cells given in [L3] and [Ca],  that 
the conjecture  holds  for
all non-exceptional, in the sense of Definition 4.1,
  principal nilpotent pairs 
known at the moment (the Conjecture is {\it false} for the
exceptional pair).\medskip

Let $\e=(e_1,e_2)$
be a non-exceptional principal nilpotent pair.
Set $\co=\Ad G(e_1)$, so that $E_1 \simeq H(\B_{e_1})$, see (4.9),
is the correspondinding special representation. Then
$\co=\Ad G(e_1)$, a (special) orbit in $\g$.
If Conjecture 4.20  holds for $\e$ then, for
every left cell $\cc_L \subset \cc_{LR}(\co)$,
we have: $E_1 \simeq E(\cc_L)$.
This isomorphism  transfers the canonical basis 
$\{c_w, w\in \cc_L\}$ in $E(\cc_L)$ to a basis in $E_1$.
We remark that, for $\g=\sln$, the basis in $E_1$ arising in this way
is independent of the
choice of the left cell, by [KL, Theorem 1.4]. It is likely that
this is true for a non-exceptional principal nilpotent pair
in an arbitrary semisimple Lie algebra.

 We fix (non-exceptional)
$\e$ and a left cell $\cc_L \subset \cc_{LR}(\co)$ as above.
Let $w_0\in W$ denote the element of maximal length.
According to [KL, Corollary 3.2], the assignment
$x\mapsto x\cdot
w_0$ gives an involution $\cc \leftrightarrow \cc\cdot w_0$
on the sets of two-sided cells, left cells, etc.

\proclaim{Proposition 4.21} Assume Conjecture 4.20  holds for $\e$.
Then 

$\op{(i)}$ There is a $W\times W$-module isomorphism:
$E_1\otimes E_2 \simeq E(\cc_L)\otimes E(\cc_L\cd
w_0)$. 

$\op{(ii)}$
The polynomial $\Delta_\e \in E_1\otimes E_2$ goes under the isomorphism
in (i) to the \par
element:
$\;\Delta_{\e} = 
\sum_{x\in \cc_L}\; \eps(x)\cdot c_x\otimes c_{x\cdot w_0}\;.$
\endproclaim

{\it Proof.} It is known that in general, for any left cell
$\cc_L$, there is a natural isomorphism:
$E(\cc_L\cdot w_0) \simeq \sign\otimes E(\cc_L)^*\,,$ 
see e.g [J2, 3.7]. Part (i) now follows from Theorem 4.4.

To prove (ii) we will use a more precise information
about the isomorphism in (i). Specifically, for each $w\in W$,
we let $\|m_{x,y}(w)\|_{x,y \in W}$ be
the matrix in the canonical basis of
the  operator $m(w): \C[W] \to \C[W]$
given by left multiplication 
by $w$, that is:
$$ w\cdot c_x = \sum_{y\in W}\; m_{x,y}(w)\cdot c_y\;.$$
The inversion formula for the Kazhdan-Lusztig polynomials,
see [KL, Theorem 3.1], implies the following identity, see
[J1, 4.7]: 
$$m_{x,y}(w^{-1}) = \eps(x\cd y\cd w)\cdot 
m_{x\cdot w_0, y\cdot w_0}\quad,\quad \forall x,y,w \in W\;.\tag 4.22$$
Alternatively, one may deduce (4.22) from [KL, Corollary 3.2]
as follows. Recall that by {\it loc.cit.}, the
$W$-graph (see [KL, pp.165, 167]) attached to the
left
cell $\cc_L\cdot w_0$ coincides with $W$-graph attached
to $\cc_L$, except that the function: $x\mapsto I_x$
on its set of vertices gets replaced by the function:
$x\mapsto S\smallsetminus I_x$, where $S$ denotes the set
of simple reflections in $W$.
In the special case where $w=s=w^{-1},$
is a simple reflection the identity (4.22) follows
readily from the previous discussion and
formula  [KL, (1.0.a)]. The latter reads:
$$m(s) :\; c_x \;\mapsto\;
\cases 
-c_x {\hskip 57mm} \text{if}\quad s\in I_x\\
\enspace\, c_x\, +\, {\dis\sum_{\{y\in X\;|\; 
s\in I_y\,,\, \{y,x\}\in Y\,\}}}\!
\mu(y,x)\cd c_y
\qquad\text{if}\quad s\not\in I_x\;.
\endcases
$$

Next, in $\C[W]$ define a new basis: $\{c^\circ_x := \eps(x)\cd c_{x\cdot
w_0}\,,\; x\in W\}$.
The identity in (4.22) means that the matrix of the
operator $m(w^{-1}) : \C[W] \to \C[W]$ in the basis
$\{c^\circ_x\}$ is equal to  $\|m_{x,y}(w)\|^{\op{T}}$,
the transpose of the matrix
of the
operator $m(w) : \C[W] \to \C[W]$ in the basis
$\{c_x\}$. This implies that the element
$\sum_{x\in \cc_L}\;  \eps(x)\cdot c_x\otimes  c_{x\cdot w_0}\;$
transforms as the sign-representation under the
$W$-diagonal action, and part (ii) of the proposition follows
from Theorem 4.4(iii).\qed

\medskip

\bigskip

\head{5. Distinguished nilpotent pairs; $\sln$-case.}
\endhead
\bigskip

Given a semisimple Lie algebra $\g$,
we call an arbitrary  pair $\e=(e_1,e_2)\in\ZZ$ satisfying
condition 1.1(Nil), but not necessarily the regularity condition
1.1(Reg),
a {\it nil-pair}.

Given a nil-pair $\e$,
choose a Cartan (= maximal diagonalizable) subalgebra 
in $\z(\e)$, and let ${\frak l}$ be the centralizer of
this subalgebra in $\g$. Thus, ${\frak l}$ is a Levi subalgebra
of $\g$ containing the pair
$\e=(e_1,e_2)$. Let $\z_{_{{\frak l}}}(\e)={\frak l}\cap \z(\e) $ be
 the centralizer of $\e$ in ${\frak l}$. 
By construction, every semisimple element of $\z_{_{{\frak l}}}(\e)$
belongs to the center of ${\frak l}$. 

We claim that $\e$ is a
nil-pair in ${\frak l}$. To prove this, we note that
Lemma 1.3 still applies to any nil-pair in $\g$, hence
to $\e$.
Therefore,  we can define a semisimple pair $\bh=(h_1,h_2)$, and
the Lie subalgebra
 $\overline{\m}:=\{x\in$
$\g\;|\; \ad x(e_i) \in {\C\cd e_i}\,,\;i=1,2\}$,
in the same way as we have done in the proof of Theorem 1.2(iv).
Let $\overline{\m}_{\bold {red}}$ be a maximal reductive subalgebra in $\overline{\m}$.
The argument in the proof of Lemma 1.3 shows that
one can choose the pair $\bh=(h_1,h_2)$ inside the center of
$\overline{\m}_{\bold {red}}$ so that one has a Lie algebra direct sum
decomposition:
$\overline{\m}_{\bold {red}}=\z(\e)_{\bold {red}} \oplus \langle h_1, h_2\rangle\,.$
With this choice of $\bh$, the elements $h_1,h_2$ clearly commute with
any Cartan subalgebra in
$\z(\e)_{\bold {red}}$, hence belong to the  Levi subalgebra ${\frak l}$
defined earlier. It follows that $\e$ is a
nil-pair in ${\frak l}$, and the claim is proved.
Note further that since all Cartan subalgebras in
$\z(\e)$ are conjugate to each other, all Levi subalgebras
arising from our construction are conjugate.

In general, we call a nil-pair $\e$  in a
 reductive Lie algebra ${\frak l}$ {\it pre-distinguished}  if
every semisimple element of $\z_{{\frak l}}(\e)$
belongs to the center of ${\frak l}$. Proposition 1.9  implies that any
principal nilpotent pair is pre-distinguished. 

Thus, to each nil-pair $\e$  in $\g$ we have associated a
Levi subalgebra ${\frak l}\subset\g$ such that $\e$ is a
pre-distinguished pair in ${\frak l}$. Observe further
that the argument used
in the
proof of Theorem 1.2(iv) applies to any pre-distinguished pair $\e$,
and not only to a principal nilpotent pair.
Thus, we have proved
the following result

\proclaim{Proposition 5.1} $\op{(a)}$ There is a bijection between
the following sets:
$$\Big\{ {\text{nil-pairs}\atop {\text{in }\g}}\Big\}/\Ad G
\enspace\longleftrightarrow\enspace\Big\{
{{\text{Levi subalgebras }{{\frak l}\subset\g,}\text{ and conjugacy}}\atop
{\text{classes of pre-distinguished pairs in }{\frak l}}}\Big\}/\Ad G$$

\noindent
$\op{(b)}$ The semisimple pair $\bh$ 
associated to a pre-distinguished pair $\e$ is  unique, up to
conjugacy by $Z^{\bold {unip}}_G(\e)$. \qed
\endproclaim

\smallskip
\noindent
{\bf Definition 5.2.} A commuting pair $\e=(e_1,e_2)\in\ZZ$ will be
called
{\it distinguished} if the following holds:

(i)$\;$  Lie algebra $\z(\e)\subset \g$ consists of nilpotent elements, and

(ii) There exists a {\it regular} semisimple pair $\bh=(h_1,
h_2)\in \ZZ$ 
such that: 
\vskip 2pt

\centerline{$\
 [h_i, e_j]= \delta_{i,j}\cdot e_i\;,\;
i,j\in \{1,2\}\,.$}

\medskip

Thus the  pair $\e=(e_1,e_2)$ is distinguished if and only if
it is pre-distinguished and the associated semi-simple pair
$\bh$ (which is unique up to conjugacy, due to Proposition 5.1(b))
 is regular. We believe
that a good deal of results that we prove here for 
principal nilpotent pairs have natural generalisations
to distinguished nilpotent pairs. For example,
 we make the following 

\smallskip

\noindent
{\bf Conjecture 5.3.} If $V$ is a rational $G$-module and
$\e$ a distinguished pair, then the space $\tlim_\e\, V^{\z(\bh)}$
is annihilated by $\z(\e)$, and hence is independent of the choice
of an associated semisimple pair $\bh$.
\smallskip

\proclaim{Conjecture 5.4}
The number of $\Ad G$-conjugacy classes of all
distinguished pairs in $\g$ is  finite.
\endproclaim

\noindent
{\bf Remark 5.5.} The following example shows that, even for
$\g=\sln$, not every pre-distinguished pair is distinguished,
and moreover, the classification 
of $\Ad G$-conjugacy classes of 
nil-pairs (hence by Proposition
5.1, of all pre-distinguished pairs) is a wild problem; in particular
there exist continuous
 families of such $\Ad G$-conjugacy classes.

Given a finite dimensional vector space $E$ and
any finite collection $E_i\,,\,i=1,\ldots,r$
of its vector
subspaces, form the following commutative ladder-shaped diagram:
$$
\align
E_1 \longrightarrow & \;\,E\quad {\overset {id} \to\longrightarrow}\quad E
          \quad {\overset {id} \to\longrightarrow}\quad E\\
        & \;\uparrow \quad\enspace\enspace \, {}^{_{id}}\uparrow\quad\enspace 
\enspace\enspace \;{}^{_{id}}\uparrow\\
        & \;E_2\enspace \to \quad E \quad  
      {\overset {id} \to\longrightarrow}\;\quad E\\
        & \enspace
\qquad\quad\enspace\; \uparrow\enspace\quad\enspace\enspace\;\, 
{}^{_{id}}\uparrow\\
        & \enspace\qquad \quad\enspace \,\,E_3\enspace
\longrightarrow\enspace 
\enspace E\\
        & \qquad \quad \qquad \qquad\cdots
\endalign
$$
Let $V$ be the direct sum of all the spaces in the diagram,
let $e_1 : V \to V$ be the map induced by all horisontal arrows in the
diagram, and let $e_2 : V \to V$
be the map induced by all vertical arrows in the
diagram. This way, to any collection $(E, E_1, E_2,\ldots, E_r)$ as
above one associates the nil-pair $(e_1, e_2)$ in $\sll(V)$.
We see that the  classification problem 
of $\Ad G$-conjugacy classes of 
nil-pairs contains as a subproblem the classification
 of $r$-tuples of vector subspaces of a vector space.
This is known to be a wild problem, in general ($4$-tuples
of lines in a 2-plane have
a well-defined cross-ratio, a continuous invariant).\qed

\medskip

From now on, we assume that $\g=\sln$.
We are going to classify all
nil-pairs
$\e=(e_1, e_2)$ in $\sln$ that have a {\it regular}
associated semisimple pair $\bh$.

Let $\lambda\subset \Z\oplus\Z$ be a finite subset, 
thought of as a collection
of boxes on the 2-plane. 
We say that  $\lambda$ is connected if any 
two boxes in $\lambda$ can be joined by a sequence of boxes of $\lambda$,
such that every pair of adjacent boxes in the sequence has a common
edge.  Let $\lambda_1$ and $\lambda_2$
be two Young diagrams with the same southwest corner, and such that
$\lambda_1 \subset\lambda_2$ and $\lambda_2 \smallsetminus\lambda_1$
is a connected $n$-element set. Connected sets of the form 
$\lambda:=\lambda_2 \smallsetminus\lambda_1$ will be referred to
as
{\it skew-diagrams}. We write $(-\lambda)$ for
the set whose boxes have
coordinates opposite to those of $\lambda$.

Given any connected set $\lambda$ with $n$ boxes, we label
the standard base vectors in $\C^n$ by the boxes of the set
(in some way), and
associate to $\lambda$ a pair $\e_\lambda=(e_1, e_2)$
of nilpotent linear transformations of $\C^n$ in the same
way as it has been done in the Introduction.
In more detail,
we let $e_1$ 
act  along the rows of the set, moving one step to the right (if this is
possible
and act by zero otherwise),
and  let $e_2$ 
act along the columns of the set, moving one step up
(if this is
possible
and act by zero otherwise).

\proclaim{Theorem 5.6} 
$\op{(a)}$ A commuting pair $\e$ in $\sln$ is distinguished
if and only if it is conjugate to  a pair of the form
$\e=\e_{_{\pm\lambda}}$, where $\lambda$ is a skew-diagram
with $n$ boxes.

$\op{(b)}$ The pair
 $\e_{_{\pm\lambda}}$ corresponding to a skew-diagram $\lambda$
is a principal nilpotent pair
if and only if the skew-diagram $\lambda$ is a Young diagram.
\endproclaim

We need some notation.
Given a collection of the form
$\blambda= \{\lambda_1\,,\,\lambda_2\,,\,\ldots\,,\,\lambda_r\}$,
where each $\lambda_i$ is either a skew-diagram or 
minus skew-diagram,
 write $|\lambda_i|$ for the number of boxes in $\lambda_i$,
and put $|\blambda|=|\lambda_1|+\ldots+|\lambda_r|$.
We let $\e_\blambda:= \bigoplus_i\;\e_{\lambda_i}\,$ denote the
commuting pair of endomorphisms of
 $\C^{|\blambda|}=\bigoplus_i\;\C^{|\lambda_i|},$
 given by the direct sum.

\proclaim{Lemma 5.7}
A nil-pair
in $\sln$ has a regular
associated semisimple pair if and only if it is conjugate to a pair
of the form $\e_\blambda$, for a certain collection of skew-diagrams
$\blambda=
\{\pm\lambda_1\,,\,\pm\lambda_2\,,\,\ldots\,,\,\pm\lambda_r\}$ 
with $|\blambda|=n$.
\endproclaim

{\it Proof of Lemma.}
We observe that 
any commuting pair of diagonal matrices $\bh=(h_1,h_2)$
gives a bigrading $\C^n =\bigoplus_{p,q}\; V_{p,q}$,
where $V_{p,q}=\{v\in \C^n\;|\; h_1(v)=p\cd v,\,h_2(v)=q\cd v\}$ is a
joint 
weight space of $(h_1,h_2)$.
It is clear that the pair $\bh=(h_1,h_2)$ is regular if and only if
all bigraded components $V_{p,q}$ are at most 1-dimensional. In this
case, one breaks up the set, $\op{Spec}(\bh)$, of all pairs $(p,q)$ such that
$\dim V_{p,q}=1$ into subsets: $\op{Spec}(\bh)
=S_1\sqcup\ldots\sqcup S_m$ 
in such a way that $(p,q)$ and $(p',q')$ belong to the same
subset if and only if $p'-p\in\Z\;\&\; q'-q\in \Z$. Now,
chosing an arbitrary pair $(p_\circ,q_\circ)\in S_i$ one may
identify  each 
subset $S_i$ with the subset in the 2-plane formed by the points
with coordinates $\{(p-p_\circ, q-q_\circ)\}_{\,(p,q)\in S_i}$.
Thus we see that if the 
nil-pair $\e=(e_1, e_2)$
in $\sln$ has a regular
associated semisimple pair, then $\e=\e_\blambda$, for a certain
collection of connected subsets
$\blambda= \{\lambda_1\,,\,\lambda_2\,,\,\ldots\,,\,\lambda_r\}$.

Next, we restrict our attention to  one such  subset $\lambda=\lambda_i$
in the collection.
Observe that, for the maps $e_1$ and $e_2$ in the pair
$\e_\lambda= (e_1, e_2)$ to 
commute, every $2\times 2$-square of boxes must commute.
This automatically holds for any $2\times 2$-square contained
inside $\lambda$. Therefore one only has to check commutativity
of the $2\times 2$-squares that intersect both $\lambda$
and its complement.
It is straightforward to see that commutativity
of all such squares amounts to the requirement that
either $\lambda$ or $(-\lambda)$
has the shape of a skew-diagram. This proves the lemma.\qed\medskip

{\it Proof of Theorem 5.6.} It is clear that the centralizer
of a pair $\e=\e_\blambda$ consisting of more than one
skew-diagram contains a diagonal matrix (that resticts to an arbitrary
scalar on each skew-diagram of the collection 
$\blambda$), hence the pair $\e$
can not
be distinguished in this case. Thus,
to prove part (a)  it suffices  to show that,
for any skew-diagram $\lambda$,
the Lie algebra $\z(\e_\lambda)$ consists
of nilpotent endomorphisms of $\C^n$ only.

Let $\lambda$ be a skew diagram. 
A subset  $\nu\subset \lambda$ will be called an {\it in-subset}
(resp.  {\it out-subset})
if $\nu$ is connected, and for any box $(p,q)\in \nu$,
all boxes
$(i,j)\in\lambda$ such that  $i\leq p$ and
$j\leq q$ (resp. $i\geq p$
and $j\geq q$)
belong to $\nu$ (note that any in-subset has the shape
of a
skew-diagram, and any out-subset has the shape
of a "minus" skew-diagram).
For any pair of integers $p,q$, let
$S_{p,q}(\lambda)$ be the set of pairs $(\nu_{\!_{{\bold {in}}}}, 
\nu_{\!_{{\bold {out}}}})$, where
$\nu_{\!_{{\bold {in}}}}$ is an in-subset of $\lambda\,,$
$\nu_{\!_{{\bold {out}}}}$ is an out-subset of $\lambda$ and, moreover,
$\nu_{\!_{{\bold {out}}}}$ is obtained from $\nu_{\!_{{\bold {in}}}}$ via translation 
 $\,T: \nu_{\!_{{\bold {in}}}} \iso \nu_{\!_{{\bold {out}}}}$
by the vector $(p,q)$ (translation 
of the plane $p$ steps horisontally
 and $q$ steps vertically).
Thus, $\nu_{\!_{{\bold {in}}}}$ and $\nu_{\!_{{\bold {out}}}}$ have the same shape.

Put $n=|\lambda|$, and label base vectors 
of $\C^n$ by the boxes of $\lambda$.
Given a pair $\bnu=(\nu_{\!_{{\bold {in}}}},
 \nu_{\!_{{\bold {out}}}})\in S_{p,q}(\lambda)$, we define
a linear map $f_\bnu : \C^n\to\C^n$ as follows.
The map $f_{\bnu}$ takes every base vector labelled by
a box $(i,j)\in \nu_{\!_{{\bold {in}}}}$ to the base vector labelled by
the translated box, $(i+p,j+q)\in \nu_{\!_{{\bold {out}}}}$, 
and takes all other
base vectors
(i.e., those labelled by the boxes of $\lambda \smallsetminus
\nu_{\!_{{\bold {in}}}}$) to $0$. Now write $\z(\e_\lambda)=\bigoplus_{p,q}\;
\z_{_{p,q}}(\e_\lambda)$ for the bigrading on the centraliser of
the pair $\e_\lambda$.
Verification of the following crucial observation is straightforward
and is left to the reader.

\proclaim{Claim 5.8} For any skew-diagram (or minus skew-diagram)
$\lambda$, and any
$p,q\geq 0$, the endomorphisms $\{f_\bnu\;|\;
\bnu=(\nu_{_{{\bold {in}}}}, 
\nu_{\!_{{\bold {out}}}})\in S_{p,q}(\lambda)\}$ form a basis of 
$\z_{_{p,q}}(\e_\lambda)$\qed
\endproclaim

Next, given $\bnu=(\nu_{\!_{{\bold {in}}}},\nu_{\!_{{\bold {out}}}})\in
S_{i,j}(\lambda)$, such that
$\nu_{\!_{{\bold {out}}}}$ is obtained from
$\nu_{\!_{{\bold {in}}}}$ via translation $(i,j) \to (p+i,q+j)$,
define  an operator
$T_\bnu:
\nu_{\!_{{\bold {in}}}} \to \nu_{\!_{{\bold {out}}}}$
to be the restriction of that translation to
 the set $\nu_{\!_{{\bold {in}}}}$.
Given a second pair
$\bnu'=(\nu_{\!_{{\bold {in}}}}', 
\nu_{\!_{{\bold {out}}}}')\in S_{i',j'}(\lambda)$, one finds by a direct calculation
$$
\align f_\bnu  \cdot 
f_{\bnu'} =
f_{\bnu''}\,, \quad&\text{where}\quad \bnu''=(\nu_{\!_{{\bold {in}}}}'',\,
\nu_{\!_{{\bold {out}}}}'')\quad\text{is defined by:}\\
&\nu_{\!_{{\bold {in}}}}''=
T^{^{-1}}_{\bnu'}(\nu_{\!_{{\bold {in}}}}\cap\nu_{\!_{{\bold {out}}}}') \subset
\nu_{\!_{{\bold {in}}}}'\quad\text{and}\quad
\nu_{\!_{{\bold {out}}}}''=T^{^{-1}}_\bnu(\nu_{_{{\bold {in}}}}\cap
\nu_{\!_{{\bold {out}}}}') \subset \nu_{\!_{{\bold {out}}}}\,.
\endalign$$

Observe next that since
the skew-diagram $\lambda$ is connected, it can not contain
a subset which is both an in-subset and an out-subset at the same time.
It follows that $\nu_{\!_{{\bold {in}}}}\neq\nu_{\!_{{\bold {out}}}}'$.
 Hence,  we deduce: $\;|\nu_{\!_{{\bold {in}}}}''|<
\max(|\nu_{\!_{{\bold {in}}}}|,|\nu_{\!_{{\bold {in}}}}'|)$.
This implies that any product of the maps of the form
$f_{\bnu}$
 that contains more than $n$ factors vanishes,
and part (a) of Theorem 5.6 follows.

Let the
pair $\e_\lambda=(e_1, e_2)$ be associated to a Young diagram
$\lambda$.
 We claim that
$(e_1, e_2)$ is a principal nilpotent pair. The only non-obvious
part of (1.1) that needs to be verified is the equation
$\dim\z(\e)=\rk\g$. To prove this, note that the base
vector $v_\circ\in\C^n$ labelled by the box with coordinates $(0,0)$,
the south-west ``corner'' 
of the diagram, is a cyclic vector for the operators
$e_1, e_2$, in the sense that $\C^n=\C[e_1, e_2]\cdot v_\circ$,
where $\C[e_1, e_2]$ denotes an abstract polynomial ring in
the variables $e_1, e_2$.
Therefore, any operator $x\in \z(\e)$ is completely determined
by the vector $x(v_\circ)\in \C^n$. It follows that any element
of $\z(\e)$ can be expressed as a polynomial without constant term
in the operators
$e_1$ and $e_2$. Hence, $\z(\e)+ \C\cdot\op{Id} \, =\,\C[e_1, e_2]/I$,
where $I$ is the ideal annihilating $v_\circ$. It is easy to see that
this ideal is
spanned by all monomials $e_1^r e_2^s$, such that 
$(r,s)$ are
not coordinates of a box
of $\lambda$. Thus, $\dim (\C[e_1, e_2]/I) = n$,
hence: $\dim \z(\e) = n-1$, so that $\e$ is a regular pair.
Finally, it is easy to verify using Claim 5.8 that, for any
skew-diagram (resp. minus skew-diagram)
$\lambda$ which is {\it not} a Young diagram (resp. minus Young diagram),
one has $\dim\z(\e_\lambda) > n-1.$ Part (b) of Theorem 5.6
follows.\qed\medskip

Theorem 5.6 shows that the pairs corresponding to
skew-diagrams may be thought of as "double-analogues" of
nilpotent matrices that have
 a single Jordan block in their Jordan form.
\bigskip

\noindent
{\bf Example 5.9. } We give an example of distinguished pairs
in the orthogonal Lie algebra $\so_n$. Let $\lambda$ be a skew
diagram with $n$ boxes, and $\e_\lambda=(e_1,e_2)$ the
pair of endomorphisms of $\C^n$ constructed as above.
Assume, in addition, that the diagram $\lambda$ is
{\it centrally symmetric} with respect to the origin
$(0,0)\in\Z^2$, i.e., $(p,q)\in \lambda\; \Longrightarrow$
$\;(-p, -q)\in \lambda$. We define a symmetric bilinear
form on $\C^n$ as follows:
$$
\omega(v_{_{p,q}}\,,\,v_{_{p',q'}}) =
\cases (-1)^{p+q}\quad\text{if}\quad
p+p'=0=q+q'\\
\; 0\hskip 15mm\text{otherwise}\endcases
$$
where $v_{_{p,q}}$ denotes the base vector in $\C^n$
corresponding to a box $(p,q)\in\lambda$. It is
clear that the form $\omega$ is non-degenerate,
and that both $e_1$ and $e_2$ are skew-symmetric
relative to $\omega$, hence form a commuting
pair in the Lie algebra $\so(\C^n, \omega)$.
The proof of Theorem 5.6 implies that the centralizer
of the pair $\e_\lambda$ in $\sln$, hence in
$\so_n$,  consists of nilpotent matrices. Thus,
to any centrally symmetric skew diagram 
$\lambda$ we have attached a distinguished pair
in $\so_n$. We do not know if every
distinguished pair
in $\so_n$ is obtained in this way. Neither do we know
for which $\lambda$, apart from  rectangular ones,
the pair $\e_\lambda$ is a principal nilpotent pair
in $\so_n$ .
\bigskip

From now until the end of this section,
we will assume that $\g$ is the reductive Lie algebra
$\gl_n$,
 rather than $\sln$. This slight modification will become 
more convenient shortly.
Write $\h$ for the Cartan subalgebra of diagonal 
$n\times n$-matrices, and identify it with $\C^n$.
Set $\uu=(u_1,\ldots,u_n)\,,$
$\vv=(v_1,\ldots,v_n)\,,\,$
 and view $\SS(\h\oplus\h)$ as the
polynomial ring, $\C[\uu,\vv]$,
in $2n$ variables. We let the Symmetric group 
$W=S_n$ act  on the $n$-tuples $\uu$ and $\vv$
by permutations.

Fix a Young diagram $\lambda$ with $n$ boxes,
let $\e=\e_\lambda\in\gl_n$ be the corresponding principal nilpotent
pair.
To write
 an associated pair $(h_1,h_2)$ of diagonal matrices
explicitly,
enumerate  $n$ boxes of the diagram $\lambda$ in some order, and write
$(a_i,b_i)$
for the coordinates of the $i$-th
box, starting counting the coordinates with $(0,0)$, see 
{\it fig. }(0.1).
This way we get a collection of non-negative integers
$a_1,\, b_1,\, a_2,\, b_2, \ldots, a_n,\, b_n$. 
One may then choose an associated semisimple
pair $\bh=(h_1,h_2)$ to be: $h_1=\text{\it diag}(a_1, \ldots, a_n)$
and $h_2=\text{\it diag}(b_1, \ldots, b_n)$. Note that these
diagonal matrices do {\it not} have zero trace, hence are not
in $\sln$. The corresponding semisimple pair in $\sln$ is: 
$(h_1 - \varkappa_1, h_2 -\varkappa_2),$ where $\varkappa_i=
1/n\cdot\op{Tr}(h_i)\cdot \op{Id}$.

All the considerations of the previous sections extend
to reductive Lie algebras without troubles. In particular, for
$\g=\gl_n$,
 we have the element $\P_{\bh}\in
\widehat{\C[\uu,\vv]}$, which is
 explicitly given by the formula:
$$ \P_{\bh}(\uu,\vv)=\sum_{w\in S_n}\;
\varepsilon(w)\cdot e^{w(\bh)}=\sum_{w\in S_n}\;
\varepsilon(w)\cdot e^{\sum_i\;(a_{w(i)} u_{w(i)} + b_{w(i)}
v_{w(i)})}\,.
\tag 5.10$$

A special feature of the $\gl_n$-case is that
one can make  a change of variables by setting:
$x_i:= e^{u_i}$ and $y_i:=e^{v_i}\,,\,
i=1,\ldots,n\,.$ Then, $\P_{\bh}$ becomes
a  polynomial (as opposed to the exponential expression above)
in the new variables: $\x=(x_1,\ldots,x_n)$ and
$\y=(y_1,\ldots,y_n)$. In view of (5.10), this
 polynomial may be 
written in the form of the following  determinant:
$$
\P_{\bh}(\x,\y) = \det \vmatrix 
x^{a_1}_1 y^{b_1}_1 & x^{a_1}_2 y^{b_1}_2 & \hdots & x^{a_1}_n y^{b_1}_n \\
\text{   } \\
x^{a_2}_1 y^{b_2}_1 & x^{a_2}_2 y^{b_2}_2 & \hdots & x^{a_2}_n y^{b_2}_n \\
\hdotsfor 4 \\
\text{   } \\
x^{a_n}_1 y^{b_n}_1 & x^{a_n}_2 y^{b_n}_2 & \hdots & x^{a_n}_n
y^{b_n}_n \endvmatrix
\tag 5.11
$$
where $\d_1=\sum a_i= \op{Tr}(h_1)$ and $\d_2=\sum b_i= \op{Tr}(h_2)$.
The determinant on the RHS has been first introduced by
P.A. MacMahon  {\it "Combinatory Analysis"},
{\sl Cambridge University Press} (1916),
as a "double analogue" of the Vandermonde determinant;
see also the paper by A.M. Garsia and  M. Haiman
{\it "A graded representation model for Macdonald's polynomials"},
{\sl Proc. Nat. Acad. Sci. USA}, {\bf {90}} (1993), 3607--3610.
\medskip

We claim next that the $W$-harmonic polynomial $\P_\e\in 
\C^{\d_1}[\uu] \bigotimes \C^{\d_2}[\vv]$ is given by 
the same determinant (5.11), where now $\x$ is replaced by $\uu$,
and $\y$ by $\vv$, that is we have

\proclaim{Proposition 5.12} In the $\gl_n$-case the elements
$\P_\bh$ and $\P_\e$ can be obtained from each other
by a change of variables, i.e.:
$\dis\P_\e(e^\uu,e^\vv)=\frac{1}{\d_1!\cdot \d_2!}\cdot
\P_\bh(\uu,\vv)\,.$
\endproclaim

This explains Corollary 4.16. We see also that
the "difficult" non-vanishing part of Theorem 4.4(ii) follows,
in the $\gl_n$-case, directly from the non-vanishing of $\P_\bh$,
which is trivial.
\medskip

{\it Proof of Proposition 5.12.} Our argument will be based on the
identity :
$$
\sum_{y \in W^{i}}\,\eps(y)\cdot e^{y(h_i)}\,\Big|_{e^\uu=\x}
\!=
\prod_{\alpha\in R^{_i}_+}\langle\alpha^\vee\!,\x\rangle =\pi_i(\x)=
\frac{1}{\d_i!}\sum_{y \in W^{i}}\,
\eps(y)\cdot \langle y(h_i), \x\rangle^{\d_i}.
\tag 5.13$$
In this identity, $i=1,2$, the LHS is the standard Weyl denominator for the
Levi subalgebra $\g^{_i}$ which,
 for $\g=\gl_n$, is equal to a product of Vandermonde determinants
corresponding to the columns (if $i=1$), resp. rows (if $i=2$),
 of the
Young diagram. It is clear that, upon substitution
$e^\uu=\x$, each Vandermonde determinant becomes the
corresponding polynomial $\pi$, the product of positive coroots that
stands
at the middle of formula (5.13). The equation on the right of
(5.13)
is nothing but  identity (4.5).

In the formulas below we abuse the notation
slightly so that $w({\langle
h ,\uu\rangle})$ should be understood as 
$w$ acting on either of the two variables $h$ or $\uu$, that is as:
$\langle w(h), \uu\rangle = \langle h, w(\uu)\rangle$,
and $w(e^{\langle
h ,\uu\rangle})$ has a similar  meaning. With this understood, we
calculate
$$
\align
\P_\bh & = \sum_{w\in W}\,\eps(w)\cdot w(e^{\langle h_1,\uu\rangle}\otimes
e^{\langle h_2,\vv\rangle})\\
&=\sum_{w\in W}\eps(w)\cd\frac{1}{\d_1!}\cdot
w\Bigl(\bigl(\!\!\sum_{w_1\in W^{1}}\,\eps(w_1)\cd 
e^{\langle w_1(h_1),
\uu\rangle}\bigr)\otimes e^{\langle
h_2,\vv\rangle}\Bigr)\,\Big|_{e^\uu=\x} =
 \quad\text{by (5.13)}\\
&=\frac{1}{\d_1!}\cdot\sum_{w\in W}\eps(w)\cd w
\Bigl(\bigl(\frac{1}{\d_1!}\cdot\sum_{w_1 \in W^{1}}\,
\eps(w_1)\cdot \langle w_1(h_1), \x\rangle^{\d_1}\bigr)
 \;\bigotimes\; e^{\langle
h_2,\vv\rangle}\Bigr)\\
&=\frac{1}{\d_1!}\cdot\sum_{w\in W}\,\eps(w)\cdot
\langle w(h_1), \x\rangle^{\d_1}\bigotimes e^{\langle
w(h_2),\vv\rangle}\\
&=\frac{1}{\d_1!}
\sum_{w\in W}\eps(w)\cdot w
\Bigl(\langle h_1, \x\rangle^{\d_1}\bigotimes \,
\bigl(\frac{1}{\d_2!}\sum_{w_2\in W^{2}}\,\eps(w_2)\cdot e^{\langle
w_2(h_2),\vv\rangle}\bigr)\Bigr)\,\Big|_{e^\vv=\y}\\
& =\quad
\text{by (5.13)}\quad =\frac{1}{\d_1!\cdot\d_2!}
\sum_{w\in W}\eps(w)\cdot
\langle w(h_1), \x\rangle^{\d_1}\bigotimes \,
\langle w(h_2), \y\rangle^{\d_2}\\
& =
\frac{1}{\d_1!\cdot\d_2!}\cdot \P_\e(\x,\y)\;.\qed
\endalign
$$
\medskip

\bigskip

\head{6. Some cohomology and generating functions.}
\endhead
\bigskip

Fix a principal nilpotent pair $\e=(e_1, e_2)$,
an associated semisimple
pair $\bh=(h_1,h_2)$, and let $\g=\bigoplus_{p,q}\, \g_{_{p,q}}$
be the corresponding bigrading.

We will be concerned with the following 3-term complex:
$$\g\; \overset {\partial'} \to \longrightarrow\;
\g\oplus\g\; \overset {\partial''} \to \longrightarrow\; \g
\qquad,\qquad
\g_{_{p-1,q-1}}\; \overset {\partial'_{^{p,q}}} \to \longrightarrow
\g_{_{p,q-1}}\oplus\g_{_{p-1,q}}\; \overset {\partial''_{^{p,q}}} 
\to \longrightarrow\; \g_{_{p,q}}\,,
\tag 6.1
$$
where the differentials are given by the formulas
$$\partial'\; :=\; \ad e_1\; \oplus\; \ad e_2\;,
\quad\text{and}\quad
\partial''\; :=\; (\ad e_2)\; \oplus\; (-\ad e_1)\,.$$
It is clear that $\partial'\,\circ\,\partial'' =0$,
and we put 
$H_{_{p,q}}:=\op{Ker}(\partial''_{^{p,q}})/\op{Im}(\partial'_{^{p,q}})$.

 Observe that the two maps below are adjoint to each other with respect
to the Killing form (see Corollary 1.8):
$$
\partial'_{^{-p+1,-q+1}}: \g_{_{-p,-q}} \to
\g_{_{-p,-q+1}}\oplus\g_{_{-p+1,-q}}\;\;
\quad,\quad\;\;
\partial''_{^{p,q}}:
\g_{_{p,q-1}}\oplus\g_{_{p-1,q}}\to \g_{_{p,q}}\,.
$$
Since
$\;\displaystyle\op{Ker}(\partial'_{^{p,q}})=\z_{_{p-1,q-1}}(\e)$,
we get by duality $\;\displaystyle
\op{Coker}(\partial''_{^{p,q}}) \simeq\z^*_{_{-p,-q}}(\e)\,,\;
$ and $\;\displaystyle
\op{Im}(\partial''_{^{p,q}}) =\z_{_{-p,-q}}(\e)^\perp
\simeq$
$\g^*_{_{-p,-q}}/\z^*_{_{-p,-q}}(\e).$
From these formulas  we deduce:
$$
H_{_{p,q}}\;=\;
\frac{\op{Ker}(\g_{_{p,q-1}}\oplus\g_{_{p-1,q}}\twoheadrightarrow
\g^*_{_{-p,-q}}/\z^*_{_{-p,-q}})}{
\op{Im}\left(\g_{_{p-1,q-1}}/ \z_{_{p-1,q-1}}
\hookrightarrow \g_{_{p,q-1}}\oplus\g_{_{p-1,q}}\right)}
\,.\tag 6.2
$$

To put the complex (6.1) in an adequate context,  we regard
$\g$  as a
module
over an abstract two-dimensional abelian Lie algebra ${\frak e}$ whose base
vectors
act on $\g$ via the operators $\{\ad e_i\}_{\,i=1,2}$. Then (6.1)
becomes the standard Koszul complex:
$\; \g\to {\frak e}^*\otimes\g \to (\wedge^2{\frak e}^*)\otimes\g\,\,,$
computing the Lie algebra cohomology 
$H^{1}({\frak e}, \g)$.

The group $H^{1}({\frak e}, \g)$ has the following geometric
interpretation in terms of the commuting variety $\ZZ$. 
First, it is easy to verify that condition 1.1(Reg) insures
that $\e$ is a smooth point of $\ZZ$. Assume further (without 
justification) that  the space of
diagonal $\Ad G$-orbits on $\ZZ$ has the structure of a 
smooth algebraic variety, $\ZZ/\Ad G$,
at least locally near $\e$. By general principles, 
the tangent space to $\ZZ/\Ad G$ at the point $\co=\Ad G(\e)$
is then given by:
$\,T_{_{\co}}(\ZZ/\Ad G)= H^{1}({\frak e}, \g)$. In other words,
$H^{1}({\frak e}, \g)$ is the normal space 
at a point $\x\in\co$ to
the submanifold $\co$ inside~$\ZZ$.
\medskip

Next, we introduce the 
following generating functions in two variables: $H_{_\e}(s,t):=$
$$ \sum_{p,q}\; s^pt^q\cdot\dim H_{p,q}\;,\quad
\g_{_\e}(s,t):=\sum_{p,q}\; s^pt^q\cdot\dim\g_{_{p,q}}\;,\quad
\z_{_\e}(s,t):=\sum_{p,q}\; s^pt^q\cdot\dim\z_{_{p,q}}(\e)\,.$$

\proclaim{Lemma 6.3} $\enspace H_{_\e}(s,t)
= st\cd\z_{_\e}(s,t) +\z_{_\e}(s^{-1},t^{-1})-
(s-1)(t-1)\cd \g_{_\e}(s,t)\,.$
\endproclaim

{\it Proof.} 
From (6.2) we find: $\displaystyle
\;\dim H_{p,q}=\dim\g_{_{p,q-1}}+  \dim\g_{_{p-1,q}}-$
$$\hskip 30mm -\dim\g_{_{p-1,q-1}} - \dim\g_{_{-p,-q}}+  
  \dim  \z_{_{p-1,q-1}}(\e)+\dim\z_{_{-p,-q}}(\e)\,.$$
Multiply each side by $s^pt^q$, and use 
that $\dim\g_{_{-p,-q}} =\dim\g_{_{p,q}}$, by Corollary 1.8.
 Summing up over all $p,q\in \Z$, we
find
$$H_{_\e}(s,t)
=t\cd \g_{_\e}(s,t) +s\cd  \g_{_\e}(s,t)  - st\cd \g_{_\e}(s,t)- \g_{_\e}(s,t)+
st\cd \z_{_\e}(s,t) + \z_{_\e}(s^{-1},t^{-1})\,.\;\; \qed
$$\medskip

Recall
that the centralizer of $e_i$, $(i=1,2)$, is a bigraded
subalgebra $\z(e_i)= \bigoplus_{p,q}\;\z_{_{p,q}}(e_i)\,$
$\subset\g$. The structure of the cohomology 
$H^{1}({\frak e}, \g)$ is described
by

\proclaim{Theorem 6.4} $\op{(i)}\quad
\dim  H^{1}({\frak e}, \g) =2\rk\g\,,$ and $H_{p,q}=0$ whenever $p\cdot q\ge 0$.

$\op{(ii)}$ Furthermore, we have
$$
H_{p,q}\,\simeq\,\cases \op{Coker}\left(\ad e_1: \z_{_{p-1,q-1}}(e_2) \to
\z_{_{p,q-1}}(e_2)\right)
\quad\op{if}\quad p<0\,\, \&\,\, q\geq 0\\
\op{Coker}\left(\ad e_2: \z_{_{p-1,q-1}}(e_1) \to
\z_{_{p-1,q}}(e_1)\right)
\quad\op{if}\quad p\geq 0\,\, \&\,\, q<0
\endcases
$$
\endproclaim

{\it Proof.} 
Let $z\in\z_{_{p,q-1}}(e_2)$. Then, for the element
$h_{p,q} := (z,0)\in \g_{_{p,q-1}}\oplus\g_{_{p-1,q}}$ we clearly have:
$\partial''_{^{p,q}}(h_{p,q})=\ad e_2(z)-\ad e_1(0)=0\,.$
Hence, $h_{p,q}$ represents a class in $H_{p,q}$, see (6.1).
The element $h_{p,q}$ is  a coboundary if and only if
there exists
$y\in \g_{_{p-1,q-1}}$ such that:
$\ad e_2(y)=0$ and $\ad e_1(y)=z$.
The first equation above implies that $y\in \z_{_{p-1,q-1}}(e_2)$.
Therefore, the second equation implies that the assignment
$z\mapsto (z,0)$ gives an injection 
$$\text{Coker}\left(\ad e_1: \z_{_{p-1,q-1}}(e_2) \to
\z_{_{p,q-1}}(e_2)\right) \into H_{p,q}\,.\tag 6.7$$

We now 
study in more detail the case when  $p<0\,\, \&\,\, q\geq 0$.
We observe that the operator $\ad e_1$ in (6.7) indeed maps $\z_{_{p-1,q-1}}(e_2)$
into $\z_{_{p,q-1}}(e_2)$,
since $e_1$ commutes with $e_2$. Moreover, for all
$p<0$, this map is injective,
due to the weak Lefschetz. This way we see, applying $\ad e_1$ several
times
that, for any $i=-p>0$, we have an injective map
$ \ad^ie_1: \z_{_{p,q-1}}(e_2) \into
\z_{_{0,q-1}}(e_2)\subset\z(h_1,e_2)\,.$
Thus, the sequence of spaces
$F^i:= \bigoplus_{q\geq
0}\;\ad^ie_1\bigl(\z_{_{-i,q-1}}(e_2)\bigr)\,,$
${\,i=0,1,\ldots}\,,$
gives a {\it decreasing} filtration of the subalgebra $\z(h_1,e_2)$,
which is, in a sense, dual to the {\it increasing} filtration $F_\bullet$
considered in $\S2$.
Observe further that the weak Lefschetz insures that
$\z_{_{0,q}}(e_2)=0$, for all $q<0$. Hence, summation over all $q$
yields, for any $i>0$, a vector 
space isomorphism 
$$
\ad^ie_1: \text{Coker}\left(\ad e_1: \z_{_{-i-1,*}}(e_2) \to
\z_{_{-i,*}}(e_2)\right) \iso F^i/F^{i-1}\,. 
$$
Thus, writing
$\gr\,\z(h_1,e_2) :=\bigoplus_i\;F^i/F^{i-1}\,,$
we obtain a graded space isomorphism 
$$
\bigoplus_i\;\Bigl(\text{Coker}\left(\ad e_1: \z_{_{-i,*}}(e_2) \to
\z_{_{-i,*}}(e_2)\right)\Bigr)\;\; \simeq\;\;
\gr\,\z(h_1,e_2)\,.\tag 6.8
$$

Recall now that $\dim\z(h_1,e_2)=\dim\z(\e)=\rk\g$, by
(2.3).  Hence, we get $\dim\bigl(\gr\,\z(h_1,e_2)\bigr)=\dim\z(h_1,e_2)$
$=\rk\g\,.$ It follows now from (6.7) and (6.8) that we have:
$\rk\g=\dim\bigl(\gr\,\z(h_1,e_2)\bigr)\leq\dim
\bigl(\bigoplus_{p<0,q\geq 0}\;H_{p,q}\bigr)\,.$
In the same way one obtains a similar inequality
 for the southeast quadrant:
$p\geq 0\,\, \&\,\, q<0$, i.e., $\rk\g\leq\dim
\bigl(\bigoplus_{p\geq 0,\, q<0}\;H_{p,q}\bigr)$.
The two inequalities combined together imply that
$$2\cdot\rk\g\enspace\leq\enspace
\dim\left(\bigoplus\nolimits_{p\cdot q <0}\;\;H_{p,q}\right)\,.
\tag 6.9
$$

On the other hand, setting $s=t=1$ in the formula of Lemma 6.3,
we find
$\dim  H^{1}({\frak e}, \g) = 2\cd \dim\z(\e) = 2\rk\g\,.$
It follows that the inequality in (6.9) has to be an equality.
This implies that the inclusion in (6.7),
and a similar inclusion for the quadrant
 $p\geq 0\,\, \&\,\, q<0$, are  in effect both isomorphisms.
Furthermore, the groups $H_{p,q}$ must vanish for the
two other quadrants, that is for all $p\cdot q\geq 0$.
This completes the proof of the theorem.\qed\medskip

Part (i) of the theorem implies, in view of
formula (6.7) that, for any $p,q \geq 0$,
the LHS of (6.7)  vanishes. This gives

\proclaim{Corollary 6.10} The maps $
\ad e_1: \z_{_{p,q}}(e_2) \to
\z_{_{p+1,q}}(e_2)$ and $\ad e_2: \z_{_{p,q}}(e_1) \to
\z_{_{p,q+1}}(e_1)$ are both surjective,
for any $p,q \geq 0.$\qed
\endproclaim

Write $\z_{_{p,q}}(e_1\cdot e_2):=\op{Ker}\bigl(\ad e_1\cd\ad e_2:
\g_{_{p,q}}\to\g_{_{p+1,q+1}}\bigr)\,.$
Corollary 6.10 may be reformulated in the following way

\proclaim{${\boldsymbol\partial}{\overline{\boldsymbol\partial}}$-lemma 6.11}
$\enspace \z_{_{p,q}}(e_1\cdot e_2)=\z_{_{p,q}}(e_1)\,+\,\z_{_{p,q}}(e_2)\;,\,\forall
p,q \geq 0.$
\endproclaim

{\it Proof.} Let $x\in \g_{_{p,q}}$ be such that 
$\ad e_1\cd\ad e_2(x)=0$. 
Put $y:= \ad e_2(x)$. Then $y\in \z_{_{p,q+1}}(e_1)$.
Hence, by Corollary 6.10, there exists 
$z\in\z_{_{p,q}}(e_1)$ such that $y=\ad e_2(z)$.
Therefore we have: $\ad e_2(x-z)=y-y=0.$
Thus, $x= z + (x-z) \in \z_{_{p,q}}(e_1)\,+\,\z_{_{p,q}}(e_2)\,,$
and the result is proved. 

Conversely, one may show that  the Lemma implies Corollary 6.10.\qed
\medskip

By Theorem 6.4 we can write $H^{1}({\frak e}, \g) = H_{\!_{{\bold {nw}}}}
\,\bigoplus\,
H_{\!_{{\bold {se}}}},$
where {\bf nw} = northwest, {\bf se} = southeast,
and we put
$\,H_{\!_{{\bold {nw}}}} :=$ $\dis\bigoplus_{p<0,q\geq 0}\;H_{p,q}\,,$ and
$\,H_{\!_{{\bold {se}}}}:=$ $\dis\bigoplus_{p\geq 0, q<0}\;H_{p,q}.$

\proclaim{Corollary 6.12} The Killing form on $\g$ induces a perfect
pairing: $H_{\!_{{\bold {nw}}}} \;\times\; H_{\!_{\bold {se}}} 
\to \C$. In particular, the space
$H^{1}({\frak e}, \g) \simeq  H_{\!_{{\bold {nw}}}}\,\bigoplus\,
 H^*_{\!_{{\bold {nw}}}}$ has a natural
symplectic
structure.
\endproclaim

{\it Proof.} Corollary 1.8 insures that the two squares below
are obtained from each other by duality induced by the Killing form:
\vskip 2pt
\centerline{$\hphantom{x}\enspace\enspace
\g_{_{p-1,q}}  \;\overset {{}_{e_1}} \to 
\longrightarrow \;\;
\g_{_{p,q}}
\qquad\qquad\qquad\enspace\enspace
\g_{_{-p,-q+1}}  \;\overset {{}_{e_1}} \to 
\longrightarrow \;
\g_{_{-p+1,-q+1}}
$}
\centerline{$\uparrow {}_{^{e_2}}\qquad\,\,\quad \;
\;  \uparrow {}_{^{e_2}}
\qquad\qquad\qquad\;\;\,
\uparrow {}_{^{e_2}}\qquad\;\;\;\qquad   \uparrow {}_{^{e_2}}$}
\centerline{$\hphantom{x}\;\enspace
\g_{_{p-1,q-1}} \;\;\; \overset {{}_{e_1}} \to 
\longrightarrow \;\;\; \g_{_{p,q-1}}
\qquad\quad\qquad\;
\g_{_{-p,-q}} \;\; \overset {{}_{e_1}} \to 
\longrightarrow \;\;\;\; \g_{_{-p+1,-q}}
$}
\vskip 2pt
\noindent
where "$e_i$" stands for the map given by the  $\ad e_i$-action.

The duality of the above squares induces a perfect duality between the
following two spaces: 
$$
\text{Coker}\bigl(\text{Ker}_{_{p-1,q-1}}(e_2)\, \overset {{}_{e_1}} \to 
\longrightarrow\,
\text{Ker}_{_{p,q-1}}(e_2)\bigr)\enspace,
\enspace
\text{Coker}\bigl(\text{Ker}_{_{-p,-q}}(e_1)\, \overset {{}_{e_2}} \to 
\longrightarrow\,
\text{Ker}_{_{-p,-q+1}}(e_1)\bigr)\,.$$
The result follows.\qed

\bigskip

We complete this section with a few numerical identities
involving bi-exponents of $\g$ relative to~$\e$.

\proclaim{Proposition 6.13} The following identities hold
 $$\dim\g^{_1}_{_+}=\sum_{p,q}\,p\cdot\dim\z_{p,q}(\e)\quad\text{and}\quad
\dim\g^{_2}_{_+}=\sum_{p,q}\,q\cdot\dim\z_{p,q}(\e)\,.$$
\endproclaim

{\it Proof.} Recall isomorphism (2.3). The spaces $\z(h_1,e_2)\,,\,
\tlim_\e\,\z(h_1,e_2)$,
and $\z(\e)$ entering (2.3) are all stable under the $\ad h_2$-action.
Further, since $e_1$ commutes with $h_2$, the limit construction in
(2.3), hence the isomorphism itself, commutes with the $\ad h_2$-action.
Separating different
weight components of $\ad h_2$, and using Lemma 2.2(ii) we get
an isomorphism:
$$\gr^{\,}\z_{_{0,q}}(h_1,e_2) 
\;\overset {\sim} \to 
\longrightarrow\; \z_{*,q}(\e)
\enspace,\enspace \forall  q\in\Z.$$
Hence,
$\,\dim\z_{*,q}(\e) =\dim \bigl(\gr\,\z_{_{0,q}}(h_1,e_2)\bigr)
=\dim\z_{0,q}(e_2)\,,$
$\forall q.$ Summing up over all $q$ we obtain:
$\,\sum_q\,q\cdot\dim\z_{*,q}(\e) =\sum_q\,q\cdot\dim \z_{0,q}(e_2).$
Clearly, we may concentrate on the terms with $q>0$.
Observe that $\z_{0,q}(e_2)$
is by definition the $q$-eigenspace of the operator
$\ad h_2$ acting on
$\z_{\g^{_2}}(e_2)$. If
$q>0$ we may replace $\g^{_2}$ here by $\s^{_2}$,
and also replace the $\ad h_2$-action by that of
$\ad s_2$, c.f. Proposition 3.2. By parts (iii)-(iv) of that Proposition,
$e_2$ is a principal
nilpotent in $\s^{_2}$,   and the theory of principal $\ll2$-triples [K1]
applied to the Lie algebra $\s^{_2}$
says that one has an identity: $\sum_{q>0}\,q\cdot\dim \z_{0,q}(e_2)
=\dim \g^{_2}_{_+}.$ The first equation is proved similarly. \qed
\medskip
\bigskip

\proclaim{Proposition 6.14} One has:
$$
\prod_{(p,q)\in\text{Exp}_\e(\g)}
\frac{1- s^{p+1}t^{q+1}}{1-s\cdot t}\;=\;
\prod_{\alpha\in R_{\ne}}
\frac{1- s^{\langle\alpha , h_1\rangle+1}t^{\langle\alpha , h_2\rangle+
1}}{1-
 s^{\langle\alpha , h_1\rangle+1}t^{\langle\alpha , h_2\rangle}}
\cdot\frac{
1- s^{\langle\alpha , h_1\rangle}t^{\langle\alpha , h_2\rangle}}{1-
 s^{\langle\alpha , h_1\rangle}t^{\langle\alpha , h_2\rangle+1}}\;.$$
\endproclaim

{\it Proof.} The LHS of the formula above equals:
$$\prod_{i,j\geq 0}\;(1- s^it^j)^{\dim\z_{_{i,j}}(\e)}\;. \tag 6.15$$
Put $g_{_{i,j}} := \dim\g_{_{i,j}}$.
The RHS of the formula above can be rewritten as
$$\prod_{i,j\geq 0}
\Bigl(\frac{(1- s^{i+1}t^{j+
1})\cdot (1- s^it^j)}{(1- s^{i+1}t^{j})\cdot 
(1- s^it^{j+1})}\Bigr)^{g_{_{i,j}}}=
\prod_{i,j\geq 0}
\frac{(1- s^it^j)^{g_{_{i+1,j+1}}}\cdot
(1- s^it^j)^{g_{_{i,j}}}}{
(1- s^it^j)^{g_{_{i+1,j}}}\cdot
(1- s^it^j)^{g_{_{i,j+1}}}}\;.\tag 6.16
$$

Theorem 6.4 yields an exact sequence:
$$ 0\to\z_{_{i,j}}(\e)\to \g_{_{i,j}}\to
\g_{_{i+1,j}} \oplus \g_{_{i,j+1}} \to
\g_{_{i+1,j+1}}\to 0$$
Hence,
$\dim\g_{_{i,j}} -\dim\g_{_{i+1,j}} -\dim\g_{_{i,j+1}}
+\dim\g_{_{i+1,j+1}} = \dim \z_{_{i,j}}(\e)\,,$
and the RHS of (6.16) reduces to: $\prod_{i,j\geq 0}\;(1-
s^it^j)^{\dim\z_{_{i,j}}(\e)}$,
which is (6.15).
\qed
\medskip

The proposition above reduces in the special
case $\e=(e,0)$, where $e\in\g$ is the principal nilpotent,
to the following classical identity, see e.g. [K1]:
$$ \prod_{1\leq i\leq r}\;
\frac{1- t^{m_i+1}}{1-t}\;=\;
\prod_{\alpha\in R_+}
\frac{1- t^{{\text{height}}(\alpha)+
1}}{1- t^{{\text{height}}(\alpha)}}\;,
$$
where $m_1,\ldots,m_r\,,\, r=\rk\g$ denote the exponents of $\g$.
\bigskip

Further, choose a bi-homogeneous base
$z_1,\ldots,z_r\,,\,
r=\rk\g$, of the space $H_{\!_{{\bold {nw}}}}
:=$ $\bigoplus\limits_{p<0,q\geq 0}\;H_{p,q}\,,$
To each $i= 1,\ldots,\rk\g$, we assign
the pair $(p_i, q_i)\in \Z^2$
such that $z_i\in H_{p_i,q_i}(\e)$. 
The subset ${\op{Exp}}^{\!^{{\bold {nw}}}}_\e(\g)=
\{(p_1, q_1), \ldots, (p_r, q_r)\,,\,
r=\rk\g\}$ of the second quadrant will 
be referrred to as the collection of {\it higher biexponents}
of $\g$ relative to $\e$. Also, write
$R_{\!_{{\bold {nw}}}}\subset R$ for the set of those roots of $(\g,\h)$
that occur in the root-decomposition of $\g_{\!_{{\bold {nw}}}}=
\bigoplus_{p<0,q\geq 0}\;\g_{_{p,q}}$. One can similarly prove
the following identity:
$$\prod_{(p,q)\in\text{Exp}^{\!^{{\bold {nw}}}}_\e(\g)}
\frac{1- s^{p+1}t^{q+1}}{1-s\cdot t}\;=\;
\prod_{\alpha\in R_{\!_{{\bold {nw}}}}}
\frac{1- s^{\langle\alpha , h_1\rangle+1}t^{\langle\alpha , h_2\rangle+
1}}{1-
 s^{\langle\alpha , h_1\rangle+1}t^{\langle\alpha , h_2\rangle}}
\cdot\frac{
1- s^{\langle\alpha , h_1\rangle}t^{\langle\alpha , h_2\rangle}}{1-
 s^{\langle\alpha , h_1\rangle}t^{\langle\alpha , h_2\rangle+1}}\;.$$

\medskip
\bigskip

\head{7. Partial slices.}
\endhead
\bigskip

Given a set $\Sigma \subset \Z^2$ and 
a $\Z^2$-graded vector space  ${\frak a}=\bigoplus_{p,q}\;
{\frak a}_{_{p,q}}
$, we will use the notation ${\frak a}_{_\Sigma}$
for $\bigoplus_{p,q\in \Sigma}\;{\frak a}_{_{p,q}}$; for example,
${\frak a}_{_{{\{p\leq -1\}}}} = \bigoplus_{p\leq -1,\,q\in \Z}\;
{\frak a}_{_{p,q}}\,.$

Fix a principal nilpotent pair
$\e=(e_1,  e_2)$, an
associated semisimple pair $\bh=(h_1,h_2)$, and let
$\g=\bigoplus\;\g_{_{p,q}}$ be
the corresponding bi-grading.
We put
$\pmb{\bgm}:= \g_{_{\{p\leq 0\}}} \oplus \g_{_{\{q\leq 0\}}}
\subset \g\oplus\g
\,$, and consider
the affine subspace $\e+\pmb{\bgm} \subset \g\oplus\g$. 

The following result and its proof are double-analogues of
[K2, Lemma 10].

\proclaim{Lemma 7.1}  The set $\ZZ\cap (\e+\bgm)$ consists
of regular pairs.
\endproclaim

{\it Proof.} Let $\x\in \ZZ\cap (\e+\bgm)$.
The bi-grading on $\g$ gives rise to
 an increasing bi-filtration:
$\g_{_{\leq p,q}} :=\g_{_{\{k\leq p,\,l\leq q\}}}$. 
The  bi-filtration on $\g$ induces a similar bi-filtration
$\z(\x)_{_{\leq p,q}}:=  \z(\x) \cap \g_{_{\leq p,q}}$
on the Lie subalgebra $\z(\x)$. Write $\gr\,\g$ and $\gr\,\z(\x)$
for the corresponding associated bi-graded spaces.
It is clear that the 
natural projection: $\g_{_{\leq p,q}} \twoheadrightarrow
\g_{_{p,q}}$ yields an isomorphism
$\gr_{_{p,q}}\g \iso\g_{_{p,q}}\,$, hence, induces an imbedding
$\sigma: \gr_{_{p,q}}\z(\x)\,\hookrightarrow\,\g_{_{p,q}}\,$.

We claim that $\sigma\bigl(\gr_{_{p,q}}\z(\x)\bigr)\subset
\z_{_{p,q}}(\e)$. To see this, for $i=1,2$,
write $x_i=e_i+a_i$, where
$a_1\in \g_{_{\{p\leq 0\}}}$ and $a_2\in \g_{_{\{q\leq 0\}}}$.
Choose $y\in \z(\x)_{_{\leq p,q}}$,
so that $[x_i,y]=0\,,\,i=1,2$. Write 
$y=\sum_{k\leq p, l\leq
q}\;y_{_{k,l}}$,
and $[x_i,y]=\sum_{k,l}\;[x_i, y]_{_{k,l}}\,,$ for the
corresponding decompositions
into graded components: $y_{_{k,l}}\,,\,[x,y]_{_{k,l}}
\in \g_{_{k,l}}\,$. Then we find:
$$
0=[x_1,y]_{_{p+1,q}}\;=\;[e_1+a_1\,,\,y_{_{p,q}}]_{_{p+1,q}}\;=\;
[e_1\,,\,y_{_{p,q}}]\,.$$
Similarly, computing the $(p,q+1)$-th component of 
$[x_2,y]$, we deduce: $\,0=[x_2,y]_{_{p,q+1}}$
$=[e_2\,,\,y_{_{p,q}}]\,.$
Thus $\sigma(y)=y_{_{p,q}}\in\z(\e)\,$, and our claim is proved.

The claim implies that the map $\sigma$ gives an imbedding:
$\gr\,\,\z(\x)\,\hookrightarrow\,\z(\e)\,$. It follows that
$\dim\z(\x)=\dim\bigl(\gr\,\,\z(\x)\bigr)\leq$
$\dim\z(\e)=\rk\g$.
On the other hand, since $\x\in \ZZ$, the Richardson inequality
yields  $\dim\z(\x) \geq \rk\g$, and the lemma follows.
\qed\medskip

The argument above yields the following result

\proclaim{Corollary 7.2} For any $\x\in \ZZ\cap (\e+\bgm)$, we have
a natural isomorphism: $\gr\,\,\z(\x)\,\iso\,\z(\e)\,$. In particular,
$\z(\x)\cap \g_{_{\{p\leq 0\;\text{or}\;q\leq 0\}}}
= 0.$
\endproclaim

{\it Proof.} The first claim is clear. 
To prove  the second, fix $y \in \z(\x)\cap 
\g_{_{\{p\leq 0\;\text{or}\;q\leq 0\}}}$.
By our choice of $y$,
 there exist a pair of integers $p,q$ such that at least one of them
is non-positive and such that
 $\sigma(y)\neq 0$, where
$\sigma: \gr_{_{p,q}}\z(\x)\,\to\,\g_{_{p,q}}\,$ is the symbol-map.
But we know that $\sigma\bigl(\gr_{_{p,q}}\z(\x)\bigr)
\subset \z(\e)$. Thus, $\sigma(y)\in \z(\e)$, which contradicts
the "positive quadrant" property: $\z(\e)\cap
\g_{_{\{p\leq 0\}}} \cap\g_{_{\{q\leq 0\}}}=0$.\qed\bigskip

We come to the main point of this section.
Recall the notation of Theorem 6.4.
For each $p,q$ such that $p<0\,\, \&\,\, q\geq 0$, choose
an arbitrary subspace $S_{p,q} \subset \z_{_{p,q-1}}(e_2)$
complementary to $\text{\it Image}\left(\ad e_1:\right.$
$\left. \z_{_{p-1,q-1}}(e_2) \to
\z_{_{p,q-1}}(e_2)\right)$, and 
form the subspace $S_{\!_{{\bold {nw}}}} :=$ $
\bigoplus_{p <0, q\geq 0}\;S_{p,q} \subset \g$.
Let $\dis{}_{_{}}\bs_{\!_{{\bold {nw}}}}$
$ := S_{\!_{{\bold
{nw}}}} \bigoplus\, \{0\}$ denote the corresponding subspace
in $\g\oplus\g$, and
define $\bs_{\!_{{\bold {se}}}}=\{0\} \,\bigoplus\,
S_{\!_{{\bold {se}}}} \subset \g\oplus\g$ similarly.
Observe that
the affine
linear spaces: $\dis{}^{^{}}$ $ \,\e + 
\bs_{\!_{{\bold {nw}}}}\,,\,\e + 
\bs_{\!_{{\bold {se}}}}\,\subset\g\oplus\g$ are both
contained in $\ZZ$, the commuting variety.
These affine spaces  play the role
of "{\it partial slices}" to the orbit $\Ad G(\e)$ in $\ZZ$.
For example,
 if $\langle e,\,h,\,f\rangle$ is an ${\frak s}{\frak
l}_2$-triple associated to the regular nilpotent $e$, then for
the principal nilpotent pair $\e=(e,0)$ we have:
$$\e + 
\bs_{\!_{{\bold {nw}}}} = \bigl(e\,, \,\z_\g(e)\bigr)
\quad\text{and}\quad
\e + 
\bs_{\!_{{\bold {se}}}} = \bigl(e + \z_\g(f)\,,\, 0\bigr)\;.$$

Further, assume $\g=\sln$, and let $\e=\e_\lambda$ be the principal
nilpotent pair associated to a Young diagram $\lambda$. Then the
spaces $S_{\!_{{\bold {nw}}}}$ and 
$S_{\!_{{\bold {se}}}}$ can be described as follows.
Given a box $(p,q)$ of the diagram $\lambda$,
let $(p, q_{_{\bold {max}}})\in\lambda $ denote the top box in the same column as
$(p,q)$, and $(p_{_{\bold {max}}},q)\in\lambda $ denote the
rightmost box in the same row as
$(p,q)$. In the notation of Claim 5.8 set $\bnu(p,q):=(\nu_{_{{\bold {in}}}}, 
\nu_{\!_{{\bold {out}}}})$, where
$$
\nu_{_{{\bold {in}}}}=\{(0, q_{_{\bold {max}}}),\ldots,
(p, q_{_{\bold {max}}})\}\;,\;
\nu_{\!_{{\bold {out}}}}=\{(p_{_{\bold {max}}}\!\!-p,q),\,
(p_{_{\bold {max}}}\!\!-p+1,q),\ldots,(p_{_{\bold {max}}},q)\}\,.$$
\smallskip\noindent
Then, one verifies that $n-1$ matrices 
$\{f_{\bnu(p,q)}\;|\;(p,q)\in \lambda\,,\,(p,q) \neq 
(0,q_{_{\text{top}}})\}$
form a basis of $S_{\!_{{\bold {se}}}}$.
Here $(0,q_{_{\text{top}}})$ stand for the coordinates
of the top box in the leftmost column of $\lambda$;
the operator $f_{\bnu(0,q_{_{\text{top}}})}$ corresponding
to this box is a rank one diagonal matrix with trace 1,
which is therefore not in $\sln$.
Fliping the roles of rows and columns of $\lambda$ one similarly
obtains a basis of $S_{\!_{{\bold {nw}}}}$.
\bigskip

We  will now show
that, for an arbitrary semisimple Lie algebra $\g$,
the "partial slices"
have quite
remarkable properties, similar to the properties
of the standard transversal slice to the regular nilpotent orbit
$\Ad G(e)\subset\g$, established by Kostant  [K1]-[K3].
We restrict our attention to the north-west quadrant and
the slice $\bs_{\!_{{\bold {nw}}}}$, the situation
with $\bs_{\!_{{\bold {se}}}}$ being entirely similar.

We remark first that  every pair
$\x\in \e\, + \bs_{\!_{{\bold {nw}}}}$ is regular, due to Lemma 7.1.
Furthermore, there is a $\C^*$-action on
$\e\, + \bs_{\!_{{\bold {nw}}}}$ with a single fixed point,
$\e$, that contracts the space $\e\, + \bs_{\!_{{\bold {nw}}}}$
to this point. Specifically, let $\gamma:\C^*\to G$ be
the homomorphism such that 
$\,\frac{d\gamma}{dt}\big|_{t=1}=h_1\,,$
see Corollary 3.6(i). Define a $\C^*$-action on $\g\oplus\g$
by the formula: $\C^*\ni t: (x_1,\,x_2) \mapsto 
\left(t\cdot\Ad \gamma(t^{-1})\,x_1\,,\,
\Ad \gamma(t^{-1})\,x_2\right)$. It is clear that
this $\C^*$-action preserves the commuting
variety $\ZZ$, and takes any $\Ad G$-diagonal orbit in
$\g\oplus\g$ into another $\Ad G$-diagonal orbit.
Furthermore, it
keeps the point $\e$ fixed, and
contracts the space $\g_{_{\{p\leq 0\}}}\oplus\{e_2\}$
to $(0, e_2)$, hence,  contracts $\e\, + \bs_{\!_{{\bold {nw}}}}$
to $\e$.\medskip

We introduce the nilpotent Lie subalgebra
 $\nnw:= \g_{_{{\{p\leq -1\,,\,
q\geq 0\}}}}\, \subset \,\g$, and write
$N\nw$ for the corresponding unipotent group. Further,  set
$$ \gnw = \gnww1\bigoplus\gnww2\;
\subset\;\g\bigoplus\g\quad,\quad
\gnww1:=\g_{_{{\{{p\leq 0\,,\,
q\geq 0}\}}}}\enspace,\enspace\gnww2:=
\g_{_{{\{{p\leq -1\,,\,
q\geq 1}\}}}}\,.\tag 7.3$$

We consider the affine subspace $\e+\gnw \subset\;\g\oplus\g$,
and observe that it is stable under $\Ad N\nw$-diagonal action.
Let $\bigl(\ZZ\cap (\e+\gnw)\bigr)_{\text{red}}$ denote
the intersection $\ZZ\cap (\e+\gnw)$ as an algebraic variety with
reduced scheme structure.
The following result is a double analogue of [K3, Theorem 1.2].
Note that  even in the classical case
 the proof given below is  simpler than the argument
in [K3].

\proclaim{Theorem 7.4} The action-map gives an isomorphism
of algebraic varieties: 
$$
 N\nw\times (\e\, + \bs_{\!_{{\bold {nw}}}})
\,\iso\, \bigl(\ZZ\cap (\e+\gnw)\bigr)_{\text{red}} \,.
$$
\endproclaim

{\it Proof.} 
We know, by Corollary 7.2,
 that for any point in $\ZZ\cap (\e+\gnw)$, the Lie algebra
of isotropy group of the $\Ad N\nw$-action is trivial.
It follows, since the group $N\nw$ is unipotent,
that  $\Ad N\nw$-diagonal action on
$\ZZ\cap (\e+\gnw)$ is free. We claim further that
the  action-map: $N\nw\times (\e\, + \bs_{\!_{{\bold {nw}}}})
\,\to\, \left(\ZZ\cap (\e+\gnw)\right)_{\text{red}}$ 
is an injective open morphism.
To prove injectivity it suffices to show that 
each $\Ad N\nw$-orbit meets the space $\e+\bs\nw$ in at most one
point. This is clear locally near $\e$ since the
space $\e+\bs\nw$ was chosen to be a transverse slice to the
$\Ad N\nw$-orbit of $\e$. Hence, the same holds
globally on $\e+\gnw$, because the  $\C^*$-action on
$\g\oplus\g$ defined two paragraphs before the theorem
preserves the variety $\ZZ\cap(\e+\gnw)$, takes
$\Ad N\nw$-diagonal
orbits into $\Ad N\nw$-diagonal
orbits and, moreover,
contracts the partial slice $\e+\bs\nw$  to the point $\e$.
Finally, by definition of a transverse slice
(see the geometric meaning of the space $H^1({\frak e},\,\g)$
explained in \S6)
the morphism: $N\nw\times (\e\, + \bs_{\!_{{\bold {nw}}}})
\,\to\, \ZZ\cap (\e+\gnw)$ has a surjective
differential at the point $(\pmb{1}_{_{N\nw}},\,\e)$,
hence is an open  morphism. 

Next, consider the affine subspace
$e_2 + \gnww2=e_2+
\g_{_{{\{{p\leq -1\,,\,q\geq 1}\}}}}\subset\g$.
 It is clear from (7.3) that:
$[\nnw\,,\,\gnww2]\subset \gnww2$, and that
$[\nnw,\,e_2]\subset \gnww2$. It follows that the set
$e_2 + \gnww2 \subset\g$ is
 $\Ad N\nw$-stable, see e.g.
[CG, Lemma 1.4.12]. Furthermore, the last inclusion is
actually an equality, due to the weak Lefschetz.
Hence, the $\Ad N\nw$-orbit of $e_2$ is Zariski open in
$e_2 + \gnww2$. The group $N\nw$ being unipotent,
the orbit has to be closed, and we conclude that
$e_2 + \gnww2=\Ad N\nw(e_2)$. Thus, any point
of $\ZZ\cap (\e+\gnw)$ is $\Ad N\nw$-conjugate to a point
of the form $(e_1+x,\, e_2)$, where $x\in\z(e_2)_{_{\{{p\leq 0\,,\,
q\geq 0}\}}}$.
Note that the condition $q\geq 0$ is superfluous, since
$\z(e_2)_{_{p,q}}=0$ for all $q<0$ by the weak Lefschetz.

Consider the Lie algebra 
$\znw(e_2):= \z(e_2)_{_{\{{p\leq -1}\}}} \subset \,\nnw$, and write
$Z\nw(e_2)$ for the corresponding unipotent group.
We have shown in the previous paragraph
that the second projection: $(x_1, x_2)\mapsto x_2$ induces, set theoretically,
an  isomorphism of $N\nw$-equivariant
fibrations:
$$\Bigl(\ZZ\cap (\e+\gnw)\, \to\Ad N\nw(e_2)\Bigr)
\;\simeq\;
\Bigl(N\nw \times_{_{_{Z_{_{\bold {nw}}}(e_2)}}} (e_1\,+ 
\z(e_2)_{_{\!\{p\leq 0\}}})\,\to N\nw/Z_{_{\bold {nw}}}(e_2)\Bigr)\,.
\tag 7.6$$
The space on the RHS here is an affine bundle over the base
$N\nw/Z\nw(e_2)$, which is a smooth connected affine variety.
It follows that $\left(\ZZ\cap (\e+\gnw)\right)_{\text{red}}$
is itself  a smooth connected affine variety. 
Hence, proving the theorem amounts to showing
that the action-map: $\varphi:
N\nw\times (\e\, + \bs_{\!_{{\bold {nw}}}})
\,\to\, \ZZ\cap (\e+\gnw)$ is bijective. 
We already know, 
by the first paragraph of the proof, that $\varphi$ is an injective
morphism with Zariski open image. It is  clear that 
the algebraic variety
$N\nw\times (\e\, + \bs_{\!_{{\bold {nw}}}})$ is isomorphic
 to a vector space $\C^k$. 
The group $N\nw$ being unipotent, the variety
$N\nw/Z_{_{\bold {nw}}}(e_2)$ is also isomorphic
 to a vector space $\C^k$. Hence, $\ZZ\cap (\e+\gnw)$,
the total space of
the fibration (7.6), is isomorphic topologically
(in fact algebraically) to $\C^k$.
Thus, we are reduced to proving the following:\par
{\it Let $\varphi: U \hookrightarrow X$
be a
Zariski open imbedding of irreducible affine 
algebraic varieties
such that both $U$ and $X$ are topologically isomorphic to
$\C^k$. Then $U=X$.}

To prove this claim  we use the standard
long exact sequence of Borel-Moore homology, see e.g.
[CG, ch.2]:
$$\ldots\to H^{^{\text{BM}}}_{i+1}(X) 
\,{\overset {\varphi^*} \to\longrightarrow}\, H^{^{\text{BM}}}_{i+1}(U)
\to H^{^{\text{BM}}}_i(X\smallsetminus U) \to H^{^{\text{BM}}}_i(X) 
\,{\overset {\varphi^*} \to\longrightarrow}\,
H^{^{\text{BM}}}_i(U)\to\ldots
$$
By assumption, the restriction map $\varphi^*$ here is
an isomorphism between two spaces of dimension
1 if $i=2k$, and of dimension zero otherwise. It follows from
the exact sequence that {\it all} Borel-Moore homology
groups of $X\smallsetminus U$ vanish. But this
 contradicts the
fact that the algebraic variety $X\smallsetminus U$
has a non-zero fundamental class in 
$H_{_{\text{top}}}^{^{\text{BM}}}(X\smallsetminus U)$,
unless $X\smallsetminus U=\emptyset$.\qed
\bigskip

Next, we study the south-west quadrant, which turs out to be
much simpler. \smallskip

Introduce the Lie subalgebra $\nsw := 
\g_{_{\{p\le -1\,,\,q\le -1\}}}$, and let
$N\sw\subset G$ denote the corresponding (connected) unipotent
 subgroup. We also put 
$\pmb{\frak g}_{_{{\bold {sw}}}}:=
\g_{_{\{p\le 0\,,\,q\le 0\}}}\oplus \g_{_{\{p\le 0\,,\,q\le 0\}}}
\subset \g\oplus\g$.

\proclaim{Proposition 7.7} The group $N\sw$ acts freely on $\ZZ\cap
(\e+\pmb{\frak g}_{_{{\bold {sw}}}})$,
and  $\ZZ\cap (\e+\pmb{\frak g}_{_{{\bold {sw}}}}) = \Ad N\sw(\e)$
is a single $\Ad N\sw$-diagonal orbit.
\endproclaim

{\it Proof.} Since the group $N\sw$ is unipotent,
to prove the first claim it suffices to show
that, for any $\x\in \ZZ\cap (\e+\pmb{\frak g}_{_{{\bold {sw}}}})$, one has:
$\nsw\cap \z(\x)=0$. This follows from Corollary 7.2.

To prove the second claim, 
we introduce a bi-filtration
on $\pmb{\frak g}_{_{{\bold {sw}}}}$ as follows:
$$\pmb{\g}_{_{\leq p,q}} := \{(a_1, a_2)\enspace|\enspace
a_1\in \g_{_{\{i\leq p\}}}\;,\;
a_2\in \g_{_{\{j\leq q\}}}\}\,.$$
Let $\sigma: \pmb{\g}_{_{\leq p,q}} \twoheadrightarrow
\g_{_{p,*}}\oplus\g_{_{*,q}}$ denote the corresponding
"symbol-map". We will prove by induction on $(p,q)$,
where $\Z\oplus\Z$ is viewed as a partially ordered set,
that for any ${\bold a}\in \pmb{\g}_{_{\leq p,q}}$,
such that $
\x=\e+{\bold a}\in \ZZ\cap (\e+\pmb{\frak g}_{_{{\bold {sw}}}})$,
the element $\x$ is
$\Ad N\sw$-conjugate to $\e$.

Fix $\x=\e+{\bold a}\in \ZZ\cap (\e+\pmb{\frak g}_{_{{\bold {sw}}}})$,
where $x_i=e_i+a_i$, and
${\bold a}\in \pmb{\g}_{_{\leq p,q}}$.
We may assume that
$\sigma({\bold a})\neq 0$. 
Then the leading term of the equation $[x_1, x_2]=0$ reads:
$$
0\;=\;[e_1+a_1\,,\,e_2+a_2]\;=\;[e_1, a_2] + [e_2, a_1]\,.$$
The equation shows that $(a_2, a_1)\in \g\oplus\g$ is a cocycle
giving a class in $H^1({\frak e}, \g)$, see  (6.1). Since,
${\bold a}\in \pmb{\frak g}_{_{{\bold {sw}}}}$, both $p$ and $q$
are non-positive. Hence, by Theorem 6.4 the corresponding cohomology
group vanishes, so that there exists $y\in\g$ such that
$a_i= [e_i, y]\,,\, i=1,2.$ Note that
since ${\bold a}\in \pmb{\frak g}_{_{{\bold {sw}}}}$, the 
equations $a_i= [e_i, y]$
imply that $y\in \g_{\{p\leq -1, q\leq -1\}} = \nsw$.
Therefore, we get:
$\Ad \exp(-y)\,(\e)= \e+{\bold a} + \boldsymbol{\epsilon},$ where 
$\boldsymbol{\epsilon}$
belongs to lower terms of the bi-filtration on $\pmb{\frak g}_{_{{\bold {sw}}}}$.
Hence, $(\Ad \exp y)(\boldsymbol{\epsilon})$ 
belongs to lower terms of the bi-filtration again, and
the induction hypothesis implies: $\e- (\Ad \exp
y)(\boldsymbol{\epsilon})
\in
\Ad N\sw(\e)$. Thus we find: 
$$
\align
(\Ad \exp y)(\x) & =(\Ad \exp y)(\e+{\bold a})\\
& =(\Ad \exp y)\bigl(\Ad \exp(-y)\,(\e)-
\boldsymbol{\epsilon}\bigr)=\e - (\Ad \exp y)(\boldsymbol{\epsilon})
\;\in\; \Ad N\sw(\e),
\endalign
$$
and the claim follows.\qed\bigskip

\newpage


\centerline{{\bf 
{8. Appendix: TOWARDS A CLASSIFICATION OF}}\qquad\hphantom{x}}
\medskip

\centerline{\bf {PRINCIPAL
NILPOTENT PAIRS.}}

\bigskip
\centerline{
A.~Elashvili\footnote"*"{A.E.: Razmadze Math. Institute, M.Aleksidze 1,
Tbilisi 380093, Republic of Georgia;\newline
$\hphantom{x}\;\quad$
E-mail: {\bf alela\@rmi.acnet.ge}\quad (Work 
supported in part by Grant INTAS-OPEN-97-1570.)}
 and D.~Panyushev\footnote"**"{D.P.: 
MIREA, Math. Dept., prosp. Vernadskogo 78,
Moscow 117454, Russia;
\newline$\hphantom{x}\;\quad$
E-mail: {\bf panyush\@dpa.msk.ru}\quad (Work 
supported by RFFI grant 98-01-00598.)}}
\bigskip

In connection with the preceding article, V.Ginzburg asked us
whether there existed non-trivial examples of principal nilpotent
pairs in exceptional Lie algebras. We present here a full 
description of principal nilpotent pairs (= {\it pn}-pairs)
in the exceptional case and some results 
towards a complete classification in classical Lie algebras.

As has been observed in section~1 of the main body of the paper, 
equalities like $\dim\g_{p,q}=\dim\g_{-p,q}$ and
$\dim\g_{p,q}=\dim\g_{p,-q}$ are typically
false for the $\Z^2$-grading associated with a
{\it pn}-pair. We give first a  classification of an interesting
class of {\it pn}-pairs satisfying these equalities. 

\medskip

\noindent
{\bf Definition. } A {\it pn}-pair $\e=(e_1,e_2) \in {\g}\times {\g}$
is said to be {\it rectangular} if $e_1$ and $e_2$ can be included in
commuting $\tri$-triples $\{e_1,h_1,f_1\}$ and $\{e_2,h_2,f_2\}$.

\medskip

The name is explained by the following observation. For $\g=\sel{N}$, a
{\it pn}-pair $\e$ is
rectangular if and only if the Young diagram of $e_1$ (or $e_2$)
is a {\it rectangle}, see below.
In particular, non-trivial rectangular pairs
exist if and only if $N$ is not a prime.

\medskip

\noindent{\bf Remarks.} 1.
It is easy to see that a {\it pn}-pair $\e$ is rectangular if and only if
there exists an $\tri$-triple $\{e_1,h_1,f_1\}$ such that $[e_2,f_1]=0$.
\newline 2. We assume that $[h_i,e_i]=2e_i$, hence
the pair $(h_1,h_2)$ in the rectangular case is twice the
associated semisimple pair, in the sense of section~1.
\newline 3. Note that if
$\{e_1,h_1,f_1\}$ and $\{e_2,h_2,f_2\}$ are commuting $\tri$-triples,
then condition (ii) of the Definition above is automatically satisfied. That is,
such a pair $(e_1,e_2)$ is principal if and only if $\dim\z(\e)=\rk\g$.

\medskip

A general description of rectangular {\it pn}-pairs can quickly
be obtained
without using structure theory of {\it pn}-pairs developed in sections 1 and 2.
To this end, we briefly recall the structure of the centralizer $\z(e)$ of
a nilpotent element $e\in\g$ (the Dynkin-Kostant theory, see e.g. [CM, ch.4]).
Let $\{e,h,f\}$ be an $\tri$-triple
and $\g=\bigoplus_{i\in\Z}\g_i$ the corresponding $\Z$-grading. Then
$\z(e)=\bigoplus_{i\ge 0}\z(e)_i$ and $\z(e)_0$ is a maximal reductive
subalgebra in $\z(e)$. Moreover, $\z(e)_0=\z(e,f)=\z(e,h,f)$.
One has $\z(e)\simeq\g_0\oplus\g_1$ as $\z(e)_0$-module.
The element $e$ is called {\it even} whenever all the eigenvalues of
$\ad h$ are even, i.e., $\g_i=0$ for $i$ odd. Obviously, $e$ is even
if and only if $\g_1=0$. In this case the weighted Dynkin diagram of
$e$ contains only numbers 0 and 2 [Dy].

\proclaim{Theorem 8.1}
The following conditions are equivalent:
\newline\noindent
(i) $e$ is a member of a rectangular {\it pn}-pair;
\newline\noindent
(ii) $e$ is even, and any regular nilpotent element in $\ka:=\z(e)_0$ is
regular in $\g_0$ as well.
\vskip0.5ex\noindent
Under condition {\rm (ii)}, if $e'$ is any regular nilpotent element in
$\ka$, then $(e,e')$ is a {\it pn}-pair.
\endproclaim

\noindent{\it Proof.} (i)$\Rightarrow$(ii)\quad
For any nilpotent element $e'\in\ka$, we can choose
an $\tri$-triple $\{e',h',f'\}$ lying inside $\ka$.
Application of $\ad f$ yields an isomorphism of
$\ka$-modules $\z(e)_{ev}\simeq\g_0$ and
$\z(e)_{odd}\simeq\g_1$. Therefore
$$
\dim\z(e,e')=\dim\z_{\g_0}(e')+\dim\z_{\g_1}(e')\ge\rk\,\g+
\dim\z_{\g_1}(e')\ .
$$
Thus if $\dim\z(e,e')=\rk\g$, then $\g_1=0$ and $e'$ is regular
in $\g_0$. The latter implies that $e'$ is regular in $\ka$, too.

\medskip

\noindent (ii)$\Rightarrow$(i)\quad Reverse the previous argument.
\hfill $\square$

 As an immediate consequence we obtain the following uniqueness
statement in the rectangular case:
\proclaim{Corollary 8.2}
Given an $\tri$-triple $\{e,h,f\}$ containing the first member of
a rectangular {\it pn}-pair, the second member is determined uniquely up to
conjugacy by an element of the connected group $K=Z_G(e,h,f)^o$.
\endproclaim
 It is likely that a kind of uniqueness holds for arbitrary
principal nilpotent pairs.
\vskip1ex\noindent
Using theorem 8.1, we find the rectangular {\it pn}-pairs in the exceptional
simple Lie algebras. The tables in [El1] contain the information on
$\g_0$ and $\ka$ for all nilpotent orbits. To verify condition (ii) of
theorem 8.1, we need the description
of inclusion $\ka\subset\g_0$. The latter is determined by the structure
of $\g_2$ as $\g_0$-module, since $\ka$ is a generic stabilizer.
An algorithm for describing $\g_0$-representation in $\g_2$ is
given in [El1] as well.
Usually, a nilpotent orbit in exceptional Lie algebra is denoted 
by a Cartan label. This label is said to be the {\it type} of nilpotent orbit.
The idea of such notation goes
back to Dynkin [Dy], who studied minimal regular\footnote"*"{a
subalgebra is called regular if its normalizer in $\g$ contains a
Cartan subalgebra.}  reductive subalgebras
in $\g$ containing a given simple 3-dimensional subalgebra or, what is the
same, a given nilpotent element. The type of
an orbit represents one of these subalgebras, namely, a unique,
up to conjugation, minimal Levi subalgebra of $\g$ containing
an element of the nilpotent orbit under consideration. We refer to
[CM, 8.4] for the tables of nilpotent orbits and the corresponding labels,
where some more explanation concerning this notation is found.

\proclaim{Theorem 8.3}
The following list contains all the rectangular {\it pn}-pairs in
exceptional Lie algebras and their bi-exponents:
\vskip1ex
\halign to\hsize{\indent #\quad & # &\hfil #\hfil & # \cr
$\GR{F}{4}$ &  & $(\GR{B}{3},\GRt{2})$ & $\{(1,0),(5,0),(0,1),(3,2)\}$ \cr
$\GR{E}{6}$ &  &$(\GR{D}{4},2\GR{A}{2})$
              & $\{(1,0),(5,0),(0,1),(3,1),(0,2),(3,2)\}$ \cr
$\GR{E}{7}$ &  &$(\GR{E}{6},[3\GR{A}{1}]'')$ &
              $\{(1,0),(5,0),(7,0),(11,0),(0,1),(4,1),(8,1)\}$\cr
            &  &$(\GR{A}{6},\GR{A}{2}{+}3\GR{A}{1})$ &
              $\{(1,0), (5,0), (0,1), (2,1), (4,1), (6,1), (3,2)\}$\cr
            &  &$(\GR{A}{4}{+}\GR{A}{2},\GR{A}{3}{+}\GR{A}{2}{+}\GR{A}{1})$
              & $\{(1,0),(0,1),(2,1),(4,1),(1,2),(3,2),(2,3)\}$\cr
            &  &$([\GR{A}{5}]'',\GR{D}{4})$ &
             $\{(1,0), (3,0), (5,0), (0,1), (2,3), (4,3), (0,5)\}$\cr
$\GR{E}{8}$ &  &$(\GR{E}{6},\GR{D}{4})$ &
              $\{(1,0),(5,0),(7,0),(11,0),(0,1),(4,3),(8,3),(0,5)\}$\cr}

\endproclaim
\demo{Proof}
We give a sample of our computations. Consider the nilpotent orbit
$\co_1=\Ad G(e_1)$ of type $\GR{A}{2}+3\GR{A}{1}$ in $\g=\GR{E}{7}$.
The weighted
Dynkin diagram of $\co_1$ is \quad \quad
$\Bigl($ \lower2.5pt\hbox{0--0--0--\vbox{\hbox{2}\hbox{0\strut}}--0--0}
$\Bigr)$\ .

 Therefore $[\g_0,\g_0]\simeq \sel{7}=\frak{sl}(\Bbb V)$
and $\sel{7}$-module $\g_2$ is isomorphic to $\wedge^3\Bbb V$.
Then $\ka\simeq\GR{G}{2}$, and the embedding $\GR{G}{2}\hookrightarrow
\sel{7}$ corresponds to the unique 7-dimensional representation of
$\GR{G}{2}$. It is easy to show that the restriction of this
representation to a principal $\tri$ in $\GR{G}{2}$ yields a 7-dimensional
irreducible representation of this $\tri$ (the single Jordan block).
This precisely means that a (any) regular  nilpotent element in
$\GR{G}{2}$ is still regular in $\sel{7}$. It remains to determine the type
of $\co_2=\Ad G (e_2)$, where $e_2$ is a regular nilpotent element in
$\GR{G}{2}$. Since $e_2$ is also regular in $\sel{7}$, we see that
$\g_0$ is a minimal Levi subalgebra intersecting $\co_2$. Thus $\co_2$
is of type $\GR{A}{6}$. Therefore the weighted Dynkin diagram of $\co_2$ is
\quad
$\Bigl($~\lower2.5pt\hbox{0--2--0--\vbox{\hbox{0}\hbox{2\strut}}--0--0}~%
$\Bigr).$
\vskip1ex

In each of these cases, it is not hard to determine the structure of $\g$
as $\tri\oplus\tri$-module and then find the bi-exponents of the
{\it pn}-pair.
Indeed, let $R(d)$ denote the irreducible $\tri$-module of dimension $d+1$
and
let $\g\vert_{_{\tri\oplus\tri}}=\bigoplus_{i=i}^{\rk\g} R(n_i)\otimes
R(m_i)$. Since nilpotent elements in both $\tri$-triples are even,
the integers $n_i, m_i$ are even. Then the bi-exponents are
$\{(n_i/2,m_i/2)\mid i=1,\dots,\rk\g\}$.
\hfill $\square$
\enddemo
\vskip1.5ex
Theorem 8.1 applies to the classical Lie algebras as well, but in this case
it is easier to obtain the classification of
{\it pn}-pairs in another way. Namely, 
we exploit a simple relationship between the simplest and
the adjoint representation of a classical Lie algebra.
Let $\g=\g(\V)$ be a classical Lie algebra, $\V$ being its tautological
representation. Then
\vskip0.5ex
\centerline{ $\g\simeq \V\otimes\V^*\ominus\text{\{triv. 1-dim repr.\}}$
for $\frak{sl}(\V)$,}
\centerline{ $\g\simeq S^2(\V)$ for $\frak{sp}(\V)$,}
\centerline{ $\g\simeq \wedge^2(\V)$ for $\frak{so}(\V)$.}
\vskip0.5ex
Consider the subalgebra $\s={'\tri}\oplus{''\tri}\subset\g(\V)$,
where $'\tri=\langle e_1,h_1,f_1\rangle$ and
$''\tri=\langle e_2,h_2,f_2\rangle$ and the corresponding $\tri$-triples
commute.
We wish to determine those $\s$ that correspond to {\it principal}
nilpotent pairs of rectangular type.
By the very definition, it is equivalent to that $\g(\V)$ as $\s$-module
is a sum of exactly $\rk\g$ irreducible submodules.
Therefore we assume that
$$
 \V\mid_\s=\bigoplus_{i=1}^p R(n_i-1)\otimes R(m_i-1) \tag 8.4
$$
and then compute the decomposition $\g\mid_\s$ making use of the
above 3 relations and variations of the Clebsch-Gordan formula.
Of course, appropriate constraints on parity of $n_i,m_i$
should be satisfied in the orthogonal and symplectic cases.
An advantage of this
approach is that the decomposition $\V\mid_\s$ immediately yields the
description of $e_1,e_2$ in terms of partitions, i.e., if $n_1\ge n_2\ge\dots$,
then $e_1$ corresponds
to the partition $(\underbrace{n_1,\dots,n_1}_{m_1 \text{ times}},
\underbrace{n_2,\dots,n_2}_{m_2 \text{ times}},\dots)$ and likewise for 
$e_2$. Note that $\dim\V=\sum_i n_im_i$.
The output of our computations is as follows:
\proclaim{Theorem 8.5}
Suppose the embedding $\s
\hookrightarrow\g(\V)$ is given by (8.4).
Then $(e_1,e_2)$ is a (rectangular)
{\it pn}-pair in $\g(\V)$ if and only if the set
${\frak I}=\{(n_i,m_i)\mid i=1,\dots,p\}$ is of the form

(i) for $\frak{sl}(\V)$: \quad $p=1$;

(ii) for  $\frak{sp}(\V)$: \quad $p=1$ and $n_1,m_1$ have different parity;

(iii) for  $\frak{so}(\V)$: either $p=1$ and $n_1,m_1$ have the same parity or

\noindent
$p=2$ with ${\frak I}=\{(n,m),(1,1)\}$ or
${\frak I}=\{(n,1),(1,m)\}$, where $n,m$ are both odd.

\endproclaim

This means in particular that
we obtain a single rectangle for $\frak{sl}(\V)$ and $\frak{sp}(\V)$
and  at most 2
rectangles for $\frak{so}(\V)$.
Let us give an explicit presentation of the respective {\it pn}-pair for
the most interesting (last) case.
Let $e_1$
(resp. $e_2$) be a regular nilpotent element in $\frak{so}_n$
(resp. $\frak{so}_m$), where $n,m$ are odd. Consider the embeddings
$\frak{so}_n\times \frak{so}_m\hookrightarrow\frak{so}_{nm}
\hookrightarrow\frak{so}_{nm+1}$ (the first one corresponds
to the tensor product and the second one is natural)
and $\frak{so}_n\times \frak{so}_m\hookrightarrow\frak{so}_{n+m}$
(the direct sum).
Then the image of $(e_1,e_2)$ under these embeddings is a
rectangular {\it pn}-pair.

\vskip2ex
As for the classification of arbitrary {\it pn}-pairs, one may use
the following approach.
Results in section 3 concerning `associated Levi subalgebras'
shows that if $\e$ is a {\it pn}-pair, then the respective pair of
nilpotent orbits $(\co_1,\co_2)$ is reciprocal, see Definition 4.10.
So, the initial step is to describe all reciprocal pairs in the
simple Lie algebras.
Making use of the explicit description of the induced orbits
given in [El2], one easily finds the pairs of
reciprocal orbits in the exceptional Lie algebras. 
Of course, the pairs of orbits indicated in Theorem 8.3 are reciprocal.
The total
number of non-trivial pairs of reciprocal orbits in $\GR{G}{2}$,
$\GR{F}{4}$, $\GR{E}{6}$, $\GR{E}{7}$,
$\GR{E}{8}$ is 0,\,2,\,4,\,7,\,5 respectively. For instance,
the five reciprocal pairs in $\GR{E}{8}$ correspond to the
following 5 pairs of Levi subalgebras:
$$(\GR{E}{6},\,\GR{D}{4})\,,\,
(\GR{A}{6},\,\GR{D}{4}{\bold{+}}\GR{A}{2})\,,\,
(\GR{D}{5},\,\GR{D}{5})\,,\,
(\GR{A}{6}{\bold{+}}\GR{A}{1},\,
\GR{A}{4}{\bold{+}}\GR{A}{2}{\bold{+}}\GR{A}{1})\,,\,
(\GR{D}{5}{\bold{+}}\GR{A}{2},\,
\GR{A}{4}{\bold{+}}\GR{A}{2})\,.$$
 
Below we give the complete description
of {\it pn}-pairs in the exceptional simple Lie algebras.

\proclaim{Theorem 8.6}  \newline
1.  There are no {\it pn}-pairs in
$\GR{G}{2}$.
The rectangular {\it pn}-pair of theorem 8.3 is the unique, up to conjugacy,
{\it pn}-pair in $\GR{F}{4}$; \newline
2. There are four $\Ad G$-orbits of {\it pn}-pairs in
$\GR{E}{6}$. The three orbits of non-rectangular {\it pn}-pairs
correspond to the following reciprocal pairs of orbits, where the last column
gives the bi-exponents:
\vskip1ex
\halign to\hsize{\indent #\quad & # &\hfil # \cr
$(\GR{D}{5},2\GR{A}{1})$ & -- &
$\{(1,0),(4,0),(5,0),(7,0),(0,1),(3,1)\}$ \cr
$(\GR{A}{4}{+}\GR{A}{1},\GR{A}{2}{+}2\GR{A}{1})$ & -- &
$\{(1,0),(4,0),(0,1),(2,1),(3,1),(1,2)\}$ \cr
$(\GR{A}{4},\GR{A}{3})$ & -- &
$\{(1,0),(3,0),(4,0),(0,1),(0,3),(2,2)\}$ \ . \cr}

\noindent
3. There are five $\Ad G$-orbits of {\it pn}-pairs in 
$\GR{E}{7}$. The unique orbit of non-rectangular {\it pn}-pairs
corresponds to the following reciprocal pairs of orbits:
\halign to\hsize{\indent #\quad & # &\hfil # \cr
 $(\GR{A}{4}{+}\GR{A}{1},\GR{A}{4}{+}\GR{A}{1})$ & -- &
$\{(1,0),(0,1),(0,4),(1,3),(2,2),(3,1),(4,0)\}$\ . \cr}

\noindent 4. The rectangular {\it pn}-pair of theorem 8.3 is the 
unique, up to conjugation, {\it pn}-pair in $\GR{E}{8}$.
\endproclaim
\demo{Proof} All the proofs are based on explicit calculations with
centralizers. Having a suitable candidate $e_1$ for the first member of
{\it pn}-pair, we try to select an $e_2\in\z(e_1)$ in order to meet all
the requirements of Definition 1.1.
To establish whether some commutators in $\z(e_1)$ are equal or
not equal
to zero, we have used computer program GAP [GP]. \par
1. For $\GR{F}{4}$, we need only to demonstrate that the reciprocal pair of 
orbits $(\GR{C}{3}, \GR{A}{2})$ do not produce a {\it pn}-pair. If $e_1$ lies
in the orbit of type $\GR{C}{3}$, then $\dim\z(e_1)=10$ and $\z(e_1)_0$ is
a 1-dimensional toral subalgebra. Therefore a nonzero $h\in\z(e_1)_0$ is, up
to conjugation, the only possible candidate for $h_2$, the second member of
an associated semisimple pair. However, explicit computations show
that $\dim\z(e_1,e_2)\ge 6$ for any $\ad h$-weight vector
$e_2\in\z(e_1)_{>0}$. \par
2.  For $\GR{E}{6}$, we show that each reciprocal pair of orbits yields
a unique {\it pn}-pair. For instance, consider an element $e_1$ in the orbit 
of type $\GR{D}{5}$.
Then $\dim\z(e_1)=10$ and  $\dim\z(e_1)_i$ is equal to
1,\,1,\,2,\,1,\,1,\,3,\,0,\,1 for $i= 0,\,2,\,4,\,6,\,8,\,10,\,12,\,14$
respectively, $\g_{14}$ being the greatest nonzero subspace in the 
$\Bbb Z$-grading defined by $\ad\,h_1$. 
These data are easily derived from the weighted Dynkin diagram of $e_1$.
Therefore $\z(e_1)_0$ is a 1-dimensional toral algebra and
$\z(e_1)_2=\langle e_1\rangle$. The weights of $\z(e_1)_0$ in
$\z(e_1)_4$ are nonzero (and opposite) and we take either of weight
vectors as $e_2$. Then a straightforward computation shows that
$\dim\z(e_1,e_2)_i$ is equal to 0,\,1,\,1,\,0,\,1,\,2,\,0,\,1 for 
$i=0,\,2,\,4,\,6,\,8,\,10,\,12,\,14$ respectively.
Let $t$ be a nonzero element in $\z(e_1)_0$ normalized 
so that $[t,e_2]=e_2$. Then the elements $\tilde h_1=h_1/2-2t$ and
$\tilde h_2=t$ satisfy the commutator relations $[\tilde h_i,e_j]=\delta_{i,j}
e_j$ and $[\tilde h_1,\tilde h_2]=0$. Hence $(\tilde h_1,\tilde h_2)$
is an associated semisimple pair in the sense of section 1 and $(e_1,e_2)$
is a {\it pn}-pair.
Since the eigenvalues of $\ad\,h_1$ are known, it only suffices to compute the 
eigenvalues of $\ad\,t$ on $\z(e_1,e_2)$ in order to determine the 
bi-exponents in this case. \par
3,4. Because calculations for $\GR{E}{7},\GR{E}{8}$ are not
illuminating, too, we omit them.
\hfill $\square$
\enddemo
{\bf Remark.} It could have happened, a priori, that some
reciprocal pair $(\co_1,\co_2)$  would give rise to {\it several} $\Ad G$-orbits 
of {\it pn}-pairs. We see, as a result of
straightforward calculations, that this never happens.
\vskip1ex

Finally, we give an example of non-rectangular {\it pn}-pair in a classical
Lie algebra. Let $\g=\frak{so}(\V)=\GR{D}{2n+1}$. Take the reciprocal
pair of orbits corresponding to the partitions
$(2n+1,2n+1)$ and $(2,\dots,2,1,1)$. This pair of orbits really produces a
{\it pn}-pair. An explicit matrix presentation of $(e_1,e_2)$ is the
following. Let $v_1,\dots,v_{4n+2}$ be a basis of $\V$ such that the 
$\g$-invariant quadratic form is $x_1x_{4n+2}+\dots+x_{2n+1}x_{2n+2}$. Then
$e_1$ acts on the basis vectors by $e_1(v_j)=v_{j-1}$ ($j\ge 2n+3$),
$e_1(v_j)=-v_{j-1}$ ($2{\le} j{\le} 2n+1$),
$e_1(v_1)=e_1(v_{2n+2})=0$;
$e_2$ acts by $e_2(v_j)=(-1)^{j}v_{j-2n-2}$ ($j\ge 2n+3$),
$e_2(v_j)=0$ $(j\le 2n+2)$.
The simultaneous centralizer $\z(e_1,e_2)$ has the basis consisting of
the following matrices:
$e_1,\,e_1^3,\dots,e_1^{2n-1},\,e_2,\,e_1^2e_2,\dots,e_1^{2n-2}e_2,x$.
Here $x$ is the operator sending $v_{4n+2}$ to $v_{2n+2}$ and
$v_{2n+1}$ to $-v_1$. We leave it to the reader as an exercise to write
down an associated semisimple pair and then the bi-exponents in this case.
\bigskip


\bigskip


\medskip\bigskip

\Refs\nofrills{\bf References}
\widestnumber\key{100}


\ref
\key AL
\by D. Alvis, G. Lusztig
\paper {\it {On Springer's correspondence for simple groups of type}}
 $E\sb{n}$ $(n=6,\,7,\,8)$. ({\sl With an appendix by N. Spaltenstein.})
\paperinfo Math. Proc. Cambridge Philos. Soc. {\bf {92}}
(1982), no. 1
\pages  65--78
\endref

\ref
\key BV
\by D. Barbasch, D. Vogan
\paper {\it Unipotent representations of complex semisimple groups}
\paperinfo  Ann. of Math.  {\bf {121}} (1985), no. 1
\pages 41--110
\endref

\ref
\key BM
\by W. Borho, R. MacPherson
\paper {\it Partial resolutions of the nilpotent variety}
\paperinfo  Asterisque {\bf {101-102}} (1982)
\pages 23--74
\endref


\ref
\key Br
\by  R. K. Brylinski
\paper {\it Limits of weight spaces, Lusztig's $q$-analogs, and fibering of
coajoint orbits}
\paperinfo Journ.~A.M.S., {\bf 2} (1989) 
\pages 517--534
\endref

\ref
\key Ca
\by R. Carter
\paper {\it Finite groups of Lie type}
\paperinfo Wiley, 1985
\endref

\ref
\key CG
\by N. Chriss, V. Ginzburg
\paper {\it Representation theory and complex geometry}
\paperinfo Birkh\"auser, Boston 1997
\endref

\ref
\key CM
\by D.H.~Collingwood, W.M.~McGovern
\paper {\it Nilpotent orbits in semisimple  Lie algebras}
\paperinfo   New York: Van Nostrand Reinhold, 1993
\endref

\ref
\key Dy
\by E.B.~Dynkin
\paper  {\it Semisimple subalgebras of semisimple Lie algebras}
\paperinfo Matem. Sbornik {\bf 30} (1952)
\pages 349--462 (Russian). English translation:
{\it Amer. Math. Soc. Transl.} II~Ser.,
{\bf 6}~(1957), 245--378
\endref

\ref
\key El1
\by A.G.~Elashvili
\paper {\it The centralizers of nilpotent elements in semisimple Lie algebras}
\paperinfo Trudy Tbiliss. Matem. Inst. {\bf 46} (1975)
\pages 109--132 (Russian)
\endref


\ref
\key El2
\by \bysame
\paper {\it Sheets of the exceptional Lie algebras}
\inbook in ``Issledovaniya po algebre''
\publaddr Tbilisi
\yr 1985
\pages 171--194 (Russian)
\endref


\ref 
\key Gi
\by V. Ginzburg
\paper {\it Perverse sheaves on a loop group and Langlands' duality}
\paperinfo Preprint 1995, alg-geom/9511007
\endref



\ref
\key J1
\by A. Joseph
\paper 
{\it Goldie rank polynomials in a semisimple Lie algebra $\op{III}$}
\paperinfo  J. Algebra
${\bold {73}}$ (1981)
\pages 295--326
\endref


\ref
\key J2
\by A. Joseph
\paper 
{\it Characters of unipotent representations}
\paperinfo  J. Algebra
${\bold {130}}$ (1990)
\pages 273--295
\endref

\ref
\key KL
\by D. Kazhdan, G. Lusztig 
\paper 
{\it  Representations of Coxeter groups and 
Hecke algebras}
\paperinfo Invent. Math. ${\bold {53}}$ (1979)
\pages  165--184
\endref

\ref
\key Ke
\by G.~Kempken
\paper {\it Induced conjugacy classes in classical Lie-algebras}
\paperinfo Abh. Math. Sem. Univ. Hamburg {\bf 53}~(1983)
\pages 53--83
\endref



\ref
\key K1
\by B. Kostant
\paper 
{\it  The principal three-dimensional subgroup
and the Betti numbers of a complex simple Lie group}
\paperinfo Amer.\ J. Math., ${\bold {81}}$ (1959)
\pages  973--1032
\endref



\ref
\key K2
\by B. Kostant
\paper {\it Lie group representations on polynomial rings}
\paperinfo Amer. J. Math. {\bf 85} (1963)
\pages 327--404
\endref

\ref
\key K3
\by B. Kostant
\paper {\it On Whittaker vectors and representation theory}
\paperinfo Invent. Math. {\bf 48} (1978)
\pages 101--184
\endref

\ref 
\key L1
\by  G. Lusztig
\paper {\it Singularities, Character formulas and a $q$-analog of
weight multiplicity}
\paperinfo Asterisque ${\bold {101-102}}$ (1983)
\pages 208--229
\endref

\ref 
\key L2
\by  G. Lusztig
\paper {\it Characters of Reductive groups over finite fields}
\paperinfo Princeton University Press
\endref

\ref 
\key L3
\by  G. Lusztig
\paper {\it A class of irreducible representations of a Weyl group $\op{I,II}$}
\paperinfo Indag. Math. ${\bold {41}}$ (1979), 323--335;
${\bold {44}}$ (1982), 219--226
\endref

\ref 
\key LS
\by  G. Lusztig, N. Spaltenstein
\paper {\it  Induced unipotent classes}
\paperinfo  J. London Math. Soc.  ${\bold {19}}$ (1979), no. 1
\pages 41--52
\endref



\ref
\key M1
\by I.G. Macdonald
\paper {\it Some irreducible representations of Weyl groups}
\paperinfo Bull. London Math. Soc. ${\bold 4}$ (1972)
\pages  148--150
\endref

\ref
\key M2
\by I.G. Macdonald
\paper {\it Symmetric functions and Hall polynomials}
\paperinfo 2nd Edition, Clarendon Press 1995
\endref



\ref
\key R
\by R.W. Richardson
\paper {\it Commuting varieties of semisimple
Lie
algebras and algebraic groups}
\paperinfo Comp. Math. {\bf 38} (1979)
\pages 311--327
\endref



\ref
\key Sh
\by T.~Shoji
\paper {\it On the Springer representations of the Weyl groups
of classical algebraic groups}
\paperinfo Commun. Algebra {\bf 7} (1989)
\pages 1713--1745, 2027--2033
\endref

\ref
\key Spa
\by N.~Spaltenstein
\book  {\it Classes Unipotentes et Sous-groupes de Borel}
\bookinfo Lecture notes in Math. {\bf 946}
\publ Springer
\publaddr Berlin Heidelberg New York
\yr 1982
\endref
\ref
\key GP
\by  Computer program GAP
\paperinfo available from
http:{/}\!{/}www-gap.dcs.st-andrews.ac.uk/\~{}gap/ \
\endref


\endRefs
\enddocument